\documentclass[a4paper, 10pt, oneside, reqno]{amsart}

\usepackage[utf8]{inputenc}
\usepackage{amsfonts,amssymb,amsthm,amsmath}
\usepackage{mathrsfs}
\usepackage{graphicx,graphics}
\usepackage{color}

\usepackage{tikz}
\usepackage{tikz-cd}
\usetikzlibrary{nfold}
\usepackage{comment}
\usepackage{subfiles}
\usepackage{stmaryrd}
\usepackage{hyperref}
\usepackage{todonotes}
\usepackage{cleveref}

\hypersetup{colorlinks = true, linkbordercolor = {white}, linkcolor = {blue}, citecolor = {blue}}

\theoremstyle{plain}
\newtheorem{thm}{Theorem}[subsection]
\newtheorem{lem}[thm]{Lemma}
\newtheorem{construction}[thm]{Construction}
\newtheorem{prop}[thm]{Proposition}
\newtheorem{corol}[thm]{Corollary}

\theoremstyle{definition}
\newtheorem{defi}[thm]{Definition}

\theoremstyle{remark}
\newtheorem{rque}[thm]{Remark}
\newtheorem{rques}[thm]{Remarks}
\newtheorem{example}[thm]{Example}
\newtheorem{nota}[thm]{Notation}

\def \F {\mathbb{F}}
\def \Fp {\F_p}
\def \Fq {\F_q}
\def \Fl {\F_{\ell}}

\def \Fqb {\overline{\F}_q}
\def \Flb {\overline{\F}_{\ell}}

\def \Z {\mathbb{Z}}
\def \Zl {\Z_{\ell}}
\def \Zlb {\overline{\Z}_{\ell}}

\def \Q {\mathbb{Q}}
\def \Ql {\Q_{\ell}}
\def \Qlb {\overline{\Q}_{\ell}}

\def \Shv {\mathrm{Shv}}

\def \Ind {\mathrm{Ind}}

\def \Gm {\mathbb{G}_m}
\def \Ga {\mathbb{G}_a}
\def \Spec {\mathrm{Spec}}

\def \Gbf {\mathbf{G}}
\def \Bbf {\mathbf{B}}
\def \Ubf {\mathbf{U}}
\def \Tbf {\mathbf{T}}
\def \Wbf {\mathbf{W}}

\def \coker {\mathrm{coker}}
\def \Hom {\mathrm{Hom}}

\def \id {\mathrm{id}}

\def \Tr {\mathrm{Tr}}

\def \1 {\mathbb{1}}
\def \ev {\mathrm{ev}}
\def \coev {\mathrm{coev}}

\def \End {\mathrm{End}}
\def \Mod {\mathrm{Mod}}

\def \ind {\mathrm{ind}}

\def \Frob {\mathrm{Frob}}

\def \Gal {\mathrm{Gal}}

\def \Fun {\mathrm{Fun}}
\def \cont {\mathrm{cont}}
\def \Frob {\mathrm{F}}
\def \Ad {\mathrm{Ad}}

\def \cl {\mathrm{cl}}

\def \Br {\mathrm{Br}}

\def \RGamma {\mathrm{R\Gamma}}

\def \DD {\mathrm{D}}

\def \unip {\mathrm{unip}}

\def \Stab {\mathrm{Stab}}

\def \Ext {\mathrm{Ext}}
\def \gr {\mathrm{gr}}
\def \Rep {\mathrm{Rep}}
\def \RHom {\mathrm{RHom}}

\def \Perf {\mathrm{Perf}}

\def \gen {\mathrm{gen}}

\def \rev {\mathrm{rev}}
\def \Irr {\mathrm{Irr}}

\def \Ccal {\mathcal{C}}

\def \GL {\mathrm{GL}}
\def \Ocal {\mathcal{O}}

\def \aff {\mathrm{aff}}
\def \Reg {\mathrm{Reg}}
\def \cusp {\mathrm{cusp}}

\def \tors {\mathrm{tors}}

\def \gr {\mathrm{gr}}

\def \Irr {\mathrm{Irr}}

\def \GL {\mathrm{GL}}

\def \Alg {\mathrm{Alg}}

\def \Hcal {\mathcal{H}}
\def \Ocal {\mathcal{O}}

\def \Bcal {\mathcal{B}}

\def \CS {\mathrm{CS}}

\def \SBim {\mathrm{SBim}}

\def \tr {\mathrm{tr}}
\def \ch {\mathrm{ch}}
\def \red {\mathrm{red}}

\def \DL {\mathrm{DL}}

\def \Ch {\mathrm{Ch}}
\def \ch {\mathrm{ch}}
\def \hc {\mathrm{hc}}

\newcommand{\Gbb}{\mathbb{G}}
\newcommand{\RR}{\mathrm{R}}
\newcommand{\Tbb}{\mathbb{T}}
\newcommand{\spe}{\mathrm{sp}}
\newcommand{\BS}{\mathrm{BS}}
\newcommand{\opp}{\mathrm{op}}

\newcommand{\fm}{\mathrm{fm}}
\newcommand{\bimod}{\mathrm{bimod}}
\newcommand{\C}{\mathbb{C}}
\newcommand{\Betti}{\mathrm{Betti}}
\newcommand{\HHH}{\mathrm{HH}}
\newcommand{\uncompl}{\mathrm{uncompl}}
\newcommand{\Prr}{\mathrm{Pr}}
\newcommand{\E}{\mathbb{E}}
\newcommand{\Cc}{\mathrm{C}}
\newcommand{\Seq}{\mathrm{Seq}}
\newcommand{\ofrak}{\mathfrak{o}}
\newcommand{\Add}{\mathrm{Add}}
\newcommand{\EE}{\mathbb{E}}
\newcommand{\Mbb}{\mathbb{M}}
\newcommand{\Proj}{\mathrm{Proj}}
\newcommand{\dec}{\mathrm{dec}}
\newcommand{\Frac}{\mathrm{Frac}}
\newcommand{\Obb}{\mathbb{O}}
\newcommand{\Acal}{\mathcal{A}}
\newcommand{\cyc}{\mathrm{cyc}}
\newcommand{\Fcal}{\mathcal{F}}
\newcommand{\Kos}{\mathrm{Kos}}
\newcommand{\Qcoh}{\mathrm{Qcoh}}
\newcommand{\act}{\mathrm{act}}
\newcommand{\restr}{\mathrm{res}}
\newcommand{\FT}{\mathrm{FT}}
\newcommand{\Lbb}{\mathbb{L}}
\newcommand{\Lbf}{\mathbf{L}}
\newcommand{\Pbf}{\mathbf{P}}
\newcommand{\Vbf}{\mathbf{V}}
\newcommand{\pfrak}{\mathfrak{p}}
\newcommand{\Dcal}{\mathcal{D}}
\newcommand{\unit}{\mathrm{unit}}
\newcommand{\counit}{\mathrm{counit}}
\newcommand{\cycl}{\mathrm{cycl}}

\title{Generic unipotent representations of finite reductive groups}
\author{Arnaud Eteve}

\begin{document}

\maketitle

\noindent\textbf{Abstract.} Let $\Gbb$ be a generic reductive group in the sense of Broué--Malle--Michel. We construct a category over $\Z[t,t^{-1}]$ which interpolates between the various categories of unipotent representations of $\Gbb(q)$ as we specialize $t$ to $q$. We use this category to show two properties of the unipotent decomposition matrix of $\Gbb(q)$ mod $\ell$ (not dividing $q$): that, uniformly in $(q,\ell)$, there are finitely many such matrices and that for $\ell$ large enough, depending only on $\Gbb$, this decomposition matrix depends only on the root of unity corresponding to $q$ mod $\ell$. We also prove that for $q$ and $q'$ two prime powers that are $\ell$-adically close enough, the categories of unipotent mod $\ell$ representations of $\Gbb(q)$ and $\Gbb(q')$ are equivalent.

\tableofcontents

\section{Introduction}

Let $\Gbf$ be a connected reductive group over a finite field $\Fq$. A surprising feature of the complex representation theory of the finite group $\Gbf(\Fq)$ is that a large part of it, i.e. the set of irreducible unipotent representations, their dimensions, character values, \dots, is essentially independent of $q$ and depends only on the ambient reductive group and in fact only on the abstract root datum defining $\Gbf$ together with its Galois action. These phenomena are, for complex representations, a direct consequence of Lusztig's classification \cite[Theorem 4.23]{LusztigBook} of irreducible representations of $\Gbf(\Fq)$. The core concept governing these ideas is that of a \emph{generic reductive group} $\Gbb$ which was introduced in \cite{BroueMalle}. A generic reductive group is the datum of a pair 
$$\Gbb = (\Gamma_{\Gbb}, \Wbf_{\Gbb}\phi)$$
where $\Gamma_{\Gbb} = (X, X^{\vee}, \Phi, \Phi^{\vee})$ is an abstract root datum with Weyl group $\Wbf_{\Gbb}$ and $\phi \in \mathrm{N}_{\GL(X)}(\Wbf_{\Gbb})$ is a finite-order automorphism. With these data, after fixing a finite field $\Fq$, one can define a unique reductive group $(\Gbf, \Tbf)$ with a maximal torus over $\Fq$, such that $\Gal(\Fqb/\Fq)$ acts by $\phi$ on $X^*(\Tbf) = X$. In particular, given $\Gbb$, one can define for every prime power $q$ a finite group $\Gbb(q) = \Gbf(\Fq)$. 

Fix $\ell \neq \mathrm{char}(\Fq)$. Turning to modular representation theory in non-defining characteristic for the finite group $\Gbf(\Fq)$, the subcategory $\Rep^{\unip}_{\Flb}(\Gbf(\Fq))$ of $\Rep_{\Flb}(\Gbf(\Fq))$ consisting of all mod $\ell$ unipotent representations of $\Gbf(\Fq)$ also turns out to have genericity (i.e. independence of $q$) features. Notably, blocks of $\Rep^{\unip}_{\Flb}(\Gbf(\Fq))$ depend, for $\ell$ large enough depending only on $\Gbf$, only on $d = \mathrm{ord}_{\Fl^{\times}}(q)$, the multiplicative order of $q$ mod $\ell$, by \cite{BroueMalleMichel}. The goal of this paper is to start recasting this theory in a more categorical framework.

\subsection{The main construction}

Let us start by fixing $(\Gbb, \ell, q, \Gbf)$ and $\Lambda \in \{\Fl, \Ql, \Zl\}$ as above. From now on, we assume that $\Gbb$ has connected center and  we exclude the twisted cases $({^2}B_2, {^2}F_4, {^2}G_2)$.

\subsubsection{Categorical trace of Frobenius}

 The main idea of this paper is to construct a category that interpolates between all the categories of representations $\DD(\Rep^{\unip}_{\Lambda}(\Gbf(\Fq)))$, and we will do so via categorical traces; we briefly review the formalism and notation in Appendix \ref{appendix:recollection-categorical-traces}. Attached to the tuple $(\Gbf, \Bbf = \Tbf\Ubf, \Lambda)$ there is a monoidal category $\Hcal^{\fm, \unip}_{\Gbf, \Lambda}$ of free-monodromic unipotent étale sheaves of $\Lambda$-modules on $\Ubf \backslash \Gbf/\Ubf$, initially introduced in \cite{BezYun}. This category comes equipped with a monoidal endofunctor via Frobenius
$$\Frob^* : \Hcal^{\fm, \unip}_{\Gbf, \Lambda} \to \Hcal^{\fm, \unip}_{\Gbf, \Lambda}.$$
The link with the representation theory of $\Gbf(\Fq)$ is via the notion of categorical trace. 

\begin{thm}[\cite{EteveFreeMono}, \cite{ZhuCategoricalTame}, \cite{GaitsgoryRozemblyumVarshavsky}]
There is a natural equivalence 
$$\Tr(\restr_{\Frob^*}, \Prr_{\Hcal^{\fm, \unip}_{\Gbf, \Lambda}}) \cong \DD(\Rep_{\Lambda}^{\unip}(\Gbf(\Fq))).$$
\end{thm} 

\subsubsection{Soergel realization}

Let $\widehat{\Tbf}_{\Lambda}$ be the torus dual to $\Tbf$ defined over $\Lambda$ and define 
$$\RR_{\Lambda} := \Lambda[X^{\vee}]^{\wedge}_{I_{\Lambda}} = \Ocal(\widehat{\Tbf}_{\Lambda, 1}^{\wedge}), \qquad I_{\Lambda} = (e^{\lambda} - 1)_{\lambda \in X^{\vee}} \subset \Lambda[X^{\vee}],$$
the ring of functions on the completion of $\widehat{\Tbf}$ along its unit section. Attached to the ring $\RR_{\Lambda}$, one can define a notion of Soergel bimodules; see \cite{BezrukavnikovRiche}. These form an additive monoidal subcategory
$$\SBim_{\Lambda}(\Gamma_{\Gbb}) \subset {_{\RR_{\Lambda}}}\Mod_{\RR_{\Lambda}}$$ 
of all $\RR_{\Lambda}$-bimodules. One then defines the Hecke category  
$$\Hcal_{\Lambda}(\Gamma_{\Gbb}) := \Ind \Ch^b(\SBim_{\Lambda}(\Gamma_{\Gbb})).$$

\begin{thm}[\cite{BezrukavnikovRiche}]
There is an equivalence of monoidal categories
$$\Hcal^{\fm, \unip}_{\Gbf, \Lambda} \cong \Hcal_{\Lambda}(\Gamma_{\Gbb})$$
compatible with Frobenius. 
\end{thm}

Let us denote by $\sigma : X^{\vee} \to X^{\vee}$ the morphism dual to $\phi$. The Frobenius on the Soergel bimodule side is $\psi_q \circ \sigma : \RR_{\Lambda} \to \RR_{\Lambda}$ where 
$$\psi_q :  \RR_{\Lambda} \to \RR_{\Lambda}, e^{\lambda} \mapsto e^{q\lambda}.$$

\subsubsection{Making the Frobenius generic}

The main idea of this paper is now to replace $q$ by an abstract parameter $t$. Let us denote by 
$$B = \Z\left[t,t^{-1}, \binom{t}{n}, n \geq 1 \right] \subset \Q(t), \qquad \binom{t}{n} = \frac{t(t-1)\dots(t-n+1)}{n!}.$$
Note that the subring $B^+ = \Z[t, \binom{t}{n}]$ is the subring of polynomials in $\Q[t]$ which take integer values on integers. Replacing $\Lambda$ by $B$, all the objects previously defined admit variants over $B$, i.e. we can define $\widehat{\Tbf}_B, I_B, \RR_B, \SBim_B(\Gamma_{\Gbb})$ and $\Hcal_{B}(\Gamma_{\Gbb})$. Over $B$, we can define the following automorphism of $\RR_{B}$, see Lemma \ref{lem:action-of-psi-t}
$$\psi_t : \RR_{B} \to \RR_B, \qquad e^{\lambda} \mapsto e^{t\lambda} := \sum_{n \geq 0}\binom{t}{n}(e^{\lambda} - 1)^n.$$
Furthermore, we can specialize 
$$B \to \Lambda, \qquad t \mapsto q$$ 
and under that specialization $\psi_t$ is mapped onto $\psi_q$.

\subsubsection{The category $\DD(\Rep_B^{\unip}(\Gbb(t)))$}

The definition of $\psi_t$ works for any $B$-algebra $A$. We denote all the corresponding objects defined over $A$ by a subscript $_A$. Let us fix such an algebra. 
\begin{defi}[Definition \ref{defi:generic-representation-category}]
The category of generic unipotent representations of $\Gbb$ is defined to be 
$$\DD(\Rep_A^{\unip}(\Gbb(t))) := \Tr(\restr_{\psi_t \circ \sigma}, \Prr_{\Hcal_A(\Gamma_{\Gbb})}).$$
\end{defi}

This is an $A$-linear presentable category. The Hecke category carries an action of $\RR_A^{\Wbf}$ compatible with its monoidal structure and action of $\psi_t \circ \sigma$, see Lemma \ref{lem:lifting-to-rr-invariants}. As the trace construction is symmetric monoidal, there is a canonical action of 
$$\HHH(\RR_A^{\Wbf}, \psi_t \circ \sigma) := \RR_A^{\Wbf} \otimes_{\RR_A^{\Wbf} \otimes_A \RR_A^{\Wbf}} \RR_{A, \psi_t \circ \sigma}^{\Wbf}$$
on $\DD(\Rep_A^{\unip}(\Gbb(t)))$, see Section \ref{sec:recollection-hh} for a reminder on Hochschild homology. 

\begin{defi}[Definition \ref{defi:generic-reduced-representation-category}]
The category of reduced generic unipotent representations of $\Gbb$ is defined to be
$$\DD^{\red}(\Rep_A^{\unip}(\Gbb(t))) := \DD(\Rep_A^{\unip}(\Gbb(t))) \otimes_{\DD(\HHH(\RR_A^{\Wbf}, \psi_t \circ \sigma))} \DD(A)$$
where the map $\HHH(\RR_A^{\Wbf}, \psi_t \circ \sigma) \to A$ is the map induced by the augmentation $\RR_A^{\Wbf} \to A$. 
\end{defi}

Let us now record the main properties we prove concerning these categories. It turns out that most of these properties are formal consequences of the trace formalism. 

\begin{thm}\label{thm:intro-main-properties-of-generic-representation}
\begin{enumerate}
\item There exist compact generators
\[
\tr_A(\mathbb{B}_A) \in \DD(\Rep_A^{\unip}(\Gbb(t)))
\]
and
\[
\tr_{A,\red}(\mathbb{B}_{A,\red}) \in \DD^{\red}(\Rep_A^{\unip}(\Gbb(t))).
\]
See Section \ref{sec:generator-in-trace} for the notation.
\item Set $E_A = \End(\tr_A(\mathbb{B}_A)), E_{A,\red} = \End(\tr_{A,\red}(\mathbb{B}_{A, \red}))$. Both algebras lie in $\DD^{\leq 0}(A)$ and we have 
$$E_{A,\red} = E_A \otimes_{\HHH(\RR_A^{\Wbf}, \psi_t \circ \sigma)} A,$$
and 
\begin{align*}
\DD(\Rep_A^{\unip}(\Gbb(t))) &\cong \DD(E_{A}^{\opp}) \\
\DD^{\red}(\Rep_A^{\unip}(\Gbb(t))) &\cong \DD(E_{A,\red}^{\opp}).
\end{align*}
\item For any map $A \to A'$ of $B$-algebras, we have 
$$E_{A, \red} \otimes_A A' \cong E_{A',\red}.$$
\item Let $|\Gbb| \in \Z[t,t^{-1}]$ be the polynomial order of $\Gbb$, in the sense of \cite{BroueMalle}. Then, for a $B[|\Gbb|^{-1}]$-algebra $A$, the map 
$$E_{A} \to E_{A,\red}$$ 
is an isomorphism. 
\item The algebra $E_{B, \red}$ descends to a $B_N$-algebra where $B_N = \Z[t,t^{-1}, \binom{t}{n}, n \leq N]$ and $N = \ell(w_0) + 1$.
\item Fix $(q, \ell, \Lambda)$ as before. There is a natural equivalence 
$$\DD(\Rep^{\unip}_{\Lambda}(\Gbb(q))) \cong \DD(\Rep^{\unip}_{\Lambda}(\Gbf(\Fq))).$$ 
\end{enumerate}
\end{thm} 

The assertions of this theorem are proved in Section \ref{sec:generic-representation}; let us indicate where each is proved. 
\begin{enumerate}
\item Point $(1)$ is proved in Lemma \ref{lem:generator-of-trace},
\item points $(2), (3)$ and $(5)$ are proved in Lemma \ref{lem:main-lemma-generator-of-trace}, 
\item point $(4)$ is proved in Corollary \ref{cor:generically-category-is-reduced}, 
\item point $(6)$ is proved in Lemma \ref{lem:etale-specialization-for-representations}. 
\end{enumerate}

We further introduce the following finiteness condition on the ring $A$ in Section \ref{sec:serre-functor} 
\begin{enumerate}
\item[(Fin)] The ring $A$ satisfies the property (Fin) if for all $w \in \Wbf$, the ring 
$$\HHH(\RR_A, \psi_t \circ w \circ \sigma)$$ 
is concentrated in degree $0$ and finite locally free over $A$ and we have, see Section \ref{ssec:HHcalculation}	
$$\HHH(\RR_A, \psi_t \circ w \circ \sigma) \cong \HHH^{\cont}(\RR_A, \psi_t \circ w \circ \sigma).$$ 
\end{enumerate}
This condition is satisfied in the following two examples:
\begin{enumerate}
\item if $|\Gbb|$ is invertible in $A$, 
\item if $A = \Zl, \Fl, \Ql$ and $B \to \Lambda$ is obtained by $t \mapsto q$. 
\end{enumerate}

\begin{thm}\label{thm:intro-generic-representations-finiteness}
\begin{enumerate}
\item Assume that $A$ satisfies (Fin), then the algebra $E_A$ is a finite $A$-projective module in degree $0$. 
\item Assume that $A$ satisfies (Fin), then the algebra $E_A$ is canonically self-dual. 
\item The algebra $E_{\Q(t)}$ is reduced and semisimple. Furthermore, if $\Gbb$ is quasi-simple, it is split over $\overline{\Q}(t^{\frac{1}{2}})$. 
\end{enumerate}
\end{thm}

The various points are proved as follows: 
\begin{enumerate}
\item The first two points are proved in Lemma \ref{lem:self-duality}.
\item Semisimplicity is proved in Lemma \ref{lem:generic-semisimplicity} and the splitting field is computed in Theorem \ref{thm:splitting-field}. 
\end{enumerate}

\begin{rque}
The main properties of Theorem \ref{thm:intro-generic-representations-finiteness} come from the analysis of the Serre functor on the generic version of the Deligne--Lusztig category $\Ocal$ which we define in terms of parabolic character sheaves as in \cite{EteveTilting} and the relative Serre duality on Hecke categories due to \cite{HoLiRelativeSerre} and \cite{SandvikRelativeSerreDuality} and parabolic character sheaves due to \cite{GaitsgoryRozemblyumVarshavsky} which we adapt to our setting in Theorem \ref{thm:relative-Serre-duality}. 
\end{rque} 

\begin{rque}
The way we compute some information about the splitting field of $E_{\Q(t)}$ is through Harish--Chandra theory. In Section \ref{sec:generic-parabolic-char-sheaves}, we define a generic version of Deligne--Lusztig induction (and thus a generic version of parabolic induction) and we prove a generic variant of the Mackey formula in Section \ref{ssec:mackay-formula}. The semisimplicity of the algebra $E_{\Q(t)}$ comes from Tits' deformation theorem as we can realize this algebra as a deformation of a corresponding algebra in the category of representations of $\Gbf(\Fq)$ for some $\Fq$ (in fact any $q$ works). Once we have both semisimplicity and a Mackey formula, we are reduced to controlling the splitting field of the cuspidal representations which we can isolate in the generic version of the cohomology of certain Deligne--Lusztig varieties using results of Lusztig \cite{LusztigCoxeterOrbit, LusztigRationality}. 
\end{rque} 

\subsubsection{Parabolic character sheaves} 

In Section \ref{sec:generic-parabolic-char-sheaves} we define generic variants of parabolic character sheaves of Lusztig \cite{LusztigParabCS} and Deligne--Lusztig induction. Let $(q,\ell)$ be as above. We use the fact that we have this category $\DD(\Rep_{B}^{\unip}(\Gbb(t)))$ to deduce certain properties of the category of representations of $\Gbf(\Fq)$. For any finite subset $I \subset S$, let us denote by $\Gamma_I \subset \Gamma_{\Gbb}$ the subroot datum attached to $I$ and by $\Wbf_I \subset \Wbf_{\Gbb}$ the corresponding parabolic subgroup of the Weyl group. There is a natural inclusion of Hecke categories 
$$\Hcal_A(\Gamma_I) \xrightarrow{\iota_I} \Hcal_A(\Gamma_{\Gbb}).$$

\begin{defi}[Definition \ref{defi:generic-parabolic-character-sheaves}]
The category of generic parabolic character sheaves of type $I$ is defined by 
$$\DD^{\unip}(\CS_I) := \Tr([\Hcal_A(\Gamma_{\Gbb})_{\psi_t \circ \sigma}], \Prr_{\Hcal_A(\Gamma_I)}).$$
There are well-defined character and horocycle functors 
$$\ch_I : \DD^{\unip}(\CS_I) \to \DD(\Rep_A^{\unip}(\Gbb(t))) : \hc_I^!, \hc_I^*$$
where we have adjunctions $\hc_I^! \dashv \ch_I \dashv \hc_I^*$. 
\end{defi} 

\begin{lem}[Lemma \ref{lem:etale-realization-gen-parab-cs}]
Let $(q, \ell, \Lambda)$ be as before; then there is an equivalence between $\DD^{\unip}(\CS_I)$ and Lusztig's category of parabolic character sheaves \cite{LusztigParabCS} together with the induction and restriction functors. This also recovers the usual Deligne--Lusztig induction and restriction functors.
\end{lem}

\subsection{Applications to modular representation theory}

\subsubsection{Local constancy over $\Zl^{\times}$}

The first application concerns the dependence of the whole category $\Rep_{\Fl}^{\unip}(\Gbb(q))$ on $q$. Our main theorem concerning this is the following. 

\begin{thm}[Theorem \ref{thm:modular-local-constancy-of-category}]\label{thm:local-constancy-intro}
There exists $n_q \geq 1$ such that for all $(q',\ell)$ with $v_{\ell}(q-q') \geq n_q$, there exists an equivalence of categories 
$$\Rep_{\Fl}^{\unip}(\Gbb(q)) \cong \Rep_{\Fl}^{\unip}(\Gbb(q'))$$
preserving Deligne--Lusztig cohomology complexes.
\end{thm}

The main idea is to introduce the algebra $A = \Ccal(\Zl^{\times}, \Fl)$ of continuous functions $\Zl^{\times} \to \Fl$ and view it as a $B$-algebra via 
$$B \to \Ccal(\Zl^{\times}, \Fl), t \mapsto \iota$$
where $\iota : \Zl^{\times} \to \Fl^{\times} \to \Fl$ is the natural quotient map. Over this ring, we can analyze the algebra $E_A$ to deduce Theorem \ref{thm:local-constancy-intro}.

\subsubsection{On the decomposition matrix}

The second set of applications concerns the unipotent decomposition matrix of $\Gbf(\Fq)$. We denote this matrix by $D_{q,\ell}^{\unip}$. Let us introduce the algebra 
$$\overline{E}_{\red} := H^0(E_{B,\red})/B-\mathrm{tors},$$
as the $B$-torsion-free quotient of $H^0(E_{B,\red})$. As $B$ is a Prüfer domain, this algebra is a direct summand of $H^0(E_{B,\red})$ and is in particular projective over $B$, see Lemma \ref{lem:structural-properties-reduced-ring} and as such there exists some $N' \geq N$ such that
$$\overline{E}_{N', \red} = H^0(E_{\red, B_{N'}})/B_{N'}-\mathrm{tors}$$
descends $\overline{E}_{\red}$ to a projective $B_{N'}$-module. We can furthermore relate the decomposition matrix of this algebra along suitable specialization maps of $\Spec(B_{N'})$ to the matrices $D_{q,\ell}^{\unip}$ as $(q,\ell)$ vary, see Lemma \ref{lem:comparison-decomposition-matrices}. 

As a first application of the existence of this flat algebra, we apply a theorem of Thiel \cite{ThielDecompMatrices} and obtain:

\begin{thm}[Theorem \ref{thm:finiteness-decomposition-matrices}]
There exist finitely many possible decomposition matrices $D_{q,\ell}^{\unip}$ uniformly in $(q,\ell)$.
\end{thm}

In a second application, following the strategy to control the decomposition matrix $D_{q,\ell}^{\unip}$ in type $A$ in terms of Schur algebras, see \cite{GeckHeckeAlgRootOfUnity}, we obtain:

\begin{thm}[Theorem \ref{thm:structure-decomposition-matrix}]
Let $\ell \geq N'$, and assume that we have fixed $(q,\ell)$ as above. Let $\xi$ be the root of unity corresponding to $q$ mod $\ell$. Then there exist two matrices $D_{\xi, \ell}$ and $D_{\xi}$ depending only on $(\xi, \ell)$ and $\xi$, respectively, such that 
$$D^{\unip}_{q,\ell} = D_{\xi,\ell}D_{\xi}.$$
\end{thm}

\subsection{Relation to other work and future directions}

\subsubsection{The work of \cite{BZCHN}}
One of the main sources of inspiration for this paper comes from \cite{BZCHN} and the theory of Block--Getzler complexes. In their situation, they consider the affine Hecke category $\Hcal^{\aff}_k$ over some field $k$ of characteristic $0$. In that situation, this category comes equipped with an action of $\Gm$ which comes from forgetting the grading on usual Soergel bimodules, or geometrically from forgetting some mixed structure on étale sheaves. Now the theory of Block--Getzler complexes, see Section 2.1.2 of \emph{loc. cit.}, is used to compute the enhanced trace 
$$\underline{\Tr}^{k-\mathrm{lin}}(\id, \Hcal^{\aff, \Gm}_k) \in \Qcoh(\Gm/\Gm).$$ 
It turns out that the complex that computes this trace is given by 
$$\Tr^{\Qcoh(\Gm)-\mathrm{lin}}(\act_t, \Hcal^{\aff} \otimes \Qcoh(\Gm))$$
where $\act_t : \Hcal^{\aff}_k \otimes \Qcoh(\Gm) \to \Hcal^{\aff}_k \otimes \Qcoh(\Gm)$ is induced by the universal action of a coordinate $t \in \Gm$. 

When passing from $k$ to $\Zl$, it is an easy exercise to observe that the group $\Zl^{\times}$ acts on the ring $\RR_{\Zl}$ with $\gamma \in \Zl^{\times}$ acting by $\psi_{\gamma}$. Unraveling the construction analogous to the one in \emph{loc. cit.} with $\Zl^{\times}$ in place of $\Gm$ leads to the category $\DD(\Rep_{\Ccal(\Zl^{\times}, \Zl)}^{\unip}(\Gbb(t)))$. 

\subsubsection{Schur algebras}

A very fruitful technique in modular representation theory of finite groups of Lie type involves constructing an algebra over $\Z[t,t^{-1}]$ which interpolates, under specializations $t \mapsto q$, certain parts of $\Rep_{\Zl}\Gbb(q)$. Two of the main examples of such constructions are, of course, Hecke algebras and the Schur algebras of Dipper--James \cite{DipperJamesQSchur} which provide in type $A_n$ a lot of control over $D_{q,\ell}^{\unip}$; we refer to \cite[Section 5.3.4]{RouquierModularLie} for a discussion. We view the construction of the algebras $\overline{E}_{\red}$ and $\overline{E}_{N',\red}$ as variants for abstract reductive groups (up to Morita equivalence) of this kind of construction.

\subsubsection{Local invariance} 

One of the main motivations for this work was Theorem \ref{thm:local-constancy-intro} which started as a question from R. Rouquier to the author; we refer to \cite[Section 6.1]{RouquierModularLie} for a discussion. We compare this result with the following fact. The principal block of $\Rep_{\Flb}\Gbb(q)$ is controlled by the dg-algebra $\RGamma(\Gbb(q), \Flb)$ which by results of Friedlander \cite[Theorem 12.2]{FriedlanderEtaleHomotopy} depends only on the homotopy type of $B\Gbb(\C)$ and $q$; as a consequence, the principal block depends itself only on $\Gbb$ and $\overline{\langle q \rangle} \subset \Zl^{\times}$, see \cite[Section 6.1.2]{RouquierModularLie}, \cite[Section 3]{BrotoMollerOliver}.

\subsubsection{Blocks}

One of the main results of \cite{BroueMalleMichel} is the generic theory of blocks which we have not addressed in this paper. The main input in their paper is through $d$-Harish--Chandra theory and while we have worked out the classical version of Harish--Chandra theory here, the $d$-version is still open using the language of categorical traces.

\subsubsection{Lusztig's classification}

Choosing any specialization $t \mapsto q$ yields a canonical bijection 
$$\Irr_{\overline{\Q(t)}}^{\unip}\Gbb(t) \cong \Irr_{\overline{\Q}}^{\unip}(\Gbb(q)).$$
This yields a proof of Lusztig's theorem \cite{LusztigBook} that the set of irreducible unipotent representations of $\Gbb(q)$ depends only on $(\Wbf, \sigma)$. Here our construction also depends on the root datum of $\Gbb$; however, this dependency is not hard to remove over $\Q$. Recently, Reeves \cite{ReevesUnipotentTrace} has also given a proof of Lusztig's classification in terms of categorical traces and shows that the category $\Rep^{\unip}(\Gbb(q))$ is in fact independent of $q$ and depends only on the underlying Weyl group. The main argument of \emph{loc. cit.} consists in showing that 
$$\Tr(\Frob^*, \Prr_{\Shv(\Bbf\backslash \Gbf/\Bbf, \Qlb)}) \cong \Tr(\id, \Prr_{\mathcal{J}})$$
where $\mathcal{J}$ is the asymptotic Hecke category introduced by Lusztig in \cite{LusztigCellsAndTensorCat}. Combining with the equivalence 
$$\Tr(\Frob^*, \Prr_{\Shv(\Bbf\backslash \Gbf/\Bbf, \Qlb)}) \cong \DD(\Rep_{\Qlb}^{\unip}\Gbf^{\Frob})$$ 
shows that there is an equivalence 
$$\Tr(\id, \Prr_{\mathcal{J}}) \cong \DD(\Rep_{\Qlb}^{\unip}\Gbf^{\Frob})$$
which exhibits the right hand side as being independent of $q$. This last equivalence was already constructed by Lusztig in \cite{LusztigCategoricalCentre} up to switching between categorical trace and center which are equivalent here as the category $\mathcal{J}$ is equipped with a pivotal structure, see \cite[Theorem 5.2]{EliasWilliamson}.

\subsection{Notation}

\subsubsection{Generic reductive groups}\label{ssec:generic-reductive-groups}
We fix $\Gbb = (\Gamma_{\Gbb}, \Wbf_{\Gbb}\phi)$ a generic reductive group in the sense of Broué--Malle--Michel \cite{BroueMalleMichel}. This is the datum of a root datum $\Gamma_{\Gbb} = (X, X^{\vee}, \Phi, \Phi^{\vee})$ and a finite-order automorphism $\phi \in \GL(X)$ which stabilizes $\Phi$. We furthermore assume that $\Gbb$ has connected center, i.e. $X/\Z\Phi$ is torsion-free. We shall write $\sigma : X^{\vee} \to X^{\vee}$ for the morphism dual to $\phi$.

For now, we exclude the twisted cases $({^2}B_2, {^2}F_4, {^2}G_2)$. It is likely that the theory can be extended to these cases as well. 

We choose a base $S$ of $\Phi$ and a representative $\phi$ of the coset
$\Wbf_{\Gbb}\phi$ preserving this base. We write $\Wbf = \Wbf_{\Gbb}$, denote its longest element by $w_0$, and write $\ell(w)$ for the length of $w \in \Wbf$. We also use $S$ for the corresponding set of simple reflections in $\Wbf$.

\subsubsection{\'{E}tale realizations}\label{ssec:etale-realization}

An étale realization datum is a triple $(q, \ell, \Lambda)$, where $q = p^r$ with $r \geq 1$, $p$ and $\ell$ are distinct primes and $\Lambda \in \{\Fl, \Zl, \Ql\}$. Let $(\Gbf, \Bbf, \Ubf)$ be the pinned reductive group with Borel pair and based root datum given by $\Gamma_{\Gbb}$ defined over $\Fqb$. The group $\Gbf$ comes equipped with a Frobenius endomorphism 
$$\Frob : \Gbf \to \Gbf$$ 
which preserves $(\Tbf, \Bbf)$ and acts by $q\phi$ on $X^*(\Tbf)$. 

We fix a choice of a topological generator $\gamma \in \pi_1^t(\Gm) \cong \widehat{\Z}^{(p)}(1)$. This choice yields an isomorphism 
\begin{equation}\label{eq:trivialization-pi-1}
\Lambda\llbracket \pi_1(\Tbf)_{\ell} \rrbracket \cong \RR_{\Lambda}
\end{equation}

\subsubsection{Betti realization}\label{ssec:betti-realization}

Attached to the generic reductive group $\Gbb$, there is a unique pinned connected reductive group $(\Gbf, \Bbf, \Tbf)$ defined over $\C$.

\subsubsection{Categories}

We shall work with $\infty$-categories in the sense of Lurie \cite{HTT,HA}. This formalism is required for the language of categorical traces used in this paper. We shall denote by $\Prr$ the category of presentable stable $\Z$-linear categories with functors given by \emph{left} adjoints. All the presentable categories appearing in this paper will belong to $\Prr$. This category is symmetric monoidal and equipped with Lurie's tensor product. For some category $\Ccal \in \Prr$, we shall denote by $\Ccal^{\omega}$ the subcategory of compact objects. Given a monoidal category $\Ccal$, we shall write $\Prr_{\Ccal}$ for the category of presentable $\Ccal$-categories, i.e.
$$\Prr_{\Ccal} = \Mod_{\Ccal}(\Prr).$$
If $\Ccal$ is symmetric monoidal, then $\Prr_{\Ccal}$ acquires a canonical relative tensor product $\otimes_{\Ccal}$. 

For a commutative ring $A$, we shall denote by $\Prr_A$ the category of presentable, stable, $A$-linear categories, explicitly 
$$\Prr_A := \Prr_{\DD(A)}$$
where $\DD(A)$ denotes the full unbounded derived category of $A$-modules equipped with its monoidal structure given by $\otimes_A$.

For an additive category $\mathcal{A}$, we denote by $\Ch^b(\mathcal{A})$ the category of bounded chain complexes on $\mathcal{A}$ which we view as a small stable exact $\infty$-category. Its homotopy category is $K^b(\mathcal{A})$, the usual homotopy category of $\mathcal{A}$. 

\subsection{Acknowledgements} The author thanks Raphaël Rouquier, Cédric Bonnafé, Olivier Dudas, Jean Michel, Gerhard Hiss, Jens Eberhardt, Sean Cotner, Thibaud van den Hove, Daniel Juteau, Wyatt Reeves, Dennis Gaitsgory and Juan Esteban Rodriguez Camargo for discussions during the preparation of this paper. The author thanks the Max Planck Institute for Mathematics for its support during the preparation of the paper. The author was partially supported by the European Research Council (ERC) under the European Union's Horizon 2020 research and innovation programme (grant agreement n° 950326).

\subsection{AI Use} We used ChatGPT 6 Astra during the writing of this paper to search the literature, help with some technical arguments and for the final spelling and proofreading checks. We discovered the reference \cite{CahenChabert} thanks to AI and the proofs of Lemmas \ref{lem:comparison-of-topologies}, \ref{lem:generic-splitting} and \ref{lem:structural-properties-reduced-ring} were obtained by AI. In particular, the precise bound $N = \ell(w_0) + 1$ in Lemma \ref{lem:existence-descent-datum} was worked out by AI.

\section{Hecke categories}

\subsection{Rings}\label{subsec:definition-polynomial-ring}
We introduce the following notation:
\begin{enumerate}
\item We denote by $\binom{t}{i} \in \Q[t]$ the polynomial 
$$\binom{t}{i} := \frac{t(t-1)\dots (t-i+1)}{i!},$$
which is a polynomial such that $\binom{t}{i}(n) \in \Z$ for all $n \in \Z$. 
\item We denote by 
$$B^+ := \Z\left[t, \binom{t}{i}\right] \subset \Q[t] \qquad \text{and by} \qquad B := B^+\left[\frac{1}{t}\right] \subset \Q[t,t^{-1}].$$
\end{enumerate}
We note that $\Frac(B^+) = \Frac(B) = \Q(t)$, i.e. the map $\Spec(B) \to \Spec(\Z[t,t^{-1}])$ is birational. 

For $k \geq 1$, we denote  by
$$B^+_k = \Z\left[t,\binom{t}{i}, i = 1, \dots, k\right], \qquad B_k = B^+_k\left[\frac{1}{t}\right]. $$
We have $B=\varinjlim_k B_k$. Each $B_k$ is a finitely presented
$\Z[t,t^{-1}]$-algebra, and
\[
 B_k[1/k!]=\Z[1/k!,t,t^{-1}].
\]
\begin{rque}
Note that the map 
$$\Spec(B_k) \to \Spec(\Z[t,t^{-1}])$$
is neither finite nor surjective for $k\geq 2$. Indeed, as we have $2\binom{t}{2} = t(t-1)$, the element $\binom{t}{2}$ cannot be integral over the integrally closed domain $\Z[t,t^{-1}]$.
\end{rque} 

\begin{rque}
The ring $B^+$ is the full subring of $\Q[t]$ consisting of polynomials $P \in \Q[t]$ such that $P(n) \in \Z$ for all $n \in \Z$. We expect that the theory constructed in this paper will admit a generalization to the ``very twisted cases'' which give rise to the Suzuki and Ree groups ${^2}B_2, {^2}F_4, {^2}G_2$ after replacing the pair $(\Z,B)$ by a pair $(\Ocal_K,B_K)$, where $\Ocal_K$ is the ring of integers of a suitable number field $K$ and $B_K$ is the ring of polynomials in $K[t]$ such that $P(x) \in \Ocal_K$ for all $x \in \Ocal_K$. Such rings of polynomials have been introduced and studied in the literature previously; we refer to \cite{CahenChabert} for a complete exposition. 
\end{rque}

From now on, unless specified otherwise, $A$ denotes a commutative ring equipped with a map $B^+\to A$.

\subsubsection{Some specializations}\label{subsubsec:specializations}
There are several morphisms out of $B^+$ and $B$ that we will use: 
\begin{enumerate}
\item Given some $q \in \Z\setminus \{0\}$, evaluation at $q$ yields a surjective map  
$$\spe_{q} : B \to \Z\left[\frac{1}{q}\right] \qquad t \mapsto q, \binom{t}{n} \mapsto \binom{q}{n}.$$ 
\item Let $\gamma \in \Zl$. There is a well-defined map 
$$\spe_{\ell,\gamma} : B^+ \to \Zl$$ 
induced by $t \mapsto \gamma$. Note that this is well-defined as for all elements $\gamma \in \Zl$, the $\ell$-adic valuation $v_{\ell}(\binom{\gamma}{i})$ is nonnegative. Indeed, this is true if $\gamma \in \Z$ and extends by density to all of $\Zl$. 
Note that if $\gamma$ is invertible, then $\spe_{\ell, \gamma}$ induces a map 
$$\spe_{\ell, \gamma} : B \to \Zl.$$
For $\gamma\in\Zl^{\times}$, we shall denote by 
\begin{equation}\label{eq:definition-point-bql}
b_{\gamma, \ell} \in \Spec(B)
\end{equation}
the point corresponding to the map $B \xrightarrow{\spe_{\ell, \gamma}} \Zl \to \Fl$.
\item Let $\Ccal(\Zl, \Zl)$ be the ring of continuous functions $\Zl \to \Zl$. There is a well-defined map 
$$\spe_{\Zl} : B^+ \to \Ccal(\Zl, \Zl), \qquad t \mapsto \id_{\Zl}.$$
As for the previous point, we can restrict to $\Zl^{\times}$ and we obtain a map 
$$\spe_{\Zl^{\times}} : B \to \Ccal(\Zl^{\times}, \Zl), t \mapsto \iota$$ 
where $\iota : \Zl^{\times} \to \Zl$ is the canonical inclusion. Note that after choosing $\gamma \in \Zl$, we have a natural map $\ev_{\gamma} : \Ccal(\Zl, \Zl) \to \Zl$ given by evaluation at $\gamma$ and there is an equality $\ev_{\gamma} \circ \spe_{\Zl} = \spe_{\ell, \gamma}$. 

\end{enumerate}

\subsection{The dual torus}
For an affine base $S = \Spec(A)$, we denote by $\widehat{\Tbf}_A$ the torus whose character lattice is $X^{\vee}$ over $S$. We think of it as the torus dual to the generic torus $\Tbb$. We denote by $\widehat{\Tbf}^{\wedge}_1$ its formal completion along its unit section and we set 
$$\RR_A = \Ocal(\widehat{\Tbf}^{\wedge}_1).$$ 
For any $\lambda \in X^{\vee}$, there is a canonical element $e^{\lambda} \in \Ocal(\widehat{\Tbf})$. We shall also denote its image in $\RR_A$ by $e^{\lambda}$. With this notation, $\RR_A$ is the completion of $\Ocal(\widehat{\Tbf})$ along the augmentation ideal denoted by $I_A$. We have 
$$I_A = \left\langle (e^{\lambda} - 1), \lambda \in X^{\vee} \right\rangle.$$
We shall also write $I_A \subset \RR_A$ for the corresponding ideal so that we have 
$$\RR_A/I_A = A.$$ 

\begin{rque}
The ring $\RR_A$ is noncanonically isomorphic to 
$$\RR_A \cong A\llbracket t_1, \dots, t_n \rrbracket.$$
\end{rque} 

\begin{rque}
The ring $\RR_A$ carries an action of the group $\Wbf_{\Gbb} \rtimes \langle \phi \rangle$.
\end{rque}

\subsection{The canonical action}
We take as a base the ring $B^+$ and we consider the ring $\RR_{B^+}$. We introduce the following notation: 
$$e^{t\lambda} := \sum_{n \geq 0}\binom{t}{n}(e^{\lambda} - 1)^n  \in \RR_{B^+}.$$
Note that this is well-defined as $(e^{\lambda} - 1)^n \in I_{B^+}^n$ for all $n$. 

\begin{lem}\label{lem:action-of-psi-t}
There is a well-defined endomorphism $\psi_t : \RR_{B^+} \to \RR_{B^+}$ which is $B^+$-linear and 
$$\psi_t(e^{\lambda}) = e^{t\lambda}.$$
For any $B^+ \to A$, there is a well-defined endomorphism of $\RR_A$ given by the same formula. Furthermore, the endomorphism $\psi_t : \RR_B \to \RR_B$ is invertible. 
\end{lem}

\begin{proof}
Let us first prove that 
$$\psi_t(e^{\lambda})\psi_t(e^{\mu}) = \psi_t(e^{\lambda + \mu}).$$
We expand both sides in $\Q[t]\llbracket e^{\lambda}-1, e^{\mu}-1\rrbracket$ as 
$$\psi_t(e^{\lambda})\psi_t(e^{\mu}) = \sum_{i,j\geq 0} P_{i,j}(t)(e^{\lambda}-1)^i(e^{\mu}-1)^j$$
and 
$$\psi_t(e^{\lambda + \mu}) = \sum_{i,j\geq 0} Q_{i,j}(t)(e^{\lambda}-1)^i(e^{\mu}-1)^j,$$
where $P_{i,j}(t), Q_{i,j}(t) \in \Q[t]$ are some polynomials. It is therefore enough to prove that for all $i,j \geq 0$, we have $P_{i,j}(t) = Q_{i,j}(t)$. Clearly, for all $n \geq 0$, we have $P_{i,j}(n) = Q_{i,j}(n)$, hence $P_{i,j} = Q_{i,j}$. 

It follows that there is a well-defined $B^+$-linear map 
$$\Ocal(\widehat{\Tbf}_{B^+}) \to \RR_{B^+}, e^{\lambda} \mapsto e^{t\lambda}.$$
As $\psi_t(e^{\lambda}) - 1 = (e^{\lambda} - 1)\left(\sum_{n \geq 1}\binom{t}{n}(e^{\lambda} - 1)^{n-1}\right) \in I_{B^+}$, it follows that the map 
$\Ocal(\widehat{\Tbf}_{B^+}) \to \RR_{B^+}$ extends continuously to $\RR_{B^+}$ and induces a map $\psi_t : \RR_{B^+} \to \RR_{B^+}$. The same argument holds over any commutative ring $A$ over $B^+$. 

It remains to prove that over $B$ the morphism $\psi_t$ is invertible. Let us choose a basis $\lambda_1, \dots, \lambda_n$ of $X^{\vee}$. This yields an isomorphism 
$$\RR_B \cong B\llbracket e^{\lambda_1} - 1 \rrbracket \widehat{\otimes}_B \dots \widehat{\otimes}_B B\llbracket e^{\lambda_n} - 1 \rrbracket,$$
where $\widehat{\otimes}_B$ denotes the completed tensor product. As $\psi_t(e^{\lambda} - 1) \in B\llbracket e^{\lambda} - 1 \rrbracket$, we can reduce to the case where $X^{\vee}$ has rank one. Now we have 
$$\psi_t(e^{\lambda} - 1) = t(e^{\lambda} - 1) + (e^{\lambda} - 1)^2\left(\sum_{n \geq 2}\binom{t}{n}(e^{\lambda} - 1)^{n-2}\right)$$
and since $t$ is invertible in $B$, it follows that the map $\psi_t : B\llbracket e^{\lambda} - 1 \rrbracket \to B\llbracket e^{\lambda} - 1 \rrbracket$ is invertible. 
\end{proof}

\begin{rque}
Note that the endomorphism $\psi_t$ commutes with the action of $\Wbf_{\Gbb} \rtimes \langle \phi \rangle$.
\end{rque} 

\subsubsection{\'{E}tale realization}

Fix $q = p^r$ with $r \geq 1$ and $\ell \neq p$ as in \ref{ssec:etale-realization}. There is a well-defined torus $\Tbf$ over $\Fqb$ equipped with a Frobenius $\Frob : \Tbf \to \Tbf$ whose action on $X^*$ is given by $q\phi$. Consider now the category $\widehat{\DD}^{\fm}(\Tbf, \Lambda)$ of free-monodromic unipotent sheaves of $\Lambda$-modules for $\Lambda \in \{\Zl, \Fl, \Ql\}$ as defined in \cite[Appendix A]{BezYun}, \cite[Sections 2--5]{BezrukavnikovRiche}
and \cite[Sections 2.3--2.7]{EteveFreeMono}.

\begin{lem}[\protect{\cite[Lemma 2.3.1]{EteveFreeMono}}]\label{lem:compat-with-frobenius-on-rt}
There is an equivalence $\widehat{\DD}^{\fm}(\Tbf, \Lambda) \xrightarrow{1^*} \DD(\RR_{\Lambda})$, depending only on the isomorphism \eqref{eq:trivialization-pi-1} 
and a commutative diagram
\[\begin{tikzcd}
	{\widehat{\DD}^{\fm}(\Tbf, \Lambda) } & { \DD(\RR_{\Lambda})} \\
	{\widehat{\DD}^{\fm}(\Tbf, \Lambda) } & { \DD(\RR_{\Lambda})}
	\arrow["{1^*}", from=1-1, to=1-2]
	\arrow["{\Frob^*}"', from=1-1, to=2-1]
	\arrow["{\restr_{\psi_q \circ \sigma}}", from=1-2, to=2-2]
	\arrow["{1^*}"', from=2-1, to=2-2]
\end{tikzcd}\]
\end{lem}

\subsection{The free-monodromic Hecke category}

\begin{defi}
We define the category of (free-monodromic) Soergel bimodules over $B^+$ to be the category 
$$\SBim_{B^+}(\Gamma_{\Gbb}) := \left\langle \RR_{B^+} \otimes_{\RR_{B^+}^s} \RR_{B^+}, s \in S \right\rangle_{\otimes_{\RR_{B^+}}, \oplus, \text{direct summand}}.$$
Correspondingly, we define the Hecke category 
$$\Hcal_{B^+}(\Gamma_{\Gbb}) := \Ind(\Ch^b(\SBim_{B^+}(\Gamma_{\Gbb})))$$
where $\Ch^b$ denotes the small stable $(\infty,1)$-category of bounded chain complexes. 
\end{defi}

\begin{nota}
As usual for Hecke categories, for a word $\underline{w} = (s_1, \dots, s_n)$ with letters in $S$, we shall write 
$$\BS_{B^+}(\underline{w}) := \RR_{B^+} \otimes_{\RR^{s_1}_{B^+}} \dots \otimes_{\RR_{B^+}^{s_n}} \RR_{B^+},$$
for the Bott--Samelson bimodule.  We denote the set of all finite words in $S$ by $\Seq(\Wbf,S)$.
\end{nota} 

\begin{rque}
For a ring $B^+ \to A$, we shall denote by $\SBim_{A}(\Gamma_{\Gbb})$ and $\Hcal_A(\Gamma_{\Gbb})$ the corresponding objects defined with $A$ in place of $B^+$. We shall also drop $A$ and/or $\Gamma_{\Gbb}$ if they are clear from context. 
\end{rque}

\subsubsection{Rouquier complexes}

\begin{defi}[Demazure operators, \cite{Demazure}]
Let $s$ be a simple reflection and $\alpha^{\vee}$ the corresponding simple coroot. 
We define the Demazure operator $D_s : \RR_A \to \RR_A$ as 
$$D_s(e^{\lambda}) = \frac{e^{\lambda} - e^{s(\lambda) - \alpha^{\vee}}}{1 - e^{-\alpha^{\vee}}}.$$
\end{defi}

\begin{lem}\label{lem:frobenius-extension}
\begin{enumerate}
\item The operator $D_s$ is $\RR_A^s$-linear and its image is $\RR_A^s$.
\item The $\RR_A^s$-module $\RR_A$ is free of rank $2$.
\end{enumerate}
\end{lem}

\begin{proof}
Let $\lambda \in X^{\vee}$ be a fundamental coweight, i.e. such that $\langle \lambda, \alpha\rangle = 1$, and fix a basis $z_1, \dots, z_n$ of $(X^{\vee})^s$. We thus have 
$$A[X^{\vee}]^s = A[e^{\lambda} + e^{-\lambda}] \otimes_A A[e^{\pm z_i}].$$
On completion, we have 
$$\RR_A = A\llbracket e^{\lambda} - 1\rrbracket \widehat{\otimes}_A A\llbracket e^{z_i} - 1 \rrbracket.$$
As $s(z_i) = z_i$, we have 
$$\RR_A^s = A\llbracket e^{\lambda} - 1\rrbracket^s \widehat{\otimes}_A A\llbracket e^{z_i} - 1 \rrbracket.$$
Set $x = e^{\lambda} - 1$ and $u = e^{\lambda} + e^{-\lambda} - 2$. We claim that 
$$A\llbracket x\rrbracket^s = A\llbracket u \rrbracket.$$ 
First note that $s(u) = u$, hence the RHS is contained in the LHS. As $u = \frac{x^2}{x-1}$, we have $A\llbracket x \rrbracket = A\llbracket u \rrbracket \oplus xA\llbracket x \rrbracket$. Let $f \in A\llbracket x \rrbracket^s$ and write it as $f(x) = a(u) + xb(u)$. As we have $s(x) = x-u$, we get 
$$a(u) + xb(u) = a(u) + (x-u)b(u)$$
and therefore $ub(u) = 0$. As $u$ is not a zero divisor in $A\llbracket x\rrbracket$, we have $b(u) = 0$ and the claim is proved. 

We have $A[e^{\pm\lambda}]^s = A[e^{\lambda} + e^{-\lambda}]$ and by the above calculation, we have $A\llbracket e^{\lambda} - 1 \rrbracket^s = A[e^{\lambda} + e^{-\lambda}]^{\wedge}_{(e^{\lambda} + e^{-\lambda} - 2)}$, the completion along the ideal $(e^{\lambda} + e^{-\lambda} - 2)$. Now both points of the lemma hold for the extension $A[X^{\vee}]^s \subset A[X^{\vee}]$ by direct calculation, see also \cite[Section 2.5]{EteveEberhardt} and thus hold for $\RR^s_A \subset \RR_A$ after completion. 
\end{proof}

We shall also write $D_s : \BS_A(s) \to \RR_A$ for the map induced by tensoring $D_s$ with $\RR_A$ over $\RR_A^s$.

\begin{thm}[\protect{\cite[Proposition 2.4]{KostantKumar}}]\label{thm:braid-relations-Demazure-operators}
The operators $(D_s)$ satisfy the braid relations.
\end{thm}

\begin{defi}\label{defi:demazure-operator-general}
It follows from Theorem \ref{thm:braid-relations-Demazure-operators} that for all $w \in \Wbf$, there is a well-defined operator 
$$D_w = D_{s_1} \circ \dots \circ D_{s_n} : \RR_A \to \RR_A$$
where $w = s_1\dots s_n$ is any reduced expression. 
\end{defi}

\begin{defi}[Rouquier complexes]
For a simple reflection $s \in S$, we define 
$$\Delta_s := [\BS_A(s) \xrightarrow{\mu_s}\RR_A] \in \Ch^b(\SBim_A(\Gamma_{\Gbb}))$$
to be the corresponding Rouquier complex, with $\BS_A(s)$ in degree $0$, $\RR_A$ in degree $1$ with $\mu_s$ the multiplication map. For a word $\underline{w} = (s_1, \dots, s_n)$, we write 
$$\Delta_{\underline{w}} := \Delta_{s_1} \otimes_{\RR_A} \dots \otimes_{\RR_A} \Delta_{s_n}.$$

Similarly for $s \in S$, we define 
$$\nabla_s := [\RR_A \to \BS_A(s)] \in \Ch^b(\SBim_A(\Gamma_{\Gbb}))$$
to be the dual of $\Delta_s$. For a word $\underline{w} = (s_1, \dots, s_n)$, we write 
$$\nabla_{\underline{w}} := \nabla_{s_1} \otimes_{\RR_A} \dots \otimes_{\RR_A} \nabla_{s_n}.$$
\end{defi}

\begin{rque}
We will prove that Rouquier complexes satisfy braid relations in \ref{sssec:Betti-realization-hecke-cat} by comparing them to the Betti version of the free-monodromic Hecke category. It seems likely that one could reproduce Rouquier's original argument in this context \cite{RouquierCategorification}. 
\end{rque} 

\begin{lem}\label{lem:projectivity-over-rr}
Any object $M \in \SBim_A(\Gamma_{\Gbb})$ is projective as a left and right $\RR_A$-module. 
\end{lem}

\begin{proof}
This is true for $\BS_A(s)$ by Lemma \ref{lem:frobenius-extension}. If $M,N$ are projective as left and right modules, then so is $M \otimes_{\RR_A} N$. It follows that all Bott--Samelson bimodules are projective over $\RR_A$ as left and right modules. Finally, being projective is stable under taking direct summands, hence all objects of $\SBim_A(\Gamma_{\Gbb})$ are projective over $\RR_A$.
\end{proof}

\subsubsection{The universal Frobenius on the Hecke category}

\begin{lem}
For $A$ defined over $B$, the functors $\restr_{\psi_t \otimes \psi_t}, \restr_{\sigma \otimes \sigma} : \RR_A-\bimod \to \RR_A-\bimod$ preserve the subcategories of Soergel bimodules. 
\end{lem}

\begin{proof}
There are natural maps of $\RR_A$-bimodules
\begin{enumerate}
\item $\BS_A(s) \xrightarrow{r \otimes r' \mapsto \psi_t(r) \otimes \psi_t(r')} \restr_{\psi_t \otimes \psi_t}(\BS_A(s))$
\item $\BS_A(s) \xrightarrow{r \otimes r' \mapsto \sigma(r) \otimes \sigma(r')} \restr_{\sigma \otimes \sigma}(\BS_A(\sigma(s))).$
\end{enumerate}
The second one is always an isomorphism; the first one is also an isomorphism because $A$ is defined over $B$ and $\psi_t$ is an isomorphism by Lemma \ref{lem:action-of-psi-t}. 
Furthermore, both functors $\restr_{\psi_t \otimes \psi_t}$ and $\restr_{\sigma \otimes \sigma}$ are compatible with convolution of bimodules, hence both functors preserve the subcategories of Bott--Samelson bimodules and thus of Soergel bimodules. 
\end{proof}

\subsubsection{\'{E}tale realization}

Let $(q, \ell)$ be as in \ref{ssec:etale-realization}. Let $\widehat{\DD}^{\fm}(\Ubf \backslash \Gbf/\Ubf, \Lambda)$ be the corresponding free-monodromic Hecke category as in \cite{BezYun, BezrukavnikovRiche}. 

\begin{thm}[\protect{\cite{BezYun}}, \protect{\cite[Theorem 11.9]	{BezrukavnikovRiche}}]\label{thm:equiv-specialization-with-hecke-category}
After choosing a nontrivial $p$-th root of $1$ in a finite extension $\Lambda'$ of $\Lambda$, there is an equivalence of categories
\begin{equation}\label{eq:specialization-hecke-category-to-geometry}
\widehat{\DD}^{\fm}(\Ubf \backslash \Gbf/\Ubf, \Lambda') \cong \Hcal_{\Lambda'}(\Gamma_{\Gbb}). 
\end{equation}
\end{thm} 

\begin{rque}
The functor that yields the equivalence \eqref{eq:specialization-hecke-category-to-geometry} is constructed using both a Whittaker datum for $\Gbf$ and an Artin--Schreier sheaf on $\Ga$. The Whittaker datum is provided by the pinning of $\Gbf$ and the Artin--Schreier sheaf by the $p$-th root of $1$ in $\Lambda'$. We refer to \emph{loc. cit.} for details. 
\end{rque} 

\begin{lem}\label{lem:specialization-etale-frobenius}
Under the equivalence of Theorem \ref{thm:equiv-specialization-with-hecke-category}, the functor $\restr_{\psi_q \circ \sigma \otimes \psi_q \circ \sigma}$ corresponds to the functor $\Frob^*$. 
\end{lem}

\begin{proof}
This is a combination of Lemma \ref{lem:compat-with-frobenius-on-rt} and \cite[Lemma 3.8.14]{EteveJordan}.
\end{proof}

\subsubsection{Betti realization}\label{sssec:Betti-realization-hecke-cat}

Let $(\Gbf, \Bbf, \Ubf)$ be the pinned reductive group defined over $\C$ and let $\Hcal_{\Betti, A} := \widehat{\DD}^{\fm}(\Ubf \backslash \Gbf / \Ubf, A)$ be the free-monodromic Hecke category as in \cite{EteveEberhardt} and \cite{TaylorUniversalMonodromic}. Let us denote by $\SBim_{A}^{\uncompl}(\Gamma_{\Gbb})$ and $\Hcal_A^{\uncompl}(\Gamma_{\Gbb}) = \Ind\Ch^b(\SBim_{A}^{\uncompl}(\Gamma_{\Gbb}))$ the \emph{uncompleted} categories of Soergel bimodules and Hecke categories as in \emph{loc. cit.}. These are defined in the same way as $\Hcal_{A}(\Gamma_{\Gbb})$ by using $\Ocal(\widehat{\Tbf}_{A})$ in place of $\RR_A$.

\begin{thm}[\cite{EteveEberhardt}, \cite{TaylorUniversalMonodromic}]\label{thm:betti-realization-hecke-cat}
There is an equivalence of categories, depending only on the Whittaker datum of $\Gbf$ coming from its pinning: 
\begin{equation}
\Hcal_{\Z}^{\uncompl}(\Gamma_{\Gbb}) \cong \Hcal_{\Betti, \Z}.
\end{equation}
\end{thm}

\begin{lem}
Let $A$ be a commutative ring. There is a natural equivalence 
$$\Hcal^{\uncompl}_A(\Gamma_{\Gbb}) \cong \Hcal_{\Betti, A}.$$
\end{lem}

\begin{proof}
There are natural equivalences 
$$\Hcal^{\uncompl}_A(\Gamma_{\Gbb}) = \DD(A) \otimes_{\DD(\Z)} \Hcal^{\uncompl}_{\Z}(\Gamma_{\Gbb}) \qquad \Hcal_{\Betti, A} = \DD(A) \otimes_{\DD(\Z)} \Hcal_{\Betti, \Z}.$$
The lemma is then a consequence of Theorem \ref{thm:betti-realization-hecke-cat}.
\end{proof}

\begin{lem}\label{lem:Rouquier-complexes-braid-relations}
The Rouquier complexes in $\Hcal_A(\Gamma_{\Gbb})$ satisfy the braid relations. 
\end{lem} 

\begin{proof}
Given two expressions $\underline{w}, \underline{w}'$ which agree in $\Br_{\Wbf}^+$, there is a canonical isomorphism 
$$\Delta_{\underline{w}}^{\uncompl} \cong \Delta_{\underline{w}'}^{\uncompl}$$ 
in $\Ch^b(\SBim_{A}^{\uncompl})$. Since the functor 
\begin{align*}
\Mod_{\Ocal(\widehat{\Tbf}_A)} &\to \Mod_{\RR_A}\\
M &\mapsto M^{\wedge},
\end{align*}
induced by completion is exact on finite free modules, there is a canonical isomorphism 
$$\Delta_{\underline{w}}^{\uncompl, \wedge} \cong \Delta_{\underline{w}}.$$
This establishes the lemma.
\end{proof}

\begin{defi}
We shall write $\Delta_b$ for the Rouquier complex attached to an element in $\Br_{\Wbf}$. Here $\Br_{\Wbf}$ is the braid group associated with $\Wbf$, and $\Br_{\Wbf}^{+}$ is its positive braid monoid.
\end{defi}

\begin{lem}\label{lem:rank-one-properties}
Let $s$ be a simple reflection.
\begin{enumerate}
\item There are natural triangles in $\Hcal_A(\Gamma_{\Gbb})$: 
$$\Delta_1 \to \BS(s) \to \nabla_s$$
and 
$$\Delta_s \to \BS(s) \to \Delta_1.$$
\item We have 
$$\BS(s) * \BS(s) \cong \BS(s) \oplus \BS(s).$$ 
\item We have $\Delta_s * \Delta_s \in \left\langle \Delta_1, \Delta_s \right\rangle$. 
\end{enumerate}
\end{lem}

\begin{proof}
Item $(3)$ is a formal consequence of the first two. The first two hold in $\Hcal_{\Z}^{\uncompl}(\Gamma_{\Gbb})$ by \cite{EteveEberhardt} and are compatible with base change along $\Z \to A$ and then with completion, which establishes the lemma.
\end{proof}

\subsubsection{Generators}

\begin{defi}
We put
\[
 \mathbb{B}_A := \bigoplus_{w\in \Wbf_{\Gbb}}
       \bigoplus_{\underline w\in\Reg(w)}\BS_A(\underline w),
\]
where $\Reg(w)$ is the set of reduced expressions of $w$. For $w = 1$, we set $\BS_A(\emptyset) = \RR_A$.
\end{defi}

\begin{lem}\label{lem:compact-generator-of-hecke-cat}
The object $\mathbb{B}_A$ is a compact generator of $\Hcal_{A}$. 
\end{lem}

\begin{proof}
We prove that the collection of objects $(\Delta_w)_{w \in \Wbf_{\Gbb}}$ generates $\Hcal_{A}(\Gamma_{\Gbb})$. The subcategory it generates contains all Bott--Samelson bimodules of the form $\BS_A(s)$ for $s \in S$, hence it is enough to prove that 
$$\left \langle \Delta_w, w \in \Wbf_{\Gbb} \right\rangle$$
is closed under convolution. Let $w \in  \Wbf_{\Gbb}$ and $s \in S$. It is enough to prove that $\Delta_w * \Delta_s \in \left \langle \Delta_w, w \in \Wbf_{\Gbb} \right\rangle$. There are two cases 
\begin{enumerate}
\item if $\ell(w) + \ell(s) = \ell(ws)$, then $\Delta_w * \Delta_s = \Delta_{ws}$, 
\item if $\ell(ws) + \ell(s) = \ell(w)$, then $\Delta_w * \Delta_s = \Delta_{ws} * \Delta_s * \Delta_s$. Furthermore, by Lemma \ref{lem:rank-one-properties}, we have $\Delta_s * \Delta_s \in \left\langle \Delta_1, \Delta_s \right\rangle$, hence $\Delta_w * \Delta_s \in \langle \Delta_{w}, \Delta_{ws} \rangle$.
\end{enumerate}
By construction, each $\Delta_w$ is isomorphic to a complex of Bott--Samelson bimodules, hence 
$$\left\langle \Delta_w, w \in \Wbf_{\Gbb} \right\rangle \subset \left\langle \mathbb{B}_A \right\rangle$$
and the lemma is proved. 
\end{proof}

\begin{lem}\label{lem:additive-generator}
There is an equivalence 
$$\SBim_A(\Gamma_{\Gbb}) \cong \Add(\mathbb{B}_A)$$ 
where $\Add(\mathbb{B}_A)$ denotes the additive subcategory of $\Hcal_{A}$ generated by finite direct sums and direct summands of copies of $\mathbb{B}_A$. 
\end{lem}

\begin{proof}
We use the theory of weight structures in the sense of Bondarko \cite{Bondarko}. By \cite[Theorem 4.3.2]{Bondarko}, as the category of Soergel bimodules is a full generating subcategory of $\Hcal_A^{\omega}$ composed of objects which are negative (i.e. $\Ext^i(M,N) = 0$ for $i > 0$ and $M,N \in \SBim_A(\Gamma_{\Gbb})$), there exists a unique weight structure whose heart is $\SBim_{A}(\Gamma_{\Gbb})$. As $\mathbb{B}_A$ is a generator of $\Hcal_A^{\omega}$, by the same theorem, there exists a unique weight structure on $\Hcal_A^{\omega}$ whose heart is the Karoubi completion of the category generated by finite direct sums of copies of $\mathbb{B}_A$. By uniqueness, these two weight structures coincide and we have 
$$\SBim_A(\Gamma_{\Gbb}) \cong \Add(\mathbb{B}_A).$$
\end{proof}

\begin{lem}\label{lem:hecke-category-is-rigid}
The Hecke category $\Hcal_{A}(\Gamma_{\Gbb})$ is rigid. 
\end{lem}

\begin{proof}
It follows from the proof of Lemma \ref{lem:compact-generator-of-hecke-cat} that the Hecke category is compactly generated by the collection of $(\Delta_{w,A})_{w \in \Wbf}$. These objects are invertible by Lemma \ref{lem:Rouquier-complexes-braid-relations} and hence dualizable. The lemma follows.
\end{proof}

\subsection{Semi-orthogonal decomposition}

\begin{lem}\label{lem:orthogonality}
Let $w,v \in \Wbf$. 
\begin{enumerate}
\item if $w \neq v$, we have 
$$\Hom_{\Hcal_{A}(\Gamma_{\Gbb})}(\Delta_{A,w}, \nabla_{A,v}) = 0,$$
\item if $w = v$, then we have 
$$\Hom_{\Hcal_{A}(\Gamma_{\Gbb})}(\Delta_{A,w}, \nabla_{A,v}) \cong \RR_{A}.$$
\end{enumerate}
\end{lem}

\begin{proof}
Since the objects $b \mapsto \Delta_{b}$ for $b \in \Br_{\Wbf}$ satisfy the braid relations, it is enough to prove that 
$$\Hom(\Delta_1, \nabla_w) = 0$$
if $w \neq 1$ and 
$$\Hom(\Delta_1, \nabla_1) = \RR_A.$$
The second case is immediate as $\Delta_1 = \nabla_1 = \RR_A$ as an $\RR_{A}$-bimodule. Let us now prove the first case. Using the notation of Section \ref{sssec:Betti-realization-hecke-cat}, we have 
$$\Hom(\Delta_{1,A}^{\uncompl}, \nabla_{w,A}^{\uncompl}) = 0$$
if $w \neq 1$. Now this object in $\DD(A)$ is represented by the complex 
$$\Hom^{\bullet}_{A[X^{\vee}] \otimes A[X^{\vee}]}(A[X^{\vee}], \nabla_{w,A}^{\uncompl})$$ 
which is a finite bounded complex of finite free $A[X^{\vee}]$-modules, by inspection. The functor $M \mapsto \RR_A \otimes_{A[X^{\vee}]} M$ sends this complex onto the same complex computing $\Hom^{\bullet}(\Delta_1, \nabla_w)$. Moreover, this functor is exact on finite free $A[X^{\vee}]$-modules. The lemma follows.  
\end{proof}

\begin{defi}
Let $Z \subset \Wbf$. We say that $Z$ is closed if for all $z \in Z$ and $w \in \Wbf$ such that $w \leq z$ in the Bruhat ordering, we have $w \in Z$. We say that $U \subset \Wbf$ is open if its complement is closed. For a closed subset $Z\subset \Wbf$, we denote by 
$$\Hcal_{A,Z}(\Gamma_{\Gbb}) = \left\langle \Delta_w, w \in Z \right\rangle \subset \Hcal_{A}(\Gamma_{\Gbb})$$
and for an open subset $U$ with complement $Z$, we denote by $\Hcal_{A,U}(\Gamma_{\Gbb}) := \Hcal_{A}(\Gamma_{\Gbb})/\Hcal_{A,Z}(\Gamma_{\Gbb})$. We also denote by 
$$\iota_Z : \Hcal_{A,Z}(\Gamma_{\Gbb}) \to \Hcal_A(\Gamma_{\Gbb}), \qquad \iota_U^* : \Hcal_A(\Gamma_{\Gbb}) \to \Hcal_{A,U}(\Gamma_{\Gbb})$$ 
the natural inclusion and quotient map. 
\end{defi} 

\begin{thm}[Semi-orthogonal decomposition]\label{thm:semi-orthogonal-decomposition-hecke-cat}
\begin{enumerate}
\item Each of the families $(\Delta_w)_{w\in Z}$, $(\nabla_w)_{w\in Z}$, and $(\BS_A(\underline{w}))_{\underline{w}\in\Reg(w),\,w\in Z}$ generates $\Hcal_{A,Z}(\Gamma_{\Gbb})$.
\item Both functors $\iota_Z$ and $\iota_U^*$ have continuous left and right adjoints that preserve compact objects sitting in adjunction triples 
$$\iota_Z^* \dashv \iota_Z \dashv \iota_Z^!, \qquad \iota_{U,!} \dashv \iota_U^* \dashv \iota_{U,*}.$$
\item The essential image of $\iota_{U,!}$ (resp. $\iota_{U,*}$) is generated by $(\Delta_{x})_{x \in U}$ (resp. $(\nabla_x)_{x \in U}$). 
\item The subcategory $\Hcal_{A, \leq w}(\Gamma_{\Gbb})/\Hcal_{A, <w}(\Gamma_{\Gbb})$ is naturally equivalent to $\DD(\RR_A)$. 
\end{enumerate}
\end{thm}

\begin{proof}
This is standard. The first point is proved as in Lemma \ref{lem:compact-generator-of-hecke-cat}. The second and third points follow from exhibiting the left and right orthogonals of $\Hcal_{A,Z}$ as the categories generated by $(\Delta_w)_{w\in U}$ and $(\nabla_w)_{w\in U}$ which follows from Lemma \ref{lem:orthogonality}. The final point follows from the fact that $\Hom(\Delta_w, \nabla_w) = \RR_A$ which is also proved in Lemma \ref{lem:orthogonality}. 
\end{proof}

\subsection{The reduced Hecke category}

\subsubsection{The subring of $\Wbf$-invariants}
We consider the ring of invariants $\RR_{A}^{\Wbf}$. Theorem \ref{thm:completed-pittie-steinberg} is a strengthening of \cite[Theorem 8.1]{BezrukavnikovRiche}.

\begin{defi}
We denote by $J_A \subset A[X^{\vee}]^{\Wbf}$ the ideal
$$J_A = I_A \cap A[X^{\vee}]^{\Wbf} = \ker(A[X^{\vee}]^{\Wbf} \to A).$$ 
\end{defi}

\begin{thm}[\cite{PittieSteinberg}]\label{thm:not-completed-pittie-steinberg}
\begin{enumerate}
\item The $\Z[X^{\vee}]^{\Wbf}$-module $\Z[X^{\vee}]$ is free of rank $|\Wbf|$.
\item Choose fundamental coweights $\lambda_1, \dots, \lambda_n \in X^{\vee}$ such that for all simple roots $\alpha_j \in X$, we have 
$$\langle \lambda_i, \alpha_j \rangle = \delta_{ij},$$
and choose $z_1, \dots, z_r$ to be a basis of $(X^{\vee})^{\Wbf}$. Then we have 
$$\Z[X^{\vee}]^{\Wbf} \cong \Z[e^{\pm z_i}, f_i]$$
where $f_i = \sum_{\mu \in \Wbf.\lambda_i} e^{\mu}$.
\end{enumerate}
\end{thm}

\begin{thm}[Completed Pittie--Steinberg theorem]\label{thm:completed-pittie-steinberg}
The ring $\RR_{A}$ is finite free of rank $|\Wbf|$ over $\RR_{A}^{\Wbf}$.
\end{thm}

\begin{proof}
The proof is a generalization of \cite[Theorem 8.1]{BezrukavnikovRiche}. As in \emph{loc. cit.}, we have 
$$\Z[X^{\vee}] \otimes_{\Z} A = A[X^{\vee}], \qquad \Z[X^{\vee}]^{\Wbf} \otimes_{\Z} A = A[X^{\vee}]^{\Wbf}.$$
Indeed, $\Z[X^{\vee}]^{\Wbf}$ has a $\Z$-basis whose terms are $e_{\ofrak} = \sum_{\lambda \in \ofrak} e^{\lambda}$ where $\ofrak \subset X^{\vee}$ ranges through the set of $\Wbf$-orbits. The same holds for $A[X^{\vee}]^{\Wbf}$ and the basis over $\Z$ visibly reduces to the one over $A$. 

By Lemma \ref{lem:comparison-of-topologies}, there exists $N \geq 1$ such that 
$$I_A^NA[X^{\vee}] \subset J_AA[X^{\vee}] \subset I_AA[X^{\vee}].$$
This implies that the $I_A$-adic and $J_A$-adic topologies on $A[X^{\vee}]$ coincide, hence $\RR_A$ is also the completion of $A[X^{\vee}]$ for the $J_A$-adic topology. Consequently, we have that $\RR_A$ is free of rank $|\Wbf|$ over $A[X^{\vee}]^{\Wbf, \wedge}_{J_A}$, the $J_A$-adic completion of $A[X^{\vee}]^{\Wbf}$. Since $\Wbf$ acts trivially on $A[X^{\vee}]^{\Wbf}$, there is a canonical injective map 
$$A[X^{\vee}]^{\Wbf,\wedge}_{J_A} \to \RR_A^{\Wbf}.$$
We claim that this map is an isomorphism. 

Consider $D_{w_0} : A[X^{\vee}] \to A[X^{\vee}]^{\Wbf}$ the Demazure operator attached to $w_0$, see Definition \ref{defi:demazure-operator-general}. This operator is $A[X^{\vee}]^{\Wbf}$-linear and passing to $J_A$-adic completion we get an $A[X^{\vee}]^{\Wbf, \wedge}_{J_A}$-linear map
$$\RR_A \to A[X^{\vee}]^{\Wbf, \wedge}_{J_A}.$$

By Lemma \ref{lem:frobenius-extension}, the operator $D_{w_0}$ acts by the identity on $\RR_A^{\Wbf}$, hence the map $\RR_A^{\Wbf} \xrightarrow{D_{w_0}} A[X^{\vee}]^{\Wbf, \wedge}_{J_A}$ is injective. The composition 
$$A[X^{\vee}]^{\Wbf, \wedge}_{J_A} \to \RR_A^{\Wbf} \xrightarrow{D_{w_0}} A[X^{\vee}]^{\Wbf, \wedge}_{J_A}$$
is the identity; it follows that both maps are isomorphisms.
\end{proof}

\begin{lem}\label{lem:comparison-of-topologies}
With the notation of the proof of Theorem \ref{thm:completed-pittie-steinberg}, there exists $N \geq 1$ such that 
$$I_A^NA[X^{\vee}] \subset J_AA[X^{\vee}] \subset I_AA[X^{\vee}].$$ 
Furthermore, we can take $N = \ell(w_0)+1$. 
\end{lem}

\begin{proof}
First note that, as in the proof of Theorem \ref{thm:completed-pittie-steinberg}, all of $J_A, I_A, A[X^{\vee}]$ and $A[X^{\vee}]^{\Wbf}$ are defined over $\Z$, hence it is enough to prove that $I_{\Z}^N \subset J_{\Z}\Z[X^{\vee}]$. It is enough to prove that in 
$$\Cc_{\Wbf, \Z} = \Z[X^{\vee}]/J_{\Z}\Z[X^{\vee}]$$
the image of $I_{\Z}^N$ vanishes, i.e.  for $N$ cocharacters $\lambda_1, \dots, \lambda_N \in X^{\vee}$, we have 
$$\prod_i (e^{\lambda_i} - 1) = 0 \in \Cc_{\Wbf, \Z}.$$
As $\Cc_{\Wbf, \Z}$ is free over $\Z$, we can test that vanishing in $\Cc_{\Wbf, \Q}$. 

By \cite{Borel}, there is an isomorphism of algebras
\begin{equation}\label{eq:iso-to-betti}
\Cc_{\Wbf, \Q} = H^*_{\Betti}(\widehat{\Bcal}, \Q)
\end{equation}
where $\widehat{\Bcal}$ is the flag variety of the Langlands dual group $\widehat{\Gbb}$ of $\Gbb$ defined over $\C$. For any $\lambda \in X^{\vee}$, we have a line bundle $\Ocal(\lambda)$ on $\widehat{\Bcal}$ and under the isomorphism \eqref{eq:iso-to-betti}, the element $e^{\lambda}$ corresponds to $\exp(c_1(\Ocal(\lambda)))$ where $c_1$ is the first Chern class of $\Ocal(\lambda)$. We thus have $e^{\lambda} - 1 \in H^{\geq 2}_{\Betti}(\widehat{\Bcal}, \Q)$, therefore the product of $N$ elements of this form lies in $H^{\geq 2N}_{\Betti}(\widehat{\Bcal}, \Q) = 0$ as $N = \dim(\widehat{\Bcal}) + 1$.
\end{proof}

\begin{lem}\label{lem:invariant-ring-is-a-power-series}
Using the notation of Theorem \ref{thm:not-completed-pittie-steinberg}, the ideal $J_A$ is generated by $(f_i - c_i, e^{z_j} - 1)$ and we have 
$$\RR_A^{\Wbf} = A\llbracket e^{z_i}-1, f_i - c_i \rrbracket,$$
where $c_i \in \Z$ is the image of $f_i$ under the augmentation $\Z[X^{\vee}] \to \Z$. 
\end{lem}

\begin{proof}
By Theorem \ref{thm:not-completed-pittie-steinberg}, we have $A[X^{\vee}]^{\Wbf} = A[e^{\pm z_j}, f_i]$ and it thus follows that $J_A$ is generated by $(e^{z_j}-1, f_i - c_i)$. As in the proof of Theorem \ref{thm:completed-pittie-steinberg}, we have
$$\RR_A^{\Wbf} = A[X^{\vee}]^{\Wbf, \wedge}_{J_A}$$
and the lemma follows. 
\end{proof}

\begin{rque}
The principal difference between Theorem \ref{thm:completed-pittie-steinberg} and \cite[Theorem 8.1]{BezrukavnikovRiche} is that they work when $A$ is a field. In the course of their proof they use the Artin--Rees lemma which would force $A$ to be noetherian. The rings $B$ and $\Ccal(\Zl^{\times}, \Fl)$ are not noetherian and we will use the latter as a specialization of the theory in Section \ref{ssec:modular-constancy-of-the-category}.
\end{rque}

\begin{lem}\label{lem:coincidence-of-actions}
The category $\SBim_A(\Gamma_{\Gbb})$ carries a natural action of $\RR_{A}^{\Wbf}$ compatible with the monoidal action. Moreover, all objects of $\SBim_A(\Gamma_{\Gbb})$ are projective over $\RR_A^{\Wbf}$.
\end{lem}

\begin{proof}
The first part of the lemma is a consequence of the following two points: 
\begin{enumerate}
\item For any $M \in \SBim_A(\Gamma_{\Gbb})$ the left and right actions of $\RR_A$ coincide when restricted to $\RR_A^{\Wbf}$, 
\item For any $M,N \in \SBim_A(\Gamma_{\Gbb})$ and any $r \in \RR_A^{\Wbf}, m \in M, n \in N$, we have 
$$(rm) \otimes n = m \otimes (rn) = r(m \otimes n).$$
\end{enumerate}
The first point holds for all Bott--Samelson bimodules of the form $\BS(s)$. Moreover, if point $(1)$ holds for $M$ and $N$, then point $(2)$ holds for $(M,N)$ and point $(1)$ holds for $M \otimes_{\RR_A} N$. It follows that all Bott--Samelson bimodules satisfy point $(1)$ and, after passing to the Karoubi completion, so do all objects of $\SBim_{A}(\Gamma_{\Gbb})$. 

For the second part of the lemma, by Lemma \ref{lem:projectivity-over-rr}, any object of $\SBim_A(\Gamma_{\Gbb})$ is projective over $\RR_A$ as a left module and $\RR_A$ is projective over $\RR_A^{\Wbf}$ by Theorem \ref{thm:completed-pittie-steinberg}. 
\end{proof}

\begin{lem}\label{lem:lifting-to-rr-invariants}
The category $\Hcal_A(\Gamma_{\Gbb})$ naturally lifts to $\Alg(\Prr_{\RR_A^{\Wbf}})$. 
\end{lem}

\begin{proof}
It follows from Lemma \ref{lem:coincidence-of-actions} that $\Hcal_A(\Gamma_{\Gbb})$ is linear over
\[
\Ind(\Ch^b(\Proj(\RR_A^{\Wbf}))) \cong \DD(\RR_A^{\Wbf}),
\]
where $\Proj(\RR_A^{\Wbf})$ is the category of finitely generated projective $\RR_A^{\Wbf}$-modules.
\end{proof}

Let $A  \to A'$ be a map of rings over $B$. There is a well-defined map 
$$\RR_A^{\Wbf} \to \RR_{A'}^{\Wbf}.$$
We set 
$$\Cc_{\Wbf,A} = \RR_A \otimes_{\RR_A^{\Wbf}} A.$$

\begin{lem}\label{lem:central-specialization}
\begin{enumerate}
\item The natural map 
\begin{equation}\label{eq:reduced-algebra-from-total-torus}
\Ocal(\widehat{\Tbf}_A) \otimes_{\Ocal(\widehat{\Tbf}_A\sslash \Wbf)} A \to \RR_A \otimes_{\RR_A^{\Wbf}} A
\end{equation}
is an isomorphism.
\item Let $s \in S$, there are canonical isomorphisms 
$$\RR_A \otimes_{\RR_{A}^{\Wbf}} \RR_{A'}^{\Wbf} \cong \RR_{A'}, \qquad \RR_{A}^s \otimes_{\RR_{A}^{\Wbf}} \RR_{A'}^{\Wbf} \cong \RR_{A'}^s, \qquad \Cc_{\Wbf,A} \otimes_{A} A' = \Cc_{\Wbf,A'}.$$
\item There is a canonical equivalence 
$$\Hcal_{A}(\Gamma_{\Gbb}) \otimes_{\DD(\RR_A^{\Wbf})} \DD(\RR_{A'}^{\Wbf}) \cong \Hcal_{A'}(\Gamma_{\Gbb}).$$
\end{enumerate}

\end{lem}

\begin{proof}
It follows from the proof of Theorem \ref{thm:completed-pittie-steinberg} that there is a basis of $\Ocal(\widehat{\Tbf}_A)$ over $\Ocal(\widehat{\Tbf}_A\sslash \Wbf)$ whose image in $\RR_A$ is a basis of $\RR_A$ over $\RR_A^{\Wbf}$, hence the map \eqref{eq:reduced-algebra-from-total-torus} is an isomorphism.

 The same basis is sent to itself by the base change map $\RR_A^{\Wbf} \to \RR_{A'}^{\Wbf}$ which establishes the first isomorphism in the second point. The third isomorphism follows immediately from the first and \eqref{eq:reduced-algebra-from-total-torus}. For the remaining isomorphism, by \cite{PittieSteinberg}, the ring $A[X^{\vee}]^s$ is also free over $A[X^{\vee}]^{\Wbf}$. The same argument as in the proof of Theorem \ref{thm:completed-pittie-steinberg} shows that $\RR_A^s$ is free over $\RR_A^{\Wbf}$ of the same rank and that there is a basis of $\RR_A^s$ as an $\RR_A^{\Wbf}$-module defined over $\Z$. This basis is sent onto the analogous basis in $\RR_{A'}^s$ and we get an isomorphism $\RR_A^s \otimes_{\RR_A^{\Wbf}} \RR_{A'}^{\Wbf} \cong \RR_{A'}^s$. 
 
Finally, for the third point, consider the functor
\begin{align*}
\Hcal_{A}(\Gamma_{\Gbb}) &\to \Hcal_{A'}(\Gamma_{\Gbb}) \\
M &\mapsto M \otimes_{\RR_A^{\Wbf}} \RR_{A'}^{\Wbf},
\end{align*}
This functor commutes with relative tensor products over $\RR_A^{\Wbf}$. This shows that it is, in particular, a monoidal functor and that Bott--Samelson bimodules are mapped to the corresponding Bott--Samelson bimodules. As these generate the categories $\Hcal_{A}(\Gamma_{\Gbb})$ and $\Hcal_{A'}(\Gamma_{\Gbb})$ respectively, it is enough to show that this functor is fully faithful. Let $\underline{w}_1, \underline{w}_2 \in \Seq(\Wbf,S)$. As both bimodules $\BS_{A}(\underline{w}_i)$ (for $i = 1,2$) are finite free over $\RR_{A}^{\Wbf}$, we have 
\begin{align*}
\Hom_{\Hcal_{A'}(\Gamma_{\Gbb})}&(\BS_A(\underline{w}_1) \otimes_{\RR_A^{\Wbf}} \RR_{A'}^{\Wbf}, \BS_A(\underline{w}_2) \otimes_{\RR_A^{\Wbf}} \RR_{A'}^{\Wbf}) \\
&\cong \Hom_{\Hcal_{A}(\Gamma_{\Gbb})}(\BS_A(\underline{w}_1), \BS_A(\underline{w}_2))\otimes_{\RR_A^{\Wbf}} \RR_{A'}^{\Wbf}.
\end{align*}
It follows that the functor $\Hcal_{A}(\Gamma_{\Gbb}) \otimes_{\DD(\RR_A^{\Wbf})} \DD(\RR_{A'}^{\Wbf}) \to \Hcal_{A'}(\Gamma_{\Gbb})$ given by $M \mapsto M \otimes_{\RR_A^{\Wbf}} \RR_{A'}^{\Wbf}$ is an equivalence.
\end{proof}

\subsubsection{Reduced Soergel bimodules}

\begin{defi}
We define the reduced Hecke category over $A$ to be 
$$\Hcal_{A,\red}(\Gamma_{\Gbb}) := \Hcal_{A}(\Gamma_{\Gbb}) \otimes_{\DD(\RR_A^{\Wbf})} \DD(A)$$
where $\RR_A^{\Wbf} \to A$ is induced by the augmentation map $\RR_A \to A$. 
For a word $\underline{w}$, we shall denote by $\BS_{A,\red}(\underline{w})$ the image of $\BS(\underline{w})$ in the category $\Hcal_{A,\red}(\Gamma_{\Gbb})$.  We similarly write $\mathbb{B}_{A,\red}$ for the image of $\mathbb{B}_A$ in this category.
\end{defi}

We define \emph{reduced} Soergel bimodules in the same way as we define usual Soergel bimodules. We set 
$$\Cc_{\Wbf,A} = \RR_A \otimes_{\RR_A^{\Wbf}} A.$$
This is a ring that carries an action of $\Wbf$.
\begin{enumerate}
\item For a simple reflection $s$, we set $\Cc_{\Wbf, A}^s = \RR_{A}^s \otimes_{\RR_A^{\Wbf}} A$. This is a subring of the actual ring of $s$-invariants on $\Cc_{\Wbf,A}$.
\item We first define reduced Bott--Samelson bimodules for an expression $\underline{w} = (s_1,\dots, s_n)$ as
$$\BS_{A, \red}(\underline{w}) = \Cc_{\Wbf,A} \otimes_{\Cc_{\Wbf, A}^{s_1}} \Cc_{\Wbf,A} \otimes_{\Cc_{\Wbf, A}^{s_2}} \dots  \otimes_{\Cc_{\Wbf, A}^{s_n}} \Cc_{\Wbf,A}.$$
\item We then declare the category of reduced Soergel bimodules to be 
$$\SBim_{A,\red}(\Gamma_{\Gbb}) := \left\langle \BS_{A, \red}(\underline{w}), \underline{w} \in \Seq(\Wbf,S) \right\rangle_{\otimes_{\Cc_{\Wbf,A}}, \oplus, \text{direct summand}}.$$
\end{enumerate}

\begin{lem}\label{lem:structure-reduced-hecke-cat}
\begin{enumerate}
\item There is a canonical isomorphism 
$$\Hom(\BS_{A,\red}(\underline{w}_1), \BS_{A,\red}(\underline{w}_2)) = \Hom(\BS(\underline{w}_1), \BS(\underline{w}_2)) \otimes_{\RR_A^{\Wbf}} A,$$
\item There is a canonical weight structure, in the sense of Bondarko \cite{Bondarko}, on the category $\Hcal_{A,\red}(\Gamma_{\Gbb})^{\omega}$ whose heart is the category of reduced Soergel bimodules. 
\item There is a canonical equivalence 
$$\Hcal_{A,\red}(\Gamma_{\Gbb})^{\omega} \cong \Ch^b(\SBim_{A,\red}(\Gamma_{\Gbb})).$$
\end{enumerate}
\end{lem}

\begin{proof}
The proof is an adaptation of the argument of \cite[Corollary 4.14]{EteveEberhardt} and the preceding lemmas. 
Point $(1)$ follows from the fact that all Bott--Samelson bimodules are finite free over $\RR_A^{\Wbf}$ by Theorem \ref{thm:completed-pittie-steinberg} and $\BS_{A,\red}(\underline{w}) = \BS_{A}(\underline{w}) \otimes_{\RR_A^{\Wbf}} A$ by construction. 

It follows from the definition of $\Hcal_{A, \red}(\Gamma_{\Gbb})$ that the collection of images in $\Hcal_{A, \red}(\Gamma_{\Gbb})$ of Bott--Samelson bimodules forms a collection of compact generators and by the first point their Ext-algebras are concentrated in degree $0$. By \cite[Theorem 4.3.2]{Bondarko}, there exists a unique weight structure on $\Hcal_{A,\red}(\Gamma_{\Gbb})^{\omega}$ whose heart is the idempotent completion of the subcategory of Bott--Samelson bimodules. Hence we have 
$$\Hcal_{A,\red}(\Gamma_{\Gbb})^{\omega, \heartsuit} \cong \SBim_{A,\red}(\Gamma_{\Gbb}).$$
Furthermore, since the Bott--Samelson objects do not have negative extensions, it follows that the weight complex functor 
$$\Hcal_{A,\red}(\Gamma_{\Gbb})^{\omega} \to \Ch^b(\Hcal_{A,\red}(\Gamma_{\Gbb})^{\omega, \heartsuit})$$
is an equivalence. 
This establishes the third point. 
\end{proof}

\begin{rque}
There are well-defined endomorphisms $\psi_t$ and $\sigma$ acting on $\Hcal_{A, \red}$. 
\end{rque}

\begin{lem}\label{lem:base-change-for-reduced-hecke-cat}
Let $A \to A'$ be a map of rings over $B^+$. There is a canonical equivalence 
$$\Hcal_{A,\red} \otimes_{\DD(A)} \DD(A') \cong \Hcal_{A',\red}.$$
\end{lem}

\begin{proof}
There are natural equivalences 
\begin{align*}
\Hcal_{A, \red} \otimes_{\DD(A)} \DD(A') &\cong \Hcal_{A} \otimes_{\DD(\RR_A^{\Wbf})} \DD(A) \otimes_{\DD(A)} \DD(A') \\
&\cong \Hcal_{A} \otimes_{\DD(\RR_A^{\Wbf})} \DD(\RR^{\Wbf}_{A'}) \otimes_{\DD(\RR^{\Wbf}_{A'})} \DD(A') \\
&\cong \Hcal_{A'} \otimes_{\DD(\RR^{\Wbf}_{A'})} \DD(A') \\
&\cong \Hcal_{A',\red},
\end{align*}
where the third equivalence is provided by Lemma \ref{lem:central-specialization} and the other ones come from the commutative diagram 
\[\begin{tikzcd}
	{\RR_{A}^{\Wbf}} & A \\
	{\RR_{A'}^{\Wbf}} & {A'}
	\arrow[from=1-1, to=1-2]
	\arrow[from=1-1, to=2-1]
	\arrow[from=1-2, to=2-2]
	\arrow[from=2-1, to=2-2]
\end{tikzcd}\]
\end{proof}

\begin{rque}
Note that the natural base change map 
$$\Hcal_{A} \otimes_{\DD(A)} \DD(A') \to \Hcal_{A'}$$
is not an equivalence. This is already true when $\Gbb = \Gm$. 
\end{rque} 

\subsection{Descending reduced Soergel bimodules}

\begin{lem}\label{lem:existence-descent-datum}
The pair $(\Hcal_{B, \red}, \psi_t)$ consisting of the reduced Hecke category and the endofunctor $\psi_t$ descends, as a monoidal category with a monoidal endomorphism, to $B_N$ where $N = \ell(w_0)+1$, that is, there exists a monoidal category $\Hcal_{B_N, \red}$ with a monoidal automorphism $\psi_t$ and an equivalence of monoidal categories
$$\Hcal_{B_N, \red} \otimes_{\DD(B_N)} \DD(B) \cong \Hcal_{B,\red}$$
making the following diagram commute: 
\[\begin{tikzcd}
	{\Hcal_{B,\red} } & {\Hcal_{B_N,\red} \otimes_{\DD(B_N)} \DD(B)} \\
	{\Hcal_{B,\red} } & {\Hcal_{B_N,\red} \otimes_{\DD(B_N)} \DD(B).}
	\arrow["{\psi_t}"', from=1-1, to=2-1]
	\arrow["\sim"', from=1-2, to=1-1]
	\arrow["{\psi_t \otimes \id_{\DD(B)}}", from=1-2, to=2-2]
	\arrow["\sim", from=2-2, to=2-1]
\end{tikzcd}\]
\end{lem}

\begin{proof}
The Hecke category and its reduced version are both defined over $\Z$ using the dual torus defined over $\Z$ and by Lemma \ref{lem:base-change-for-reduced-hecke-cat}, we have 
$$\Hcal_{\Z,\red} \otimes_{\Z} B \cong \Hcal_{B, \red}$$
hence it is only a matter of descending $\psi_t$. 

By Lemma \ref{lem:structure-reduced-hecke-cat}, we have 
$$\Hcal_{B,\red} = \Ind \Ch^b(\SBim_{B, \red})$$
hence it is enough to descend $\psi_t$ when it is acting on $\SBim_{B,\red}$. As this is a full subcategory of all $\Cc_{\Wbf, B}$-bimodules, it is enough to descend the action of $\psi_t$ on $\Cc_{\Wbf, B}$ to $\Cc_{\Wbf, B_N}$. Recall that we have 
$$\Cc_{\Wbf, B} = \RR_{B} \otimes_{\RR_B^{\Wbf}} B = B[X^{\vee}]/J_BB[X^{\vee}]$$
where $J_B = \ker(B[X^{\vee}]^{\Wbf} \to B)$. By Lemma \ref{lem:comparison-of-topologies}, $I_B^N \subset J_BB[X^{\vee}]$ hence in $\Cc_{\Wbf, B}$, for all $\lambda \in X^{\vee}$, we have 
$$\psi_t(e^{\lambda}) = \sum_{n = 0}^N \binom{t}{n}(e^{\lambda} - 1)^n.$$
It then follows that $\psi_t$ acting on $\Hcal_{B, \red}$ descends to $B_N$. 
\end{proof}

\subsection{Relative Serre duality}

The goal of this subsection is to prove a version over some base ring $A$ of relative Serre duality for Hecke categories, Theorem \ref{thm:relative-Serre-duality}, initially proved in \cite{HoLiRelativeSerre}. Since we use a more ``algebraic'' version of the Hecke category, we will follow the proof of \cite{SandvikRelativeSerreDuality}. The proof in \emph{loc. cit.} relies essentially on certain relations between standard/costandard and the (relative) full twists and semi-orthogonal decomposition properties of the Hecke categories, all of which hold in our setting. Let $I \subset S$ be a set of simple roots and let $\Lbb_I \subset \Gbb$ be the corresponding generic Levi subgroup. We do not assume here that $\Lbb_I$ is stable under $\phi$; in this section, we only care about the parabolic root subdatum corresponding to $I$. We shall denote by $\Wbf_{I} \subset \Wbf$ the corresponding parabolic subgroup and by 

$$\SBim_{A}(\Gamma_{I}) = \left\langle \BS_A(\underline{w}), \underline{w} \in \Seq(\Wbf_I, I) \right\rangle, \qquad \Hcal_{A}(\Gamma_{I}) := \Ind\Ch^b(\SBim_{A}(\Gamma_{I})).$$
We shall denote by 
$$\iota_I : \Hcal_{A}(\Gamma_{I}) \to \Hcal_{A}(\Gamma_{\Gbb})$$
the inclusion of the parabolic Hecke category into the ambient Hecke category. 

\begin{defi}\label{defi:relative-full-twists}
\begin{enumerate}
\item The full twist $\FT_{A,\Wbf} \in \Hcal_{A}(\Gamma_{\Gbb})$ is defined to be 
$$\FT_{A,\Wbf} = \Delta_{A, w_0} * \Delta_{A, w_0}.$$
\end{enumerate}
Let $I \subset S$. The relative full twist $\FT_{A,I, \Wbf}$ is defined to be 
$$\FT_{A, I, \Wbf} = \FT_{A, \Wbf_I}^{-1} * \FT_{A, \Wbf}.$$
\end{defi}

\begin{lem}
The inclusion $\iota_I : \Hcal_{A}(\Gamma_{I}) \to \Hcal_{A}(\Gamma_{\Gbb})$ has continuous left and right $\Hcal_{A}(\Gamma_{I})$-linear adjoints. 
\end{lem}

\begin{proof}
The linearity of adjoints over $\Hcal_{A}(\Gamma_{I})$ is an immediate consequence of the rigidity of Hecke categories by Lemma \ref{lem:hecke-category-is-rigid}, see \cite[Chap. 1, Lemma 9.3.6]{GaiRozDAG1}. As $\Hcal_{I,A}$ is generated by either of the collections $(\Delta_{A,w})_{w \in \Wbf_I}$ or $(\nabla_{A,w})_{w \in \Wbf_I}$, it follows that $\Hcal_{A}(\Gamma_{I})$ has left and right orthogonals generated by $(\nabla_{A,w})_{w \in \Wbf\setminus \Wbf_I}$ and $(\Delta_{A,w})_{w \in \Wbf\setminus \Wbf_I}$ respectively, all of whose generators are compact. This guarantees the existence and continuity of the desired adjoints. 
\end{proof}

\begin{thm}[Relative Serre duality]\label{thm:relative-Serre-duality}
There are natural isomorphisms of functors
$$\iota^L_I(\FT_{A,I,\Wbf} * -) \cong \iota^R_I \cong \iota^L_I(- * \FT_{A,I,\Wbf}).$$
\end{thm}

\begin{proof}[Sketch of proof]
Let us indicate the main steps of the proof of \cite{SandvikRelativeSerreDuality, HoLiRelativeSerre} and explain why those steps hold in our situation. 
\begin{enumerate}
\item[Step 1.] Prove that both functors $\FT_{A,I,\Wbf} * -$ and $-*\FT_{A,I,\Wbf}$ induce equivalences $\ker(\iota_I\iota_I^R) \to \ker(\iota_I\iota_I^L)$. The arguments of \emph{loc. cit.} Lemma 3.1 (for $I = \emptyset$) and Lemma 3.2 (for general $I$) rely on the following main points: the semi-orthogonal decomposition of the Hecke category, established in Theorem \ref{thm:semi-orthogonal-decomposition-hecke-cat}, the braid relations for Rouquier complexes, see Lemma \ref{lem:Rouquier-complexes-braid-relations} and the fact that multiplication by $w_0$ is order reversing in the Bruhat order. 
\item [Step 2.] Prove that for $F \in \Hcal_{A,Z}$, where $Z = \{x < w_Iw_0\}$ we have $\iota_I^L(\Delta_{w_I} * F * \Delta_{w_0}) = 0$. This is Lemma 3.3 of \emph{loc. cit.}. This again follows from braid relations and topological considerations on $\Wbf$. 
\item [Step 3.] Prove that there exists a map $\alpha : \FT_{A, \Wbf} \to \FT_{A,\Wbf_{I}}$ which induces an isomorphism after applying $\iota_I^L$. This is Lemma 3.4 of \emph{loc. cit.}. The map is given by
\begin{align*}
\FT_{\Wbf} &= \Delta_{w_0} * \Delta_{w_0} \cong \Delta_{w_I} * \Delta_{w_Iw_0} * \Delta_{w_0} \\
&\to \Delta_{w_I} * \nabla_{w_Iw_0} * \Delta_{w_0} \cong \Delta_{w_I} * \Delta_{w_I} = \iota_I(\FT_{A, \Wbf_I}).
\end{align*}
The fact that this map becomes an isomorphism after applying $\iota_I^L$ follows from the previous point.
\item [Step 4.] Proof of the theorem. As in \emph{loc. cit.}, the proof of the theorem uses only the semi-orthogonal decomposition of $\Hcal_A(\Gamma_{\Gbb})$ and the previous steps. 
\end{enumerate}
\end{proof}

Consider the functor 
$$m_I : \Hcal_A(\Gamma_{\Gbb}) \otimes_{\Hcal_{A}(\Gamma_{I})}  \Hcal_A(\Gamma_{\Gbb}) \to  \Hcal_A(\Gamma_{\Gbb})$$
induced by $X \otimes Y \mapsto X * Y$. 

\begin{lem}\label{lem:relative-serre-duality-for-multiplication}
The functor $m_I^R$ is continuous with continuous right adjoint $m_I^{RR}$, and there is an isomorphism 
$$m_I^{RR}(X \otimes Y) = X * \FT_{A,I, \Wbf}^{-1} * Y.$$
\end{lem}

\begin{proof}
First note that the functor $m_I$ is continuous, preserves compact objects, and is $\Hcal_A(\Gamma_{\Gbb})$-bilinear. By rigidity, the right adjoint $m_I^R$ is also continuous and $\Hcal_A(\Gamma_{\Gbb})$-bilinear. It therefore also has an $\Hcal_A(\Gamma_{\Gbb})$-bilinear right adjoint. By linearity, we have 
$$m_I^R(X) = X * m_I^R(\Delta_{A,1})$$
for all $X \in \Hcal_A(\Gamma_{\Gbb})$. Let us analyze $m_I^R(\Delta_{A,1})$. For $w,v \in \Wbf$, we have 
$$
\Hom(\Delta_{A,w} \boxtimes \Delta_{A,v}, m_I^R(\Delta_{A,1})) \cong \Hom(\Delta_{A,w} * \Delta_{A,v}, \Delta_{A,1}) \cong \Hom(\Delta_{A,v}, \nabla_{A,w^{-1}}).
$$
It follows that $m_I^R(\Delta_{A,1})$ has a filtration by finitely many objects of the form $\Delta_{A,v} \boxtimes \Delta_{A,w}$ and is, in particular, compact. It follows that $m_I^R$ preserves compact objects; therefore $m_I^{RR}$ is also continuous. We can now compute: 
$$m_I^{RR}(X \otimes Y) = X * m_I^{RR}(\Delta_1 \otimes \Delta_1) * Y.$$
Let us now compute $m_{I}^{RR}(\Delta_1 \otimes \Delta_1)$. Consider the functor $\iota_I : \Hcal_A(\Gamma_I) \to \Hcal_A(\Gamma_{\Gbb})$ and extend scalars along $\otimes_{\Hcal_{A}(\Gamma_I)} \Hcal_{A}(\Gamma_{\Gbb})$, which defines a functor
$$\tilde{\iota}_I : \Hcal_A(\Gamma_{\Gbb}) \to \Hcal_A(\Gamma_{\Gbb}) \otimes_{\Hcal_{A}(\Gamma_{I})}  \Hcal_A(\Gamma_{\Gbb}), X \mapsto X \otimes \iota_I(\Delta_1).$$
The composition $m_I \circ \tilde{\iota_I}$ is isomorphic to $\id_{\Hcal_A(\Gamma_{\Gbb})}$, so we have 
$$m_I^{RR} \circ \tilde{\iota}_I^{RR} \cong \id_{\Hcal_A(\Gamma_{\Gbb})}$$
and evaluating on $\Delta_{A,1}$ we get 
$$\tilde{\iota}_I^{RR}(\Delta_{A,1}) = \Delta_{1,A} \otimes \FT_{A,I,\Wbf}$$ 
by Theorem \ref{thm:relative-Serre-duality} and thus 
$$\Delta_{1,A} \cong m_I^{RR}(\Delta_{1,A} \otimes \Delta_{A,1}) *  \FT_{A,I,\Wbf}$$
and the lemma follows. 
\end{proof}

\begin{corol}[\protect{\cite[Corollary 3.1]{SandvikRelativeSerreDuality}}]\label{corol:centrality-full-twist}
The full twist $\FT_{A, \Wbf}$ is central in $\Hcal_{A}(\Gamma_{\Gbb})$.
\end{corol}

\begin{rque}
As before, the setting of the paper \cite{SandvikRelativeSerreDuality} is slightly different from ours but the proof of \emph{loc. cit.} holds in greater generality. 
\end{rque}

\begin{lem}\label{lem:serre-duality-hecke-category}
Let $X,Y \in \Hcal_A(\Gamma_{\Gbb})$. There is a functorial isomorphism 
$$\Hom_{\RR_A}(\Hom_{\Hcal_A}(X,Y), \RR_A) \cong \Hom_{\Hcal_A}(Y, * \FT_{A,\Wbf}^{-1})$$
where $\Hom_{\Hcal_A}(X,Y)$ is viewed as an $\RR_A$-module via the left action of $\RR_A$.
\end{lem}

\begin{proof}
There is a sequence of isomorphisms 
\begin{align*}
\Hom_{\RR_A}(\Hom_{\Hcal_A}(X,Y), \RR_A) &\cong \Hom_{\RR_A}(\Hom_{\Hcal}(1_{\Hcal}, X^{\vee} * Y), \RR_A) \\
&\cong \Hom_{\RR_A}(\iota_{\emptyset}^R(X^{\vee} * Y), \RR_A) \\
&\cong \Hom_{\RR_A}(\iota_{\emptyset}^L(X^{\vee} * Y * \FT_{A, \Wbf}), \RR_A) \\
&\cong \Hom_{\Hcal_A}(X^{\vee} * Y * \FT_{A, \Wbf}, 1_{\Hcal}) \\
&\cong \Hom_{\Hcal_A}(Y, * \FT_{A,\Wbf}^{-1}),
\end{align*}
where the first and last ones follow from dualizing $X,Y$ and $\FT_{A,\Wbf}$, the second and fourth come from the adjunctions $\iota_{\emptyset}^{L} \dashv \iota_{\emptyset} \dashv \iota_{\emptyset}^R$ and the third one from Theorem \ref{thm:relative-Serre-duality}.
\end{proof}

\section{Categorical generic representations}

\subsection{The generic representation category}\label{sec:generic-representation}

\begin{defi}\label{defi:generic-representation-category}
We define the category 
$$\DD(\Rep_{B}^{\unip}(\Gbb(t))) := \Tr(\restr_{\psi_t \circ \sigma}, \Prr_{\Hcal_B})$$
as the categorical trace of $\psi_t \circ \sigma$ acting on the Hecke category.
\end{defi}

\begin{rques}
\begin{enumerate}
\item We write $\DD(\Rep_{B}^{\unip}(\Gbb(t)))$ to suggest that this category is to be thought of as the derived category of unipotent representations of some generic group $\Gbb(t)$ which depends on the parameter $t$. It is, at this point, not clear that this category comes equipped with a $t$-structure. 
\item More generally, for $A$ a $B$-algebra, we denote by 
$$\DD(\Rep_{A}^{\unip}(\Gbb(t))) := \Tr(\restr_{\psi_t \circ \sigma}, \Prr_{\Hcal_A})$$
the category defined with $A$ in place of $B$. 
\end{enumerate}
\end{rques}

Since $\Hcal_A$ carries an action of $\RR_{A}^{\Wbf}$, the category $\DD(\Rep_{A}^{\unip}(\Gbb(t)))$ is naturally linear over 
$$\HHH(\RR_{A}^{\Wbf}, \psi_t \circ \sigma) := \HHH(\RR_{A}^{\Wbf}, \RR_{A, \psi_t \circ \sigma}^{\Wbf}).$$

By functoriality of Hochschild homology, the map 
$$\RR_A^{\Wbf} \to A$$
induced by the augmentation yields a map 
$$\HHH(\RR_{A}^{\Wbf}, \psi_t \circ \sigma) \to A.$$

\begin{defi}\label{defi:generic-reduced-representation-category}
We define 
$$\DD^{\red}(\Rep_{A}^{\unip}(\Gbb(t))) := \Tr(\restr_{\psi_t \circ \sigma}, \Prr_{\Hcal_{A,\red}(\Gamma_{\Gbb})}).$$
This is an $A$-linear category. We call this category the reduced generic representation category. 
\end{defi}

\begin{lem}\label{lem:specialization-of-trace}
Let $A \to A'$ be a map of $B$-algebras. There are canonical equivalences of categories: 
\begin{enumerate}
\item  $\DD(\Rep_{A}^{\unip}(\Gbb(t))) \otimes_{\HHH(\RR_{A}^{\Wbf}, \psi_t \circ \sigma)}\HHH(\RR_{A'}^{\Wbf}, \psi_t \circ \sigma) \cong \DD(\Rep_{A'}^{\unip}(\Gbb(t))),$
\item $\DD^{\red}(\Rep_{A}^{\unip}(\Gbb(t)))  \cong \DD(\Rep_{A}^{\unip}(\Gbb(t))) \otimes _{\HHH(\RR_{A}^{\Wbf}, \psi_t \circ \sigma)} A,$
\item $\DD^{\red}(\Rep_{A}^{\unip}(\Gbb(t))) \otimes_A A' \cong \DD^{\red}(\Rep_{A'}^{\unip}(\Gbb(t))).$
\end{enumerate} 
\end{lem}

\begin{proof} 

All three points are induced by the same principle. Let $\RR_A^{\Wbf} \to C$ be any map of $A$-algebras equipped with a $\psi_t \circ \sigma$-compatible automorphism $\psi_C : C \to C$. Then there is a canonical equivalence 
$$\DD(\Rep_{A}^{\unip}(\Gbb(t))) \otimes_{\HHH(\RR_{A}^{\Wbf}, \psi_t \circ \sigma)}\HHH(C, \psi_C) \cong \Tr(\restr_{\psi_t \circ \sigma}, \Prr_{\Hcal_{A} \otimes_{\RR_A^{\Wbf}} C}).$$

The proof is a variation of the base change property for Hochschild homology which relies on the symmetric monoidality of the trace construction, see Lemma \ref{lem:hh-is-symmetric-monoidal}. We have 
\begin{align*}
&\DD(\Rep_{A}^{\unip}(\Gbb(t))) \otimes_{\HHH(\RR_{A}^{\Wbf}, \psi_t \circ \sigma)}\HHH(C, \psi_C) \\
&= \varinjlim_{[n] \in \Delta^{\opp}} \DD(\Rep_{A}^{\unip}(\Gbb(t))) \otimes \HHH(\RR_{A}^{\Wbf}, \psi_t \circ \sigma)^{\otimes n} \otimes \HHH(C, \psi_C) \\
&= \varinjlim_{[n] \in \Delta^{\opp}} \Tr(\restr_{\psi_t \circ \sigma}, \Prr_{\Hcal_{A}}) \otimes \Tr(\restr_{\psi_t \circ \sigma}, \Prr_{\DD(\RR_A^{\Wbf})})^{\otimes n} \otimes \Tr(\restr_{\psi_C}, \Prr_{\DD(C)}) \\
&= \varinjlim_{[n] \in \Delta^{\opp}} \Tr(\restr_{\psi_t \circ \sigma}, \Prr_{\Hcal_{A} \otimes  \DD(\RR_A^{\Wbf})^{\otimes n} \otimes \DD(C)}) \\
&= \Tr(\restr_{\psi_t \circ \sigma}, \Prr_{\Hcal_{A} \otimes_{\DD(\RR_A^{\Wbf})} \DD(C)})
\end{align*}
The last equality follows from the commutation of Hochschild homology with sifted colimits, see Lemma \ref{lem:commutation-trace-sifted-colimits}.
\end{proof}

\subsubsection{\'{E}tale realization}

Let us fix $(q, \ell, \Lambda)$ as in \ref{ssec:etale-realization}.
\begin{lem}\label{lem:etale-realization}
There is a natural equivalence
$$\DD(\Rep_{B}^{\unip}(\Gbb(t))) \otimes_{\HHH(\RR_B^{\Wbf}, \psi_t \circ \sigma)} \HHH(\RR_{\Lambda}^{\Wbf}, \psi_q \circ \sigma) \cong \DD(\Rep^{\unip}_{\Lambda}(\Gbf^{\Frob}))$$
\end{lem}

\begin{proof}
This is a combination of Theorem \ref{thm:equiv-specialization-with-hecke-category}, Lemmas \ref{lem:specialization-of-trace} and  \ref{lem:specialization-etale-frobenius} together with \cite[Theorem 1.3.1]{EteveFreeMono} which provides an equivalence 
$$\Tr(\Frob^*, \Prr_{\widehat{\DD}^{\fm}(\Ubf \backslash \Gbf/\Ubf, \Lambda)}) \cong \DD(\Rep^{\unip}_{\Lambda}(\Gbf^{\Frob})).$$
\end{proof}

It follows from Lemma \ref{lem:specialization-of-trace} that the category  $\DD(\Rep_{B}^{\unip}(\Gbb(t)))$ interpolates between the various categories of unipotent representations of $\Gbb(q)$ as $q$ varies. 

\subsection{A generator of the generic representation category}\label{sec:generator-in-trace}

Let us denote by 
$$\tr_A : \Hcal_{A} \to \DD(\Rep_{A}^{\unip}(\Gbb(t)))$$ 
the canonical trace map and by 
$$\tr_{A,\red} : \Hcal_{A, \red} \to \DD^{\red}(\Rep_{A}^{\unip}(\Gbb(t)))$$
the reduced version. 

\begin{lem}\label{lem:generator-of-trace}
The object $\tr_A(\mathbb{B}_A)$ (resp. $\tr_{A, \red}(\mathbb{B}_{A, \red})$) is a generator of $\DD(\Rep_{A}^{\unip}(\Gbb(t)))$ (resp. of $\DD^{\red}(\Rep_{A}^{\unip}(\Gbb(t)))$).
\end{lem}

\begin{proof}
By Lemma \ref{lem:trace-is-computed-via-HH}, we have
$$\DD(\Rep_{A}^{\unip}(\Gbb(t))) = \Hcal_A \otimes_{\Hcal_{A} \otimes \Hcal_A^{\rev}} \Hcal_{A, \psi_t \circ \sigma},$$
the essential image of the map $\tr_A$ (resp. $\tr_{A,\red}$) generates the target category. Since $\mathbb{B}_A$ (resp. $\mathbb{B}_{A, \red}$) is a generator of $\Hcal_A$ (resp. $\Hcal_{A, \red}$), its image under the trace functor is a generator as well. As both categories $\Hcal_{A}$ and $\Hcal_{A, \red}$ are rigid, both functors $\tr_A$ and $\tr_{A,\red}$ preserve compact objects, see \cite[Example 7.90]{ZhuCategoricalTame}.
\end{proof}

We denote by $E_A$ the algebra
$$E_A := \End_{\DD(\Rep_{A}^{\unip}(\Gbb(t)))}(\tr_A(\mathbb{B}_A)) \in \Alg(\DD(A)),$$
and by 
$$E_{A, \red} := \End_{\DD^{\red}(\Rep_{A}^{\unip}(\Gbb(t)))}(\tr_{A, \red}(\mathbb{B}_{A, \red})) \in \Alg(\DD(A)).$$
For $i \in \Z$, we denote by 
$$E^i_{A} = H^i(E_A), \qquad E^i_{A, \red} = H^i(E_{A, \red}).$$

\begin{lem}\label{lem:main-lemma-generator-of-trace}
\begin{enumerate}
\item The algebras $E_A$ and $E_{A, \red}$ are connective, i.e., contained in $\DD^{\leq 0}(A)$.
\item The functor 
$$\DD(\Rep_{A}^{\unip}(\Gbb(t))) \rightarrow \DD(E_A^{\opp}), \qquad M \mapsto \RHom(\tr_A(\mathbb{B}_A), M)$$
is an equivalence; likewise, 
$$\DD^{\red}(\Rep_{A}^{\unip}(\Gbb(t))) \rightarrow \DD(E_{A, \red}^{\opp}), \qquad M \mapsto \RHom(\tr_{A, \red}(\mathbb{B}_{A, \red}), M),$$
is an equivalence. 
\item There is a natural isomorphism 
$$E_{A, \red} = E_A \otimes_{\HHH(\RR_{A}^{\Wbf}, \psi_t \circ \sigma)} A.$$ 
\item Let $A \to A'$ be a map of rings over $B$. There is a natural isomorphism 
$$E_{A, \red} \otimes_{A} A' = E_{A',\red}.$$
\item As an $A$-module, $E_{A,\red}^0$ is of finite presentation.
\item Let $N$ be as in Lemma \ref{lem:existence-descent-datum}. The algebra $E_{B,\red}$ descends to $B_N$. 
\end{enumerate}
\end{lem}

\begin{proof}
We will compute the algebras $E_A$ and $E_{A,\red}$ in terms of categorical traces. Let $\EE_A = \End_{\Hcal_A}(\mathbb{B}_A)$ (resp. $\EE_{A,\red} = \End_{\Hcal_{A,\red}}(\mathbb{B}_{A,\red})$). Since $\mathbb{B}_A$ (resp. $\mathbb{B}_{A, \red}$) generates $\Hcal_A$ (resp. $\Hcal_{A,\red}$) we have 
\begin{equation}\label{eq:proof-structure-generic-rep-lem}
\Hcal_A \cong \DD(\EE_A^{\opp}), \qquad \Hcal_{A, \red} \cong \DD(\EE_{A, \red}^{\opp}).
\end{equation}

By Proposition \ref{prop:computing-in-the-trace}, we have canonical isomorphisms in $\DD(A)$ 
$$E_A \cong \Tr(\mathbb{B}_A * \psi_t \circ \sigma (-) * \mathbb{B}_{A}^{\vee}, \Hcal_A) \qquad E_{A,\red} \cong \Tr(\mathbb{B}_{A,\red} * \psi_t \circ \sigma (-) * \mathbb{B}_{A,\red}^{\vee}, \Hcal_{A,\red}).$$

Under the equivalence \eqref{eq:proof-structure-generic-rep-lem}, the endofunctor $\mathbb{B}_A * \psi_t \circ \sigma (-) * \mathbb{B}_{A}^{\vee}$ (resp. $\mathbb{B}_{A,\red} * \psi_t \circ \sigma (-) * \mathbb{B}_{A,\red}^{\vee}$) corresponds to a $\EE_A$-bimodule $\Mbb_A$ (resp. to a $\EE_{A, \red}$-bimodule $\Mbb_{A,\red}$). It then follows that we have isomorphisms in $\DD(A)$ 
\begin{align*}
E_A &= \HHH(\EE_A, \Mbb_A), \\
E_{A,\red} &= \HHH(\EE_{A, \red}, \Mbb_{A,\red}).
\end{align*}

We now claim that $\Mbb_A$ (resp. $\Mbb_{A, \red}$) is in degree $0$ and projective of finite presentation as a left and right $\EE_A$-module (resp. $\EE_{A,\red}$-module). It is the bimodule corresponding to the endofunctor $\mathbb{B}_A * \psi_t \circ \sigma (-) * \mathbb{B}_{A}^{\vee}$. This is the composition of the three functors 
$$X \mapsto \psi_t \circ \sigma(X), \qquad X \mapsto \mathbb{B}_A * X, \qquad X \mapsto X *\mathbb{B}_A^{\vee}.$$
Hence it is enough to show that each of them is represented by an $\EE_A$-bimodule that is in degree $0$ and projective of finite presentation as a left and right module. As $\psi_t \circ \sigma$ is an isomorphism, the first one is clear. As $\mathbb{B}_A^{\vee}$ is also a Soergel bimodule, it is enough to prove that convolution with any Soergel bimodule is represented by a degree $0$ $\EE_A$-bimodule that is projective of finite presentation on both sides. As $\SBim_A = \Add(\mathbb{B}_A)$ by Lemma \ref{lem:additive-generator} and the equivalence \eqref{eq:proof-structure-generic-rep-lem} induces an equivalence
$$\Add(\mathbb{B}_A) \cong \Proj(\EE_A^{\opp}),$$
we obtain that convolution with a Soergel bimodule is represented by an $\EE_A$-bimodule that maps projectives of finite presentation in degree $0$ to projectives of finite presentation in degree $0$, hence it is represented by a bi-projective bimodule in degree $0$ of finite presentation. The same holds for the reduced versions. 

As $\HHH(\EE_A, \Mbb_A)$ is computed by the following bar complex 
$$ \EE_A^{\otimes_A n} \otimes_A \Mbb_A \to \dots \to \EE_A \otimes_A \Mbb_A \to \Mbb_A,$$
it then follows that $E_A$ is in $\DD^{\leq 0}(A)$. The same holds for $E_{A,\red}$. This establishes point $(1)$. 

Point $(2)$ is immediate given that $\tr_A(\mathbb{B}_A)$ (resp. $\tr_{A, \red}(\mathbb{B}_{A, \red})$) is a generator by Lemma \ref{lem:generator-of-trace}.

We have canonical isomorphisms
$$\EE_A \otimes_{\RR_A^{\Wbf}} A \cong \EE_{A, \red}, \qquad \Mbb_A \otimes_{\RR_A^{\Wbf}} A \cong \Mbb_{A, \red}.$$
Point $(3)$ then follows from the base change properties for Hochschild homology, as in the proof of Lemma \ref{lem:specialization-of-trace}.

Point $(4)$ is the same base change property using that 
$$\Hcal_{A, \red} \otimes_{\DD(A)} \DD(A') \cong \Hcal_{A', \red}$$
by Lemma \ref{lem:base-change-for-reduced-hecke-cat} and we thus have
$$\EE_{A, \red} \otimes_A A' \cong \EE_{A', \red}, \qquad \Mbb_{A,\red} \otimes_A A' \cong \Mbb_{A',\red}.$$

Point $(5)$ follows from the fact that, by the bar complex calculation of $E_{A,\red}$, we have 
$$E^0_{A,\red} = \coker(\EE_{A,\red} \otimes_A \Mbb_{A,\red} \to \Mbb_{A,\red})$$
and both $\Mbb_{A,\red}$ and $\EE_{A,\red}$ are finite free $A$-modules as $\Mbb_A$ and $\EE_A$ are finite free $\RR_A^{\Wbf}$-modules by Theorem \ref{thm:completed-pittie-steinberg}. 

Finally, point $(6)$ is again an application of base change for Hochschild homology and the existence of the pair $(\Hcal_{B_N, \red}, \psi_t \circ \sigma)$. 
\end{proof}

\subsubsection{\'{E}tale specialization}

We fix $(q, \ell)$ as an étale realization datum. 

\begin{defi}
Let $(q,\ell, \Zl)$ be as in \ref{ssec:etale-realization}. We denote by $E_{q,\ell}$ the algebra 
$$E_{q,\ell} := \End_{\Rep_{\Zl}^{\unip}(\Gbb(q))}(\tr_{q,\ell}(\mathbb{B}_{\Zl})),$$
by 
$$E_{q, \ell, \red}^0 := H^0(E_{q,\ell} \otimes_{\Ocal((\widehat{\Tbf}_{\Zl}\sslash \Wbf)^{\Frob}_{(1)})} \Zl),$$
and 
$$\overline{E}^0_{q,\ell, \red} := E^0_{q,\ell, \red}/E_{q,\ell, \red}^0[\ell^{\infty}].$$
We also denote by 
$$E_{q, \ell, \Ql} := \End_{\Rep_{\Ql}^{\unip}(\Gbb(q))}(\tr_{q,\ell}(\mathbb{B}_{\Ql})).$$
\end{defi}

\begin{rque}\label{rmk:tilting-sheaves-are-projective}
Note that by \cite[Theorem 1.2.4]{EteveTilting} and \cite[Theorem 4.91]{ZhuCategoricalTame}, the object $\tr_{q,\ell}(\mathbb{B}_{\Zl})$ is a projective module and $E_{q,\ell}$ is concentrated in degree $0$. The ring $E_{q, \ell, \Ql}$ is also concentrated in degree $0$. 
\end{rque} 

\begin{lem}\label{lem:etale-specialization-for-representations}
There is a canonical equivalence for $\Lambda \in \{\Fl, \Zl,\Ql\}$
$$\DD(\Rep_{\Lambda}^{\unip}(\Gbb(q))) \cong \DD(\Rep_{\Lambda}^{\unip}(\Gbf^{\Frob})).$$
\end{lem}

\begin{proof}
This is a formal consequence of Theorem \ref{thm:equiv-specialization-with-hecke-category}, Lemma \ref{lem:specialization-etale-frobenius} and \cite[Theorem 1.3.1]{EteveFreeMono}.
\end{proof}

\begin{lem}\label{lem:control-of-ell-torsion}
\begin{enumerate}
\item As an object of $\Rep_{\Ql}(\Gbf^{\Frob})$, the object $\tr_{q,\ell}(\mathbb{B}_{\Zl})[\frac{1}{\ell}]$ decomposes as 
$$\tr_{q,\ell}(\mathbb{B}_{\Zl})[\frac{1}{\ell}] = \tr_{q,\ell}(\mathbb{B}_{\Ql}) \oplus \text{non-unipotent terms}.$$
\item There are isomorphisms 
$$E^0_{q,\ell, \red}[\frac{1}{\ell}] = E_{q, \ell, \Ql}$$
and
$$\overline{E}^0_{q,\ell, \red} = \mathrm{im}(E_{q,\ell} \to E_{q, \ell, \Ql})$$
where the map $E_{q,\ell} \to E_{q, \ell, \Ql}$ is induced by inversion of $\ell$ and projection to the unipotent part.
\end{enumerate}
\end{lem}

\begin{proof}
The first part of the lemma is proved as in \cite[Lemma 2.4.6]{EteveTilting} by observing that Bott--Samelson bimodules are mapped to Bott--Samelson bimodules (with the same label) after inversion of $\ell$. 

Recall that we have 
$$\HHH(\RR_{\Zl}^{\Wbf}, \psi_q \circ \sigma) \cong \Ocal((\widehat{\Tbf}_{\Zl}\sslash \Wbf)^{\Frob}_{(1)})$$
the ring of functions on the connected component of $(\widehat{\Tbf}_{\Zl}\sslash \Wbf)^{\Frob}$ containing the point $1 \in \widehat{\Tbf}$. Over $\Ql$, we have 
$$(\widehat{\Tbf}_{\Zl}\sslash \Wbf)^{\Frob}_{(1)}\left[\frac{1}{\ell}\right] = (\widehat{\Tbf}_{\Ql}\sslash \Wbf)^{\Frob}_{(1)} \bigsqcup \text{non-unipotent connected components}$$
and we have a commutative diagram, where the vertical maps are the augmentation maps 
\[\begin{tikzcd}
	{\Ocal((\widehat{\Tbf}_{\Zl}\sslash \Wbf)^{\Frob}_{(1)})} & {\Ocal((\widehat{\Tbf}_{\Ql}\sslash \Wbf)^{\Frob}_{(1)}) \cong \Ql.} \\
	\Zl & \Ql
	\arrow[from=1-1, to=1-2]
	\arrow[from=1-1, to=2-1]
	\arrow[from=1-2, to=2-2]
	\arrow[from=2-1, to=2-2]
\end{tikzcd}\]
It follows that we have 
$$E^0_{q,\ell, \red}[\frac{1}{\ell}] = E_{q, \ell, \Ql}.$$
There is a commutative diagram 
\[\begin{tikzcd}
	{E_{q,\ell}} & {E_{q,\ell, \Ql}} \\
	{E^0_{q,\ell, \red}} & {E_{q,\ell, \Ql}}
	\arrow[from=1-1, to=1-2]
	\arrow[from=1-1, to=2-1]
	\arrow[equals, nfold, from=1-2, to=2-2]
	\arrow[from=2-1, to=2-2]
\end{tikzcd}\]
where the bottom horizontal map becomes the first isomorphism after inverting $\ell$. The commutativity of this diagram implies that we have 
$$\overline{E}^0_{q,\ell, \red} = \mathrm{im}(E_{q,\ell} \to E_{q, \ell, \Ql}).$$
\end{proof}

\subsection{Generic parabolic character sheaves}\label{sec:generic-parabolic-char-sheaves}

Let $I \subset S$. We define a generic analog of parabolic character sheaves of Lusztig \cite{LusztigParabCS}. We refer to \cite[Section 4]{EteveJordan} and \cite{GaitsgoryRozemblyumVarshavsky} for some background on their realization in terms of categorical traces.

\begin{defi}\label{defi:generic-parabolic-character-sheaves}
Consider $\Hcal_A(\Gamma_{\Gbb})$ as an $\Hcal_{A}(\Gamma_{I})$-bimodule via the inclusion $\iota_I$. We define 
$$\DD(\CS_{I,A}^{\unip}) := \Tr(\Hcal_A(\Gamma_{\Gbb})_{\psi_t \circ \sigma}, \Prr_{\Hcal_{A}(\Gamma_I)}).$$
For $I = \emptyset$, we shall denote by $\Ocal^{\DL, \unip}_A$ the category $\DD(\CS_{I,A}^{\unip})$.
\end{defi}

\begin{rque}
For the case where $I = \emptyset$, the category $\Ocal^{\DL, \unip}_A$ is a generic variant of what we called the ``Deligne--Lusztig category $\Ocal^{\DL}$'' in \cite{EteveTilting}.
\end{rque} 

\begin{rque}\label{rque:semi-orthogonal-decompositon-category-o}
It follows from the existence of the semi-orthogonal decomposition of $\Hcal_A(\Gamma_{\Gbb})$ by pieces equivalent to $\DD(\RR_A)$, see Theorem \ref{thm:semi-orthogonal-decomposition-hecke-cat}, that the category $\Ocal^{\DL, \unip}_A$ also carries a semi-orthogonal decomposition as all the functors of the theorem are $\DD(\RR_A)$-bilinear. Furthermore, since all functors $\iota_w^?$ and $\iota_{w,?}$ for $? \in \{!,*\}$ preserve compact objects, an object is compact if and only if its $*$-fibers are compact. 
\end{rque} 

There is a natural lax commutative diagram 
\[\begin{tikzcd}
	{\Prr_{\Hcal_{A}(\Gamma_I)}} & {\Prr_{\Hcal_{A}(\Gamma_{\Gbb})}} \\
	{\Prr_{\Hcal_{A}(\Gamma_I)}} & {\Prr_{\Hcal_{A}(\Gamma_{\Gbb})}}
	\arrow["{\ind_{\iota_I}}", from=1-1, to=1-2]
	\arrow["{[\Hcal_A(\Gamma_{\Gbb})_{\psi_t\circ \sigma}]}"', from=1-1, to=2-1]
	\arrow["{\restr_{\psi_t \circ \sigma}}", from=1-2, to=2-2]
	\arrow[Rightarrow, nfold, from=2-1, to=1-2]
	\arrow["{\ind_{\iota_I}}"', from=2-1, to=2-2]
\end{tikzcd}\]
where the lax map is induced by the multiplication map 
$$m_I : \Hcal_{A}(\Gamma_{\Gbb}) \otimes_{\Hcal_{A}(\Gamma_{I})} \Hcal_{A}(\Gamma_{\Gbb}) \to \Hcal_{A}(\Gamma_{\Gbb}).$$
We get a well-defined map 
$$\ch_I : \DD(\CS_{I,A}^{\unip}) \to \DD(\Rep_A^{\unip}(\Gbb(t))).$$

Let $I \subset S$ and let $v \in \Wbf_{\Gbb}$ be an element such that $v\sigma(I) = I$. Let $\Hcal_A(\Gamma_{\Gbb})^{I,v} \subset \Hcal_A(\Gamma_{\Gbb})$ be the full subcategory generated by 
$$\Hcal_A(\Gamma_{\Gbb})^{I,v} := \left\langle \Delta_w, w \in \Wbf_Iv\Wbf_{\sigma(I)} \right\rangle \subset \Hcal_A(\Gamma_{\Gbb})$$
\begin{lem}
The category $\Hcal_A(\Gamma_{\Gbb})^{I,v}$ is naturally an $\Hcal_A(\Gamma_I)$-sub-bimodule of $\Hcal_A(\Gamma_{\Gbb})_{\psi_t \circ \sigma}$ and is isomorphic, as a bimodule, to $\Hcal_A(\Gamma_I)_{\psi_t \circ v\sigma}$.
\end{lem}

\begin{proof}
Both claims follow from the fact that $\Wbf_Iv\Wbf_{\sigma(I)} = \Wbf_I \cdot v$ as a $\Wbf_I$-bimodule with the right action twisted by $\sigma$. 
\end{proof}

\begin{defi}[Generic Deligne--Lusztig induction]
Consider the map of pairs $(\Hcal_A(\Gamma_I), \Hcal_A(\Gamma_{\Gbb})^{I,v}_{\psi_t \circ \sigma}) \to (\Hcal_A(\Gamma_{\Gbb}), \Hcal_A(\Gamma_{\Gbb})_{\psi_t \circ \sigma})$
induced by the inclusion into $\Hcal_A(\Gamma_{\Gbb})$. This defines a functor 
$$R_{I,v} : \DD(\Rep_{A}^{\unip}(\Lbb_{I,v}(t))) \to \DD(\Rep_{A}^{\unip}(\Gbb(t)))$$
which we call generic Deligne--Lusztig induction. We shall denote by ${^*}R_{I,v}$ the right adjoint of $R_{I,v}$.
\end{defi}

\subsubsection{Full twist acting on parabolic character sheaves}

Recall from Definition \ref{defi:relative-full-twists} that we have an object $\FT_{A,I,\Wbf} \in \Hcal_{A}(\Gamma_{\Gbb})$ and that both functors $l_{\FT_{A,I,\Wbf}} = \FT_{A,I,\Wbf} * - $ and $r_{\FT_{A,I,\Wbf}} = - * \FT_{A,I,\Wbf}$ are morphisms of $\Hcal_{A}(\Gamma_I)$-bimodules
$$\Hcal_{A}(\Gamma_{\Gbb}) \to \Hcal_{A}(\Gamma_{\Gbb}).$$ 
We therefore have well-defined induced endofunctors 
$$l^{\tr}_{\FT_{A,I,\Wbf}}, r^{\tr}_{\FT_{A,I,\Wbf}} : \DD(\CS_{I,A}^{\unip}) \to \DD(\CS_{I,A}^{\unip}).$$

\begin{lem}\label{lem:adjoints-to-induction}
\begin{enumerate}
\item There is a natural isomorphism of functors $l^{\tr}_{\FT_{A,I,\Wbf}} \cong r^{\tr}_{\FT_{A,I,\Wbf}}$. We shall then write $t_{\FT_{A,I,\Wbf}} = l^{\tr}_{\FT_{A,I,\Wbf}} =  r^{\tr}_{\FT_{A,I,\Wbf}}$.
\item The functor $\ch_I$ has continuous left and right adjoints, denoted by $\hc_I^!$ and $\hc_I^*$, respectively. Furthermore, there is an isomorphism 
$$\hc_{I}^! = t_{\FT_{A,I,\Wbf}} \circ \hc_I^*.$$
\end{enumerate}
\end{lem}

\begin{rque}
This lemma is a generic counterpart to \cite[Corollary 7.5.10]{GaitsgoryRozemblyumVarshavsky} which was established using Braden's hyperbolic localization theorem. 
\end{rque} 

\begin{proof}
Recall that $\FT_{A,I,\Wbf} = \FT_{A,\Wbf} * \FT_{A,\Wbf_I}^{-1}$ and left and right convolution with either factor commute with the action of $\Hcal_{A}(\Gamma_{I})$. Hence these operations induce the corresponding functors 
$$l^{\tr}_{\FT_{A,\Wbf}}, l^{\tr}_{\FT_{A,\Wbf_I}}, r^{\tr}_{\FT_{A,\Wbf}}, r^{\tr}_{\FT_{A,\Wbf_I}}: \DD(\CS_{I,A}^{\unip}) \to \DD(\CS_{I,A}^{\unip}).$$
Since $\FT_{A,\Wbf_I} \in \Hcal_{A}(\Gamma_I)$, we have 
$$l^{\tr}_{\FT_{A,\Wbf_I}} = r^{\tr}_{\FT_{A,\Wbf_I}}$$
by the universality of the trace property and since $\FT_{A,\Wbf}$ is central by Corollary \ref{corol:centrality-full-twist} in $\Hcal_{A}(\Gamma_{\Gbb})$, we have 
$$l^{\tr}_{\FT_{A,\Wbf}} = r^{\tr}_{\FT_{A,\Wbf}},$$
so the first part of the lemma follows.

Let us prove the second half of the lemma. Consider $\Hcal_{A}(\Gamma_{\Gbb}) \otimes_{\Hcal_A(\Gamma_I)} \Hcal_{A}(\Gamma_{\Gbb})$ as an $\Hcal_{A}(\Gamma_{\Gbb})$-bimodule. By the cyclicity of traces, there is an equivalence 
$$\DD(\CS_{I,A}^{\unip}) \cong \Tr([\Hcal_{A}(\Gamma_{\Gbb}) \otimes_{\Hcal_A(\Gamma_I)} \Hcal_{A}(\Gamma_{\Gbb})_{\psi_t \circ \sigma}], \Prr_{\Hcal_A(\Gamma_{\Gbb})}),$$
and under that equivalence the functor $\ch_I$ is induced by the multiplication 
$$m_I : \Hcal_{A}(\Gamma_{\Gbb}) \otimes_{\Hcal_A(\Gamma_I)} \Hcal_{A}(\Gamma_{\Gbb}) \to \Hcal_A(\Gamma_{\Gbb}).$$
The lemma follows by unraveling the comparison between the left and right adjoints of $m_I$ from Lemma \ref{lem:relative-serre-duality-for-multiplication}.
\end{proof}

\subsubsection{Generic Mackey formula for parabolic induction}\label{ssec:mackay-formula}

We fix $I,J \subset S$ which we assume to be stable under $\sigma$ and we denote by $\Lbb_I$ and $\Lbb_J$ the corresponding Levi subgroups. 

\begin{lem}\label{lem:parab-reduction}
Consider the functor 
$$\DD(\Rep_A^{\unip}(\Lbb_I(t))) \xrightarrow{R_I} \DD(\Rep_A^{\unip}(\Gbb(t))).$$
This functor admits isomorphic left and right adjoints.
\end{lem}

\begin{proof}
The left and right adjoints are induced by the left and right adjoints of both functors 
$$\DD(\Rep_A^{\unip}(\Lbb_I(t))) \xrightarrow{\iota_{I}} \DD(\CS_I) \xrightarrow{\ch_I} \DD(\Rep^{\unip}_A(\Gbb(t))).$$
Both functors have left and right adjoints that differ by the action of some relative full twists by Theorem \ref{thm:relative-Serre-duality} and Lemma \ref{lem:adjoints-to-induction} which cancel out. 
\end{proof}

\begin{prop}[Mackey formula]\label{prop:mackay-for-parab-induction}
The functor 
$$\DD(\Rep_A^{\unip}(\Lbb_I(t))) \xrightarrow{R_{I}}  \DD(\Rep_A^{\unip}(\Gbb(t))) \xrightarrow{{^*}R_J} \DD(\Rep_A^{\unip}(\Lbb_J(t)))$$
has a finite filtration indexed by $\Wbf_J^{\sigma} \backslash \Wbf^{\sigma}/\Wbf_{I}^{\sigma}$ with graded pieces isomorphic to 
$$\gr_w({^*}R_JR_{I}) \cong R_{w(I)\cap J \subset J}\Ad(w){^*}R_{I \cap {w^{-1}(J)} \subset I}$$
where $R_{I_1 \subset I_2}$ denotes the functor $\DD(\Rep_A^{\unip}(\Lbb_{I_1}(t))) \to \DD(\Rep_A^{\unip}(\Lbb_{I_2}(t)))$.
\end{prop}

\begin{proof}
The functor $R_I$ is induced by the $(\Hcal_A(\Gamma_{\Gbb}), \Hcal_A(\Gamma_I))$-bimodule given by $\Hcal_A(\Gamma_{\Gbb})$. By rigidity, see Lemma \ref{lem:right-adjoint-to-induction}, its right adjoint is again $\Hcal_A(\Gamma_{\Gbb})$ viewed as a $(\Hcal_A(\Gamma_I),\Hcal_A(\Gamma_{\Gbb}))$-bimodule. The composition ${^*}R_JR_{I}$ is thus controlled by $\Hcal_A(\Gamma_{\Gbb})$ viewed as a $(\Hcal_A(\Gamma_J), \Hcal_A(\Gamma_I))$-bimodule. 

Recall that for any $Z \subset \Wbf$, we denote by $\Hcal_{A,Z} \subset \Hcal_A(\Gamma_{\Gbb})$ the subcategory generated by all $\Delta_x$ for $x \in Z$. Let us partition $\Wbf = \bigsqcup_{w} \Wbf_J w \Wbf_I$ and for $w \in \Wbf_J \backslash \Wbf/\Wbf_I$, we denote by $Z_w$ the $\sigma$-orbit of $\Wbf_J,\Wbf_I$-double cosets containing $w$. These subsets $Z_w$ are ordered using the Bruhat ordering of $\Wbf$ and therefore provide a $\sigma$-stable filtration of $\Hcal_{A}(\Gamma_{\Gbb})$ as a $(\Hcal_A(\Gamma_J), \Hcal_A(\Gamma_I))$-bimodule with graded pieces isomorphic to $\Hcal_{A,Z_w}$. We also have 
$$\Wbf_J^{\sigma} \backslash \Wbf^{\sigma}/\Wbf_{I}^{\sigma} = (\Wbf_J \backslash \Wbf/\Wbf_{I})^{\sigma}$$ 
and we will prove that the functor $\DD(\Rep_A^{\unip}(\Lbb_I(t))) \to \DD(\Rep_A^{\unip}(\Lbb_J(t)))$ induced by $[\Hcal_{A,Z_w}]$ is zero unless $\Wbf_J w \Wbf_I \in (\Wbf_J \backslash \Wbf/\Wbf_{I})^{\sigma}$ in which case it agrees with the functor given by restriction-induction. 

Recall that the functoriality of traces, see Section \ref{sec:funtoriality-of-traces}, is induced by a morphism of pairs 
$$(\Hcal_{A}(\Gamma_I), \psi_t \circ \sigma) \xrightarrow{[\Hcal_{A,Z_w}], \phi = \psi_t \circ \sigma} (\Hcal_{A}(\Gamma_J), \psi_t \circ \sigma)$$
where $\phi : \Hcal_{A,Z_w} \to \Hcal_{A,Z_w}$ is a $\psi_t \circ \sigma$-linear morphism of bimodules. Assume that $\sigma(\Wbf_J w \Wbf_I) \neq \Wbf_J w \Wbf_I$. Then $[\Hcal_{A,Z_w}]$ has a nontrivial stratification induced by the various $(\Wbf_J, \Wbf_I)$-double cosets in $Z_w$. This defines a filtration of the functor $\phi$ with graded pieces isomorphic to 
$$\psi_t \circ \sigma \circ \iota_{i,!}\iota_{i}^*$$
where $\iota_{i,!}$ denotes the inclusion of $\Hcal_{A, \Wbf_J \sigma^i(w) \Wbf_I} \subset \Hcal_A(\Gamma_{\Gbb})$.
Consequently the functor $\Tr([\Hcal_{A,Z_w}], \phi)$ has a filtration with graded pieces given by $\Tr([\Hcal_{A,Z_w}], \psi_t \circ \sigma \circ  \iota_{i,!}\iota_{i}^*)$. Using the fact that $ \iota_{i,!}\iota_{i}^*$ is idempotent and cyclicity of traces, we have 
\begin{align*}
&\Tr([\Hcal_{A,Z_w}], \psi_t \circ \sigma \circ \iota_{i,!}\iota_{i}^*)\\
 &\cong \Tr([\Hcal_{A,Z_w}], \psi_t \circ \sigma \circ  \iota_{i,!}\iota_{i}^* \iota_{i,!}\iota_{i}^*) \\
&\cong \Tr([\Hcal_{A,Z_w}],  \iota_{i,!}\iota_{i}^* \psi_t \circ \sigma \circ  \iota_{i,!}\iota_{i}^*) \\
&\cong \Tr([\Hcal_{A,Z_w}], \psi_t \circ \sigma \circ  \iota_{i+1,!}\iota_{i+1}^* \iota_{i,!}\iota_{i}^*) \\
&= 0,
\end{align*}
as $\iota_{i+1}^* \iota_{i,!} = 0$ by disjointness of $\Wbf_J \sigma^{i+1}(w) \Wbf_I$ and $\Wbf_J \sigma^i(w) \Wbf_I$.

Assume now that $\sigma(\Wbf_J w \Wbf_I) = \Wbf_J w \Wbf_I$. The functor $R_{w(I)\cap J \subset J}\Ad(w)\allowbreak{^*}R_{I \cap {w^{-1}(J)}}$ is induced by the bimodule 
$$\Hcal_A(\Gamma_J) \otimes_{\Hcal_A(\Gamma_{J \cap w(I)})} (\Delta_w * \Hcal_A(\Gamma_{J \cap w(I)})) \otimes_{\Hcal_A(\Gamma_{w^{-1}(J) \cap I})} \Hcal_A(\Gamma_I)$$ 
and multiplication 
$$\Hcal_A(\Gamma_J) \otimes_{\Hcal_A(\Gamma_{J \cap w(I)})} (\Delta_w * \Hcal_A(\Gamma_{J \cap w(I)})) \otimes_{\Hcal_A(\Gamma_{w^{-1}(J) \cap I})} \Hcal_A(\Gamma_I) \to \Hcal_{A,Z_w}$$
is an isomorphism of bimodules compatible with the action of $\psi_t \circ \sigma$.
\end{proof}

\subsubsection{Serre functor}\label{sec:serre-functor}

We introduce the following assumption on $A$: 

\begin{enumerate}
\item[(Fin)] The ring $A$ satisfies the property (Fin) if for all $w \in \Wbf$, the ring 
$$\HHH(\RR_A, \psi_t \circ w \circ \sigma)$$ 
is concentrated in degree $0$ and finite locally free over $A$ and we have, see Section \ref{ssec:HHcalculation}	
$$\HHH(\RR_A, \psi_t \circ w \circ \sigma) \cong \HHH^{\cont}(\RR_A, \psi_t \circ w \circ \sigma).$$ 
\end{enumerate}

\begin{rque}
The condition (Fin) is for instance satisfied when we fix $A = \Lambda$ for any étale realization datum $(q,\ell, \Lambda)$.
\end{rque} 

\begin{lem}\label{lem:base-change-tori}
Assume that $A$ satisfies (Fin). Then for any surjective map $A \to A'$, we have 
$$\HHH(\RR_A, \psi_t \circ \sigma) \otimes_A A' \cong \HHH(\RR_{A'}, \psi_t \circ \sigma).$$
\end{lem}

\begin{proof}
The base change formula holds for continuous Hochschild homology and agrees with usual Hochschild homology by Lemma \ref{lem:reduction-koszul-calculation}.
\end{proof}

Let us denote by 
$$\tr_{A,\emptyset} : \Hcal_{A}(\Gamma_{\Gbb}) \to \Ocal_{A}^{\DL, \unip}$$
the canonical trace map. 

\begin{lem}\label{lem:calculation-tr-vs-category-o}
Assume that $A$ satisfies (Fin). Then for all $X,Y \in \Hcal_A(\Gamma_{\Gbb})^{\omega}$, there is a canonical isomorphism 
$$\Hom_{\Hcal}(X,Y) \otimes_{\RR_A \otimes_A \RR_A} \RR_A \cong \Hom_{\Ocal}(\tr_{A,\emptyset}(X), \tr_{A,\emptyset}(Y))$$
where $\Hom_{\Hcal}(X,Y)$ is viewed as an $\RR_A \otimes_A \RR_A$-module in the usual way and $\RR_A \otimes_A \RR_A \to \RR_A$ is given by $x \otimes y\mapsto x\,\psi_t\circ \sigma(y)$.
\end{lem}

\begin{proof}
This is immediate from the definition of $\Ocal_{A}^{\DL, \unip} = \Hcal_A(\Gamma_{\Gbb}) \otimes_{\DD(\RR_A \otimes_A \RR_A)} \DD(\RR_A)$.
\end{proof}

\begin{lem}\label{lem:serre-duality-on-category-O}
Assume that $A$ satisfies (Fin). 
\begin{enumerate}
\item The category of compact objects $\Ocal^{\DL, \unip, \omega}_A$ of $\Ocal^{\DL, \unip}_A$ is proper over $A$, i.e., for any two objects $X,Y \in \Ocal^{\DL, \unip, \omega}_A$, the mapping complex
$$\Hom_{\Ocal^{\DL, \unip, \omega}_A}(X,Y) \in \Perf(A)$$
is perfect over $A$.
\item There exists a well-defined Serre functor $S_{A,\Ocal} : \Ocal^{\DL, \unip, \omega}_A \to \Ocal^{\DL, \unip, \omega}_A$ which satisfies, for all $X,Y \in \Ocal^{\DL, \unip, \omega}_A$, 
$$\Hom_{\Ocal^{\DL, \unip}_A}(X,Y)^{*} \cong \Hom_{\Ocal^{\DL, \unip}_A}(Y, S_{A, \Ocal}(X)).$$
\item There is an isomorphism of functors 
$$S_{A, \Ocal} \cong t_{\FT_{A,\emptyset,\Wbf}}^{-1}.$$
\end{enumerate}
\end{lem}

\begin{proof}
As $\Hcal_{A}(\Gamma_{\Gbb})$ has a semi-orthogonal decomposition with strata indexed by $\Wbf$ and this semi-orthogonal decomposition is stable under the action of $\DD(\RR_A)$, it follows that $\Ocal^{\DL, \unip}$ has a semi-orthogonal decomposition indexed by $\Wbf$ with strata isomorphic to $\DD(\HHH(\RR_A, \psi_t \circ w \circ \sigma))$. For any pair of compact objects $X,Y \in \Ocal^{\DL, \unip}$, there is a filtration of $\Hom_{\Ocal^{\DL, \unip}}(X,Y)$ with pieces isomorphic to 
$$\Hom(j_w^*X,j_w^!Y)$$
and as both functors $j_w^*, j_w^! : \Hcal_{A}(\Gamma_{\Gbb}) \to \DD(\RR_A)$ preserve compact objects, the corresponding functors $\Ocal^{\DL, \unip} \to \DD(\HHH(\RR_A, \psi_t \circ w \circ \sigma))$ also preserve compact objects. It follows that both $j_w^*X$ and $j_w^!Y$ are perfect complexes of $\HHH(\RR_A, \psi_t \circ w \circ \sigma)$-modules and thus $\Hom(j_w^*X,j_w^!Y)$ is itself a perfect $A$-module by hypothesis (Fin). 

For any $X \in \Ocal_A^{\DL, \unip, \omega}$, the functor $Y \mapsto \Hom_{\Ocal}(X,Y)$ commutes with filtered colimits as $X$ is compact. Let $S_{A,\Ocal}(X) = \Hom(X,-)^{R}(A)$ be the right adjoint of $\Hom(X,-)$ evaluated at $A$. We thus have 
$$\Hom_{\Ocal}(X,Y)^* \cong \Hom_{\Ocal}(Y, S_{A,\Ocal}(X)).$$
Evaluating this at $Y = \Delta_{A,w}$ yields that for all $w$, the object $j_w^*S_{A,\Ocal}(X)$ is compact hence $S_{A,\Ocal}(X)$ is compact. 

Let $X,Y \in \Hcal_A$. We claim that there is a canonical isomorphism 
$$\Hom_A(\Hom_{\Ocal}(\tr_{A,\emptyset}(X), \tr_{A,\emptyset}(Y)), A) \cong \Hom_{\RR_A}(\Hom_{\Hcal_A}(X,Y), \RR_A) \otimes_{\RR_A \otimes_A \RR_A} \RR_A.$$
Assume that this claim holds; then we have 
\begin{align*}
&\Hom_A(\Hom_{\Ocal}(\tr_{A,\emptyset}(X), \tr_{A,\emptyset}(Y)), A) \\
&\cong \Hom_{\RR_A}(\Hom_{\Hcal_A}(X,Y), \RR_A) \otimes_{\RR_A \otimes_A \RR_A} \RR_A \\
&\cong \Hom_{\Hcal}(Y, * \FT_{A,\Wbf}^{-1}) \otimes_{\RR_A \otimes_A \RR_A} \RR_A \\
&\cong \Hom_{\Ocal}(\tr_{A,\emptyset}(Y), \tr_{A,\emptyset}(* \FT_{A,\Wbf}),
\end{align*}
where the second isomorphism comes from Lemma \ref{lem:serre-duality-hecke-category} and the third one from Lemma \ref{lem:calculation-tr-vs-category-o}. 

Let us now prove the claim. For any $X,Y \in \SBim_A(\Gamma_{\Gbb})$, there is a natural 	map 
$$\Hom_{\RR_A}(\Hom_{\Hcal_A}(X,Y), \RR_A) \otimes_{\RR_A \otimes_A \RR_A} \RR_A \to \Hom_A(\Hom_{\Ocal}(\tr_{A,\emptyset}(X), \tr_{A,\emptyset}(Y)), A)$$
as any map $\Hom_{\Hcal_A}(X,Y) \to A$ of $\RR_A$-modules factors through $\Hom_{\Hcal_A}(X,Y) \otimes_{\RR_A \otimes_A \RR_A} \RR_A \to A$. We can now filter $\Hom_{\Hcal_A}(X,Y)$ and replace $\Hom_{\Hcal_A}(X,Y)$ by its graded pieces which are of the form $\RR_{A,w}$ where the claim is now clear. 
\end{proof}

\begin{lem}\label{lem:self-duality}
Assume that $A$ satisfies (Fin).
\begin{enumerate}
\item The algebra $E_A$ is projective, locally of finite presentation and in degree $0$.
\item There is a canonical isomorphism $E_A \cong E_A^{\vee}$.
\end{enumerate}
\end{lem}

\begin{proof}
We prove that $E_A$ is self-dual. As $E_A \in \DD^{\leq 0}(A)$, we have $E_A^{\vee} \in \DD^{\geq 0}(A)$, so self-duality implies that $E_A \in \DD^0(A)$. 

Let us prove self-duality. We compute directly:
\begin{align*}
&\Hom_{\DD(\Rep_A(\Gbb(t)))}(\tr_A(\mathbb{B}_A), \tr_A(\mathbb{B}_A))^{\vee}\\
&\cong \Hom_{\DD(\Rep_A(\Gbb(t)))}(\ch \circ \tr_{A,\emptyset}(\mathbb{B}_A), \ch \circ\tr_{A,\emptyset}(\mathbb{B}_A))^{\vee} \\
&\cong \Hom_{\Ocal}(\hc^!\ch\circ \tr_{A,\emptyset}(\mathbb{B}_A), \tr_{A,\emptyset}(\mathbb{B}_A))^{\vee} \\
&\cong \Hom_{\Ocal}(\tr_{A,\emptyset}(\mathbb{B}_A), S_{A,\Ocal}\hc^!\ch \circ \tr_{A,\emptyset}(\mathbb{B}_A)) \\
&\cong \Hom_{\Ocal}(\tr_{A,\emptyset}(\mathbb{B}_A), \hc^*\ch\circ \tr_{A,\emptyset}(\mathbb{B}_A)) \\
&\cong \Hom_{\DD(\Rep_A(\Gbb(t)))}(\tr_A(\mathbb{B}_A), \tr_A(\mathbb{B}_A)).
\end{align*}
The maps are as follows: 
\begin{enumerate}
\item The first map follows from the fact that $\tr_A = \ch \tr_{A,\emptyset}$, 
\item the second one follows from the $(\hc^!, \ch)$-adjunction,
\item the third one follows from the definition of Serre duality, 
\item the fourth one follows from Lemma \ref{lem:serre-duality-on-category-O} and Lemma \ref{lem:adjoints-to-induction}, 
\item the last one follows from the $(\ch, \hc^*)$-adjunction.
\end{enumerate}

Let us now prove that $E_A$ is flat and locally of finite presentation. Let $A \to A'$ be a surjective map. By Lemma \ref{lem:base-change-tori}, the map 
$$\Ocal^{\DL, \unip}_A \otimes_{\DD(A)} \DD(A') \to \Ocal^{\DL,\unip}_{A'}$$ 
is an equivalence. Indeed, this functor commutes with the semi-orthogonal decomposition, see Remark \ref{rque:semi-orthogonal-decompositon-category-o} and the lemma asserts that it is an equivalence on strata. It then follows that there is an isomorphism 
$$\tr_{A, \emptyset}(\mathbb{B}_A) \otimes_A A' \cong \tr_{A',\emptyset}(\mathbb{B}_{A'}).$$
This implies that we have 
$$E_A \otimes_A A'  \cong E_{A'}$$ 
which is still in degree $0$, and in particular $E_A$ is flat over $A$. We therefore have a finitely presented flat $A$-module which is self-dual; this implies that it is projective. 
\end{proof}

\begin{rque}\label{rque:t-structure}
Under the equivalence of Lemma \ref{lem:main-lemma-generator-of-trace}
$$\DD(\Rep_{A}(\Gbb(t))) \cong \DD(E_A^{\opp})$$ 
and under the assumption (Fin), by Lemma \ref{lem:self-duality}, the algebra $E_A$ lies in degree $0$, we can therefore equip the category $\DD(\Rep_A(\Gbb(t)))$ with a $t$-structure satisfying 
\begin{itemize}
\item the object $\tr_A(\mathbb{B}_A)$ is a compact progenerator of $\DD(\Rep_{A}(\Gbb(t)))^{\heartsuit}$.
\end{itemize}
Whenever the condition (Fin) is satisfied, we shall write $\Rep_A(\Gbb(t)) := \DD(\Rep_{A}(\Gbb(t)))^{\heartsuit}$. This is an abelian category.
\end{rque}

\subsubsection{\'{E}tale realization}
Let us fix an étale realization datum $(q, \ell, \Lambda)$. Let us denote by $\Pbf = \Lbf\Vbf \subset \Gbf$ the standard parabolic subgroup with Levi factor $\Lbf$ and unipotent radical $\Vbf$ corresponding to $I$. 

\begin{lem}\label{lem:etale-realization-gen-parab-cs}
There is a natural equivalence 
$$\DD(\CS_{I,\Lambda}^{\unip}) \cong \DD^{\unip}\left(\frac{\Vbf\backslash \Gbf/\Frob(\Vbf)}{\Ad_{\Frob}(\Lbf)}, \Lambda\right),$$
and under that equivalence the adjunction triple $(\hc_I^!, \ch_I, \hc_I^*)$ corresponds to the adjunction triple $(\hc_{\Pbf \subset \Gbf}^!, \ch_{\Pbf \subset \Gbf}, \hc_{\Pbf \subset \Gbf}^*)$ using the notation of \cite[Section 4]{EteveJordan}.

Let $\dot{v} \in \Gbf$ be a lift of the element $v \in \Wbf$. Then there exists a commutative diagram 
\[\begin{tikzcd}
	{\DD(\Rep_{\Lambda}^{\unip}(\Lbf^{\dot{v}\Frob}))} & {\DD(\Rep_{\Lambda}^{\unip}(\Lbb_{I,v}(t)))} \\
	{\DD(\Rep_{\Lambda}^{\unip}(\Gbf^{\Frob}))} & {\DD(\Rep_{\Lambda}^{\unip}(\Gbb(t)))}
	\arrow["\cong", no head, from=1-1, to=1-2]
	\arrow["{R_{\Lbf, \dot{v}}}", from=1-1, to=2-1]
	\arrow["{R_{I,v}}", from=1-2, to=2-2]
	\arrow["\cong"', no head, from=2-1, to=2-2]
\end{tikzcd}\]
where $R_{\Lbf, \dot{v}}$ denotes the corresponding usual Deligne--Lusztig induction. 
\end{lem}

\begin{proof}
Immediate from Theorem \ref{thm:equiv-specialization-with-hecke-category} and \cite[Theorem 5.1.3]{EteveJordan}.
\end{proof}

\subsection{Generic representations over $B[|\Gbb|^{-1}]$}

We shall write $|\Gbb| \in \Z[t,t^{-1}]$ for the polynomial order of $\Gbb$ as in \cite[Définition 1.9]{BroueMalle}. This is the polynomial satisfying 
$$|\Gbb(q)| = |\Gbb|(q).$$ In this subsection, we let $A$ be any commutative ring over $B[|\Gbb|^{-1}]$.

\begin{lem}\label{lem:calculation-determinant}
We have an equality in $A$ 
$$\det(\id - \psi_t \circ \sigma, J_A/J_A^2) = \pm t^{-|\Phi^+|}|\Gbb|.$$
\end{lem}

\begin{proof}
By Lemma \ref{lem:invariant-ring-is-a-power-series}, we have $J_A/J_A^2 = \bigoplus_{i} A.(e^{z_i} - 1) \oplus \bigoplus_{j} A.(f_j - c_j)$ (using the notation of the lemma). It follows that 
$$\det(\id - \psi_t \circ \sigma, J_A/J_A^2) = \prod_{i,j} \det(\id - \psi_t \circ \sigma, A.(e^{z_i} - 1))\det(\id - \psi_t \circ \sigma, A.(f_j - c_j))$$
which by \cite[Proposition 1.8]{BroueMalle} is $\pm t^{-|\Phi^+|}|\Gbb|$. 
\end{proof}

\begin{lem}\label{lem:generic-ring-of-invariants}
There is a natural isomorphism over $A$, 
$$\HHH(\RR_{A}^{\Wbf}, \psi_t \circ \sigma) \cong A,$$
and it is, in particular, in degree $0$. 
\end{lem}

\begin{proof}
By Lemma \ref{lem:invariant-ring-is-a-power-series}, $\RR_{A}^{\Wbf}$ is a power series ring and we have 
$$\det(\id - \psi_t \circ \sigma, J_A/J_A^2) \in A^{\times}$$
by Lemma \ref{lem:calculation-determinant}. The lemma follows from Lemma \ref{lem:reduction-koszul-calculation}. 
\end{proof}

\begin{corol}\label{cor:generically-category-is-reduced}
Over $B[|\Gbb|^{-1}]$, we have 
$$\DD(\Rep^{\unip}_{B[|\Gbb|^{-1}]}(\Gbb(t))) \cong \DD^{\red}(\Rep^{\unip}_{B[|\Gbb|^{-1}]}(\Gbb(t))).$$
\end{corol}

\begin{proof}
Immediate from Lemma \ref{lem:generic-ring-of-invariants}.
\end{proof}

\begin{lem}\label{lem:finite-etale}
Assume that $A$ satisfies the condition (Fin). Then the ring 
$$\HHH(\RR_A^{\Wbf}, \psi_t \circ \sigma)\left[\frac{1}{|\Gbb|}\right]$$
is in degree $0$ and finite étale over $A$.
\end{lem}

\begin{proof}
As in the proof of Lemma \ref{lem:generic-ring-of-invariants}, identify $\RR_A^{\Wbf} = A\llbracket e^{z_i} - 1, f_i - c_i \rrbracket$, and by assumption on $A$, there exists $N > 0$ such that $J_A^N \subset (\id - \psi_t \circ \sigma)(J_A)$, hence by Lemma \ref{lem:reduction-koszul-calculation}, we find that $\HHH(\RR_A^{\Wbf}, \psi_t \circ \sigma)$ is in degree $0$ and finite over $A$. Let us compute $\Omega^1_{\HHH(\RR_A^{\Wbf}, \psi_t \circ \sigma)/A}$. By the Künneth formula, we have 
$$\Omega^1_{\HHH(\RR_A^{\Wbf}, \psi_t \circ \sigma)/A} = \boxtimes_A \Omega^1_{\HHH(A\llbracket z_i-1\rrbracket, \psi_t \circ \sigma)} \otimes_ A\Omega^1_{\HHH(A\llbracket f_i-c_i \rrbracket, \psi_t \circ \sigma)}.$$
But we then have 
$$ \Omega^1_{\HHH(A\llbracket f_i-c_i \rrbracket, \psi_t \circ \sigma)} = A\llbracket f_i-c_i \rrbracket/(f_i - \psi_t \circ \sigma(f_i), 1 - \frac{d}{d(f_i - c_i)}(\psi_t \circ \sigma(f_i-c_i)))d(f_i - c_i)$$
which vanishes over $A\left[\frac{1}{|\Gbb|}\right]$ by Lemma \ref{lem:calculation-determinant}.
\end{proof}

\begin{lem}\label{lem:compacity-as-a-bimodule}
Assume that $A$ satisfies (Fin) and consider the algebra $E_A$ as an $E_A$-bimodule. Then $E_A[\frac{1}{|\Gbb|}]$ is a compact object of $\DD(E_A[\frac{1}{|\Gbb|}] \otimes E_A[\frac{1}{|\Gbb|}]^{\opp})$. 
\end{lem}

\begin{proof}
Let $\Hcal_{A}^{\rev}(\Gamma_{\Gbb})$ be the monoidal category with underlying category given by $\Hcal_{A}(\Gamma_{\Gbb})$ and tensor product given by $x \otimes^{\rev} y = y \otimes x$. It follows from Lemma \ref{lem:dual-of-categorical-trace} that we have an equivalence in $\Prr_A$
$$\Tr(\restr_{\psi_t \circ \sigma}, \Prr_{\Hcal_{A}^{\rev}(\Gamma_{\Gbb})}) \cong \Tr(\restr_{\psi_t \circ \sigma}, \Prr_{\Hcal_{A}(\Gamma_{\Gbb})})^{\vee}.$$
In particular, we have 
$$\Tr(\restr_{\psi_t \circ \sigma}, \Prr_{\Hcal_{A}^{\rev}(\Gamma_{\Gbb})}) \cong \DD(E_A),$$
by Lemma \ref{lem:main-lemma-generator-of-trace}. 
Consider the category $\Hcal_{A}(\Gamma_{\Gbb})$ as a bimodule over itself, i.e. as an $\Hcal_{A}(\Gamma_{\Gbb}) \otimes \Hcal_{A}(\Gamma_{\Gbb})^{\vee}$-module. There is a well-defined class 
$$\cl(\Hcal_{A}(\Gamma_{\Gbb})) \in \Tr(\restr_{\psi_t \circ \sigma \otimes \psi_t \circ \sigma}, \Prr_{\Hcal_{A}(\Gamma_{\Gbb}) \otimes \Hcal_{A}(\Gamma_{\Gbb})^{\rev}}).$$
Under the equivalence 
$$\Tr(\restr_{\psi_t \circ \sigma \otimes \psi_t \circ \sigma}, \Prr_{\Hcal_{A}(\Gamma_{\Gbb}) \otimes \Hcal_{A}(\Gamma_{\Gbb})^{\rev}}) \cong \DD(E_A \otimes E_A^{\opp})$$
the object $\cl(\Hcal_{A}(\Gamma_{\Gbb}))$ corresponds to $E_A$ viewed as a bimodule by Lemma \ref{lem:dual-of-categorical-trace}. 

Now consider the functor $\ch : \Ocal^{\DL,\unip}_A \to \DD(\Rep_A^{\unip}(\Gbb(t)))$. The left adjoint $\hc^!$ is conservative and detects compactness. The conservativity follows from the fact that the essential image of $\ch$ generates the target category. The fact that this functor detects compactness follows from the continuity of $\ch$. 

The object $\hc^!(\cl(\Hcal_{A}(\Gamma_{\Gbb}))) \in \Ocal^{\DL,\unip}(\Gbb) \otimes \Ocal^{\DL,\unip}(\Gbb)$ can be computed as the class of $\Hcal_{A}(\Gamma_{\Gbb})$ viewed as a $\DD(\RR_A)$-bimodule. Consider the following diagram: 
\[\begin{tikzcd}
	{\Prr_A} & {\Prr_{\Hcal_{A}(\Gamma_{\Gbb}) \otimes \Hcal_{A}(\Gamma_{\Gbb})^{\rev}}} & {\Prr_{\DD(\RR_A) \otimes \DD(\RR_A)}} \\
	{\Prr_A} & {\Prr_{\Hcal_{A}(\Gamma_{\Gbb}) \otimes \Hcal_{A}(\Gamma_{\Gbb})^{\rev}}} & {\Prr_{\DD(\RR_A) \otimes \DD(\RR_A)}}
	\arrow["{[\Hcal_{A}(\Gamma_{\Gbb}) ]}", from=1-1, to=1-2]
	\arrow[equals, nfold, from=1-1, to=2-1]
	\arrow["\restr", from=1-2, to=1-3]
	\arrow["{\restr_{\psi_t \circ \sigma \otimes \psi_t \circ \sigma}}", from=1-2, to=2-2]
	\arrow["{[\Hcal_{A}(\Gamma_{\Gbb}) \otimes \Hcal_{A}(\Gamma_{\Gbb})^{\rev}]}", from=1-3, to=2-3]
	\arrow[Rightarrow, nfold, from=2-1, to=1-2]
	\arrow["{[\Hcal_{A}(\Gamma_{\Gbb}) ]}"', from=2-1, to=2-2]
	\arrow[Rightarrow, nfold, from=2-2, to=1-3]
	\arrow["\restr"', from=2-2, to=2-3]
\end{tikzcd}\]
where the $2$-maps are induced by $\psi_t \circ \sigma$ and the multiplication in $\Hcal_{A}(\Gamma_{\Gbb})$ respectively. By Remark \ref{rque:semi-orthogonal-decompositon-category-o}, it is enough to check that for all $w,w' \in \Wbf$ the fibers $j_{w,w'}^*\hc^!(\cl(\Hcal_{A}(\Gamma_{\Gbb})))$ are compact in $\DD(\RR_A \otimes \RR_A)$ (after inverting $|\Gbb|$). This complex is computed as the class of $\Hcal_A(\Gamma_{\Gbb})$ viewed as a $\DD(\RR_A)$-module via 
$$\DD(\RR_A) \otimes_A \Hcal_A(\Gamma_{\Gbb}) \otimes \DD(\RR_A) \to \Hcal_A(\Gamma_{\Gbb}), \RR_A \otimes X \otimes \RR_A \mapsto \Delta_{w,A} * X * \Delta_{w',A},$$
and whose equivariance data are 
$$\psi_{t,w,w'} :  \Hcal_A(\Gamma_{\Gbb}) \to  \Hcal_A(\Gamma_{\Gbb}), X \mapsto \Delta_w * \psi_t \circ \sigma (\Delta_w^{-1} * X * \Delta_{w'}^{-1}) * \Delta_{w'}$$ 
so that we have a commutative diagram 
\[\begin{tikzcd}
	{\Prr_A} & {\Prr_{\DD(\RR_{A} \otimes \RR_A)}}  \\
	{\Prr_A} & {\Prr_{\DD(\RR_{A} \otimes \RR_A)}} 
	\arrow["{[\Hcal_A(\Gamma_{\Gbb})]}", from=1-1, to=1-2]
	\arrow[equals, nfold, from=1-1, to=2-1]
	\arrow["{\restr_{\psi_t \circ w\sigma \otimes \psi_t \circ w'\sigma}}", from=1-2, to=2-2]
	\arrow["{\psi_{t,w,w'}}"{description}, Rightarrow, nfold, from=2-1, to=1-2]
	\arrow["{[\Hcal_A(\Gamma_{\Gbb})]}"', from=2-1, to=2-2]
\end{tikzcd}\]

It follows from \cite[Theorem 3.8.5]{ToyModel}  that this class is computed as 
$$\begin{aligned}\cl(\Hcal_A(\Gamma_{\Gbb}), \psi_{t,w,w'}) &= \Tr(\psi_{t,w,w'}, \Hcal_{A}(\Gamma_{\Gbb})) \\ &\in \DD(\HHH(\RR_A \otimes \RR_A, \psi_t \circ w\sigma \otimes \psi_t \circ w'\sigma)).\end{aligned}$$
As all the above isomorphisms are $A$-linear, and since by Lemma \ref{lem:finite-etale}, the algebra $\HHH(\RR_A \otimes \RR_A, \psi_t \circ w\sigma \otimes \psi_t \circ w'\sigma)[\frac{1}{|\Gbb|}]$ is finite étale over $A$, it is enough to prove that $\Tr(\psi_{t,w,w'}, \Hcal_{A}(\Gamma_{\Gbb}))[\frac{1}{|\Gbb|}]$ is compact as an $A$-module. Using Proposition \ref{prop:computing-in-the-trace}, we can recognize that 
\begin{align*}
&\HHH(\RR_A \otimes \RR_A, \psi_t \circ w\sigma \otimes \psi_t \circ w'\sigma) \\
&= \Hom_{\DD(\Rep_A^{\unip}(\Gbb(t)))}(\tr_A(\Delta_w *  \psi_t \circ \sigma (\Delta_w^{-1})), \tr_A(\psi_t \circ \sigma (\Delta_{w'}^{-1})  * \Delta_{w'})).
\end{align*}
Writing that $\tr_A = \ch_A \circ \tr_{A,\emptyset}$ and noting as before that both $\ch$ and $\hc^!$ preserve compact objects, we have 
\begin{align*}
&\Hom_{\DD(\Rep_A^{\unip}(\Gbb(t)))}(\tr_A(\Delta_w *  \psi_t \circ \sigma (\Delta_w^{-1})), \tr_A(\psi_t \circ \sigma (\Delta_{w'}^{-1})  * \Delta_{w'})) \\
&= \Hom_{\Ocal}(\hc^!\ch\tr_{A,\emptyset}(\Delta_w *  \psi_t \circ \sigma (\Delta_w^{-1})), \tr_{A,\emptyset}(\psi_t \circ \sigma (\Delta_{w'}^{-1})  * \Delta_{w'}))
\end{align*}
and therefore, it is enough to prove that for any two compact objects $M,N \in \Ocal_{A}^{\DL, \unip}$, the mapping complex $\Hom_{\Ocal}(M,N)[\frac{1}{|\Gbb|}]$ is compact in $\DD(A[\frac{1}{|\Gbb|}])$. Using the semi-orthogonal decomposition of $\Ocal$, see Remark \ref{rque:semi-orthogonal-decompositon-category-o}, this reduces to the fact that $\HHH(\RR_A, \psi_t \circ w\sigma)[\frac{1}{|\Gbb|}]$ is finite étale over $A$ by Lemma \ref{lem:finite-etale}.
\end{proof}

\begin{lem}\label{lem:generic-semisimplicity}
The category $\DD(\Rep^{\unip}_{\Q(t)}(\Gbb(t)))$ is semisimple.
\end{lem}

\begin{proof}
This is an application of Tits' deformation theorem \cite[Chap. IV, §2,
Exercise 27]{Bourbaki456}. Combining Corollary \ref{cor:generically-category-is-reduced} with Lemmas \ref{lem:base-change-for-reduced-hecke-cat} and \ref{lem:etale-specialization-for-representations}, we have that 
$$\DD(\Rep_{\Ql}^{\unip}(\Gbb(q))) = \DD(\Rep_{\Ql}^{\unip}(\Gbf^{\Frob}))$$
is semisimple. Therefore a specialization of the category $\DD(\Rep^{\unip}_{\Q(t)}(\Gbb(t)))$ to $\Ql$ is semisimple, hence the category itself is semisimple. 
\end{proof}

\subsection{Generic Harish--Chandra theory}

In this section, we work over $A = \Q[t,t^{-1}, |\Gbb|^{-1}]$. The goal is to work out the generic version of Harish--Chandra theory over this ring and in particular, to study the splitting field of $\Rep_{\Q(t)}(\Gbb(t))$ and prove the following theorem.

\begin{thm}\label{thm:splitting-field}
Assume that $\Gbb$ is semisimple and quasi-simple. 
The category $\Rep_{\Q(t)}(\Gbb(t))$ splits over $\overline{\Q}(t^{\frac{1}{2}})$. 
\end{thm}

\begin{rque}
The point of Theorem \ref{thm:splitting-field} is to control the part of the extension of $\Q(t)$ in $t$ which splits the category. It is likely that depending on the type of $\Gbb$, one can find a splitting field of the form $\Q(\zeta_n, t^{\frac{1}{2}})$ for some suitable $n$ and $\zeta_n$ a primitive $n$-th root of $1$. 
\end{rque} 

\begin{lem}
Let $I \subset S$. Both generic parabolic induction and restriction functors $R_I$ and ${^*}R_I$ are $t$-exact for the $t$-structure of Remark \ref{rque:t-structure}. 
\end{lem}

\begin{proof}
The fact that $R_I$ is $t$-exact follows from the fact that it sends the progenerator of $\Rep_A(\Lbb_I(t))$ to a summand of the progenerator of $\Rep_A^{\unip}(\Gbb(t))$. Since ${^*}R_I$ is both left and right adjoint to $R_I$ by Lemma \ref{lem:parab-reduction}, it is also $t$-exact. 
\end{proof}

\begin{defi}
We define the cuspidal and induced parts of $\Rep_A^{\unip}(\Gbb(t))$ as follows: 
\begin{enumerate}
\item the induced subcategory $\Rep_A^{\unip, \ind}(\Gbb(t))$ is the full subcategory generated by essential images of parabolic induction functors $R_I$ for $I \subsetneq S$ which are $\sigma$-stable, 
\item the cuspidal part $\Rep_A^{\unip, \cusp}(\Gbb(t))$ is 
$$\Rep_A^{\unip, \cusp}(\Gbb(t)) := \bigcap_{I \subsetneq S, \sigma(I) = I}\ker({^*}R_I).$$
\end{enumerate}
\end{defi}

\begin{lem}
There is a direct sum decomposition
$$\Rep_{\Q(t)}^{\unip}(\Gbb(t)) = \Rep_{\Q(t)}^{\unip, \ind}(\Gbb(t)) \bigoplus \Rep_{\Q(t)}^{\unip, \cusp}(\Gbb(t)).$$
\end{lem}

\begin{proof}
By Lemma \ref{lem:generic-semisimplicity}, the category $\Rep_{\Q(t)}^{\unip}(\Gbb(t))$ is semisimple, so an orthogonal decomposition is simply given by a partition of irreducible objects and the lemma is clear.
\end{proof}

\begin{lem}\label{lem:cuspidal-represenation-splitting-field}
Assume that $\Gbb$ is quasi-simple. The subcategory $\Rep_{\Q(t)}^{\unip, \cusp}(\Gbb(t))$ is split over $\overline{\Q}(t^{\frac{1}{2}})$.
\end{lem}

\begin{proof}
The specialization $\overline{\Q(t)} \rightsquigarrow \overline{\Q}$ obtained by sending $t \mapsto q$ yields a bijection 
$$\Irr_{\overline{\Q(t)}}^{\unip, \cusp}(\Gbb(t)) \cong \Irr^{\unip, \cusp}_{\overline{\Q}}(\Gbb(q)).$$
For any $w\in\Wbf$ and $\underline{w}\in\Reg(w)$, the object $\tr_A(\BS(\underline{w})) \in \Rep_A(\Gbb(t))$ is flat over $A$ and therefore the irreducible cuspidal components are the same over $\overline{\Q(t)}$ and over $\overline{\Q}$ for any specialization $t \mapsto q$. 

Consider the action of $\psi_t \circ \sigma$ on $\Rep_A(\Gbb(t))$ and let $\delta \geq 1$ be the order of $\sigma$. The object $\BS(\underline{w})$ is $(\psi_t \circ \sigma)^{\delta}$-equivariant and we have a well-defined isomorphism 
$$\restr_{(\psi_t \circ \sigma)^{\delta}}(\tr(\BS(\underline{w}))) \cong \tr(\BS(\underline{w})).$$
Denoting by $(-)^{\cusp}$ the projection onto the cuspidal part, we have a well-defined isomorphism 
$$\restr_{(\psi_t \circ \sigma)^{\delta}}(\tr(\BS(\underline{w}))^{\cusp}) \cong \tr(\BS(\underline{w}))^{\cusp}.$$
Given some absolutely irreducible generic  representation $\rho_{\gen} \in \Irr_{\overline{\Q(t)}}^{\unip, \cusp}(\Gbb(t))$ with corresponding representation $\rho_q \in \Irr^{\unip, \cusp}_{\overline{\Q}}(\Gbb(q))$, denote by $(-)_{[\rho_{\gen}]}$ (resp. $(-)_{[\rho_q]}$) the projection onto the isotypic component. We get an isomorphism 
$$\restr_{(\psi_t \circ \sigma)^{\delta}}(\tr(\BS(\underline{w}))_{[\rho_{\gen}]}) \cong \tr(\BS(\underline{w}))_{[\rho_{\gen}]}.$$
Now this is an isomorphism of $\Ocal_F$-modules where $F/\Q(t)$ is an extension over which $\rho_{\gen}$ and $\rho_q$ are defined, compatibly with the action of $(\psi_t \circ \sigma)^{\delta}$. 

We now claim that for any irreducible cuspidal unipotent generic representation $\rho_{\gen}$, there always exists some $\underline{w}$ such that $\rho_{\gen}$ appears with multiplicity one and $(\psi_t \circ \sigma)^{\delta}$ acts on the corresponding isotypic component with eigenvalues of the form $\zeta t^{N/2}$ for integers $N$. Both statements hold if they hold for all specializations. We can now proceed in a type-by-type analysis: 
\begin{enumerate}
\item it follows from Lusztig's classification \cite{LusztigBook} and the study of Deligne--Lusztig varieties attached to Coxeter elements \cite{LusztigCoxeterOrbit} that in types ${^2}A_n$, $B_n$, $C_n$, $D_n$, ${^2}D_n$, ${^3}D_4$, $E_6$, ${^2}E_6$, $E_7$, all cuspidal unipotent representations appear in $H^0(\tr(\Delta_c))$ hence in $\tr(\BS(\underline{c}))$ for any reduced expression of a Coxeter element $c$ with multiplicity one and distinct Frobenius eigenvalues. 
\item in types $F_4,G_2,E_8$, certain unipotent cuspidal representations may fail to appear in the cohomology of Deligne--Lusztig varieties attached to Coxeter elements. In general, those cuspidal representations appear in $H^{\ell(w)}_c(X(w), \Ql)$ with multiplicity one where $w$ is chosen to be minimal among those such that the cohomology of $X(w)$ contains the prescribed cuspidal representation, we refer to \cite{LusztigRationality}. 
\end{enumerate}
\end{proof}

\begin{lem}\label{lem:hecke-algebra-endomorphism}
Let $F/\Q(t)$ be a finite extension and let $I \subset S$ be a $\sigma$-stable subset. Let $\lambda \in \Irr_F^{\unip, \cusp}(\Lbb_I(t))$ be an absolutely irreducible cuspidal representation. Set $\Wbf(I, \lambda) = \Stab_{\Wbf^{\sigma}}(I,\lambda)$. We have 
\begin{enumerate}
\item $\Wbf(I, \lambda)$ is a Coxeter group, 
\item $\End(R_I(\lambda))$ is isomorphic to the Hecke algebra of $\Wbf(I, \lambda)$ with parameters given by $(t^{a_s})_{s \in S_{I,\lambda}}$ where $S_{I,\lambda}$ denotes the set of simple reflections of $\Wbf(I, \lambda)$ and the $a_s$ are positive integers.
\end{enumerate}
\end{lem}

\begin{proof}
Both points are well-known results after specialization; we refer to \cite[Theorem 3.2.5]{GeckMalleBook} for a complete exposition. In the generic situation, the first point can be tested after specialization and is thus settled. The second point follows from the Mackey formula, see Proposition \ref{prop:mackay-for-parab-induction} which holds over $\Q[t,t^{-1}, |\Gbb|^{-1}]$ and thus defines an integral basis of endomorphism algebras of parabolic inductions. The relations are therefore polynomials (or valued in some finite extension of $\Q[t,t^{-1}, |\Gbb|^{-1}]$) and can thus be tested after specializing to all $q$. 
\end{proof}

\begin{defi}
We define  the set of generic unipotent cuspidal data to be the set of pairs $(I, \lambda)$ where $I \subset S$ and $\lambda$ is an absolutely irreducible cuspidal representation of $\Lbb_I(t)$. There is a natural action of $\Wbf$ on the set of such data. 
\end{defi}

\begin{lem}
There is a natural partition 
$$\Irr_{\overline{\Q(t)}}(\Gbb(t)) = \bigsqcup_{(I,\lambda)/\Wbf}\Irr_{\overline{\Q(t)}}^{(I,\lambda)}(\Gbb(t))$$
where $\Irr_{\overline{\Q(t)}}^{(I,\lambda)}(\Gbb(t))$ denotes the set of irreducible representations in $R_I(\lambda)$. 
\end{lem}

\begin{proof}
This is standard and follows immediately from the Mackey formula of Proposition \ref{prop:mackay-for-parab-induction} and the semisimplicity of the category from Lemma \ref{lem:generic-semisimplicity}. 
\end{proof}

\begin{proof}[Proof of Theorem \ref{thm:splitting-field}]
It follows from Lemma \ref{lem:cuspidal-represenation-splitting-field} that all cuspidal data are defined over $\overline{\Q}(t^{\frac{1}{2}})$. Hence it is enough to control the splitting field of the corresponding Hecke algebras by Lemma \ref{lem:hecke-algebra-endomorphism}. By \cite[Theorem 8.2, Corollary 8.3]{CurtisRepresentationsHeckeAlgebras}, these Hecke algebras split over $\overline{\Q}(t^{\frac{1}{2}})$ and the theorem follows. 
\end{proof}

\subsection{Modular representations of $\Gbb(t)$}\label{ssec:modular-constancy-of-the-category}

The goal of this subsection is to analyze, for a fixed $\ell$, how the category $\Rep^{\unip}_{\Flb}(\Gbb(q))$ varies as $q$ varies. Recall the ring $\Ccal(\Zl^{\times}, \Fl)$ introduced in Section \ref{subsubsec:specializations}; we have a map 
$$\spe_{\Zl^{\times}} : B \to \Ccal(\Zl^{\times}, \Fl), \qquad t \mapsto \iota \  \mathrm{mod} \ \ell.$$
Note that since $\Fl$ is discrete, we have 
$$\Ccal(\Zl^{\times}, \Fl) = \varinjlim_n \Fun((\Z/\ell^n\Z)^{\times}, \Fl),$$
and the scheme $\Spec(\Ccal(\Zl^{\times}, \Fl))$ is $\underline{\Zl^{\times}} = \varprojlim_{n} \underline{(\Z/\ell^{n}\Z)^{\times}}$ where each $\underline{(\Z/\ell^{n}\Z)^{\times}}$ is the discrete scheme over $\Fl$ isomorphic to $|(\Z/\ell^{n}\Z)^{\times}|$ copies of $\Spec(\Fl)$. There is a natural open subscheme 
$$U = \underline{\Zl^{\times}} - \mu(\Zl) \subset \underline{\Zl^{\times}}$$
where we remove the finitely many closed points corresponding to roots of unity. Note that the scheme $U$ is not affine. 

Set $A = \Ccal(\Zl^{\times}, \Fl)$. We consider the category $\DD(\Rep_A^{\unip}(\Gbb(t)))$ and its associated algebra $E_A \in \DD(A)$.

\begin{lem}\label{lem:structure-modular-generic-representation}
\begin{enumerate}
\item There is a canonical equivalence 
\begin{equation}\label{eq:reduction-functor}
\DD(\Rep_{\Ccal(\Zl^{\times}, \Fl)}^{\unip}(\Gbb(t))) \otimes_{\DD({\Ccal(\Zl^{\times}, \Fl)}), \ev_q} \DD(\Fl) \cong \DD(\Rep_{\Fl}^{\unip}(\Gbb(q)))
\end{equation}
obtained by evaluation at $q \in \Zl^{\times}$. 
\item Over $U \subset \underline{\Zl^{\times}}$, the algebra $E_{\Ccal(\Zl^{\times}, \Fl)}$ is concentrated in degree $0$ and locally of finite presentation. 
\end{enumerate}
\end{lem}

\begin{proof}
In this proof, we shall write $A = \Ccal(\Zl^{\times}, \Fl)$. Given some clopen subset $V$ of $\Zl^{\times}$, we shall write $A_V = \Ccal(V, \Fl)$. Since $V$ is also closed, there is an idempotent $e_V \in A$ such that $e_VA = A_V$ and in particular, we have 
$$\Hcal_{A} \otimes_{\DD(A)} \DD(A_V) \cong \Hcal_{A_V}$$
which follows from the fact that $e_V\RR_A = \RR_{A_V}$. We now claim that for all $q \in \Zl^{\times}$, there exists a neighborhood over which the algebra $E_{A_V}$ is of finite presentation over $A_V$ and satisfies $E_{A_V} \otimes_{A_V} \Fl \cong E_q$.

Let us now prove the second point when $\Gbb$ is a torus. Let us first note that there are canonical isomorphisms
\begin{align*}
E_A &= \HHH(\RR_A, \psi_t \circ \sigma) \\
&\cong A \otimes_{A\llbracket X^{\vee} \otimes \Zl \rrbracket} A\llbracket X^{\vee} \otimes \Zl \rrbracket,
\end{align*}
where $A\llbracket X^{\vee} \otimes \Zl \rrbracket \to A$ is the augmentation and $A\llbracket X^{\vee} \otimes \Zl \rrbracket \to A\llbracket X^{\vee} \otimes \Zl \rrbracket$ is the map induced by $\id - \psi_t\circ\sigma : X^{\vee} \otimes \Zl \to X^{\vee} \otimes \Zl$ by Lemma \ref{lem:reduction-koszul-calculation}. Over $U$ the map $\id - \psi_t\circ\sigma$ is injective and has finite cokernel and we can locally on $U$ find bases $(e_i)$ and $(f_i)$ of $X^{\vee} \otimes \Zl$ such that the matrix of $\id - \psi_t\circ\sigma$ is equal to $\mathrm{diag}(t_1,\dots,t_n)$ where $t_i : U \to \Zl$ are continuous functions. It follows that 
$$ A \otimes_{A\llbracket X^{\vee} \otimes \Zl \rrbracket} A\llbracket X^{\vee} \otimes \Zl \rrbracket = \bigotimes_i A\llbracket f_i - 1 \rrbracket/((f_i - 1)^{\ell^{v_{\ell}(t_i)}})$$
which is locally of finite presentation over $U$ and in degree $0$. 

Let us now prove the claim that all cohomology groups are locally of finite presentation. Recall that by Proposition \ref{prop:computing-in-the-trace} and by the proof of Lemma \ref{lem:main-lemma-generator-of-trace} we can compute $E_{A}$ in terms of Hochschild homology of the algebra $\EE_{A}$ for the bimodule $\Mbb_{A}$. 
Let us introduce the relative bar resolution of $\EE_A$ over $\RR_A$ so that we have 
$$\EE_A = \varinjlim_{\Delta^{\opp}} \EE_A \otimes_{\RR_A} \dots \otimes_{\RR_A} \EE_A$$
where there are $(n+2)$ copies of $\EE_A$ in degree $n$. We then have 
\begin{align*}
\HHH(\EE_A, \Mbb_A) &= \varinjlim_{\Delta^{\opp}} \EE_A \otimes_{\RR_A} \dots \otimes_{\RR_A} \EE_A  \otimes_{\EE_A \otimes \EE_A} \Mbb_A  \\
&= \varinjlim_{\Delta^{\opp}} \EE_A \otimes_{\RR_A} \dots \otimes_{\RR_A} \EE_A  \otimes_{\RR_A \otimes \RR_A} \Mbb_A
\end{align*}
where there are $n$ copies of $\EE_A$ in the second line. As an $\RR_A$-bimodule, $\EE_A$ has a filtration whose graded pieces are of the form $\RR_{A,w}$ for various elements $w \in \Wbf$ and $\Mbb_A$ has a filtration as an $\RR_A$-bimodule whose graded pieces are $\RR_{A, \psi_t \circ v \circ \sigma}$ for various elements $v \in \Wbf$. In total, each of the terms in the simplicial object has a finite filtration with graded pieces of the form 
$$\RR_{A,w} \otimes_{\RR_{A} \otimes \RR_{A}} \RR_{A, \psi_t \circ v \circ \sigma} \cong \HHH(\RR_{A}, \psi_t \circ w' \circ \sigma)$$
for suitable elements $w,v,w' \in \Wbf$. By the torus case, all the terms of the simplicial object are in degree $0$ and locally of finite presentation and locally constant. And locally at a point $x \in U$, all the terms in the filtration of $\EE_A \otimes_{\RR_A} \dots \otimes_{\RR_A} \EE_A  \otimes_{\EE_A \otimes \EE_A} \Mbb_A$ are constant on a small enough neighborhood of $x$ and their formation commutes with base change, i.e.

$$(\EE_A \otimes_{\RR_A} \dots \otimes_{\RR_A} \EE_A  \otimes_{\EE_A \otimes \EE_A} \Mbb_A) \otimes_{A, x} \Fl \cong \EE_{\Fl} \otimes_{\RR_{\Fl}} \dots \otimes_{\RR_{\Fl}} \EE_{\Fl} \otimes_{\EE_{\Fl} \otimes \EE_{\Fl}} \Mbb_{\Fl}.$$
In particular, the map 
$$E_{A} \otimes_{A} \Fl \cong E_{\Fl},$$
is an isomorphism. By Remark \ref{rmk:tilting-sheaves-are-projective}, $E_{\Fl}$ is in degree $0$ and the lemma follows.

\end{proof}

\begin{thm}\label{thm:modular-local-constancy-of-category}
For every étale realization datum $(q, \ell)$, there exists $n_q \geq 1$ such that for every étale realization datum $(q', \ell)$ such that $v_{\ell}(q - q') \geq n_q$, there exists an equivalence of categories 
$$\Rep_{\Fl}^{\unip}(\Gbb(q)) \cong \Rep^{\unip}_{\Fl}(\Gbb(q'))$$
preserving Deligne--Lusztig cohomology complexes and tilting representations (in the sense of \cite{EteveTilting}).
\end{thm}

\begin{proof}
Let $(q,\ell)$ be an étale realization datum. Set $A = \Ccal(\Zl^{\times}, \Fl)$ and write $E_{\Fl,q,\ell} = E_A \otimes_{A,\ev_q} \Fl$. By Lemma \ref{lem:structure-modular-generic-representation}, the algebra $E_{\Ccal(\Zl^{\times}, \Fl)}$ is in degree $0$ and locally of finite presentation in a neighborhood of $\ev_q \in \Spec(\Ccal(\Zl^{\times}, \Fl))$. It then follows that there exists an $\ell$-adic neighborhood $\{q' \in \Zl^{\times} : v_\ell(q-q') \geq n_q\}$ of $q$ over which the algebra $E_{\Ccal(\Zl^{\times}, \Fl)}$ is constant. In particular, for any other étale realization datum $(q', \ell)$ such that $v_{\ell}(q-q') \geq n_q$, we have 
\begin{equation}\label{eq:iso-modular-algebras}
E_{\Fl,q,\ell} \cong E_{\Fl, q', \ell}.
\end{equation}
We therefore have a sequence of equivalences 
\begin{equation}\label{eq:equiv-modular-reps}
\Rep_{\Fl}^{\unip}(\Gbb(q)) \cong \Mod_{E_{\Fl,q,\ell}} \cong \Mod_{E_{\Fl,q',\ell}} \cong \Rep_{\Fl}^{\unip}(\Gbb(q')).
\end{equation}
Furthermore, the isomorphism \eqref{eq:iso-modular-algebras} preserves the idempotents corresponding to the decomposition $\mathbb{B}_A = \bigoplus_{w, \underline{w}} \BS(\underline{w})$, hence the equivalence \eqref{eq:equiv-modular-reps} sends $\tr(\BS(\underline{w}))$ to the corresponding object and thus preserves tilting representations and Deligne--Lusztig cohomology complexes. 
\end{proof}

\begin{rque}
It is likely that we can quantitatively refine the knowledge of integers $n_q$ to get more precise bounds. Note that, as a module, $E_{\Ccal(\Zl^{\times}, \Fl)}$ is locally constant over $U$ whenever each of the algebras $\HHH(\RR_{\Ccal(\Zl^{\times}, \Fl)}, \RR_{\Ccal(\Zl^{\times}, \Fl), \psi_t \circ w \circ \sigma})$ is constant. It is, however, not clear that this implies that $E_{\Ccal(\Zl^{\times}, \Fl)}$ is constant as an algebra. 
\end{rque}

\subsection{Around the decomposition matrix}

\subsubsection{Generalities on decomposition matrices}

In this subsection, we fix a noetherian base ring $k$ which we assume to be normal, integral and of finite type over $\Z$. Let $A$ be a finite flat $k$-algebra. 
\begin{defi}
\begin{enumerate}
\item When $k$ is a field, we say that $A$ is split if for all simple $A$-modules $M$, we have $\End_A(M) = k$,
\item When $k$ is a general commutative ring, we say that $A$ is split if for all prime ideals $\mathfrak{p} \in \Spec(k)$, the $\kappa(\mathfrak{p})$-algebra $A \otimes_k \kappa(\mathfrak{p})$ is split. 
\end{enumerate}
\end{defi}

\begin{defi}
\begin{enumerate}
\item A geometric point of $\Spec(k)$ is a map $s : \Spec(\Omega) \to \Spec(k)$ where $\Omega$ is an algebraically closed field.
\item A specialization map $s \to t$ of geometric points is a lift $s : \Spec(\Omega) \to \Spec(k)_{(t)}$ to the strict henselization of $\Spec(k)$ at $t$. 
\end{enumerate}
\end{defi}

\begin{defi}
Let $s \to t$ be a specialization map in $\Spec(k)$. There is a well-defined decomposition map
$$\dec_{s \to t} : G_0(A_s) \to G_0(A_t)$$
defined as usual, see \cite[Section 2]{GeckHeckeAlgRootOfUnity}. 
\end{defi}

\begin{lem}[\protect{\cite[Proposition 2.6]{GeckHeckeAlgRootOfUnity}}]\label{lem:composition-decomp-matrices} 
Let $s \to t \to u$ be specialization maps in $\Spec(k)$. The following triangle is commutative: 
\[\begin{tikzcd}
	& {G_0(A_s)} & \\
	{G_0(A_u)} && {G_0(A_t).}
	\arrow["{\dec_{s \to u}}"', from=1-2, to=2-1]
	\arrow["{\dec_{s \to t}}", from=1-2, to=2-3]
	\arrow["{\dec_{t \to u}}", from=2-3, to=2-1]
\end{tikzcd}\]
\end{lem}

\begin{thm}[\protect{\cite[Theorem 2.3]{ThielDecompMatrices}}]\label{thm:stratification-splitting-dec-map}
Assume that $A$ is split over $k$. Then there exists a finite stratification 
$$\Spec(k) = \bigsqcup_{i \in I} X_i$$
by locally closed subschemes $X_i$ such that for all $i \in I$ and every specialization map $s \to t$ in $X_i$, the map $\dec_{s \to t}$ is the identity.
\end{thm}

\begin{lem}\label{lem:generic-splitting}
Assume furthermore that $k$ is normal. Then there exists a finite Galois extension $F/\Frac(k)$ of the fraction field of $k$ such that $A_{k'}$ is split, where $k'$ is the normalization of $k$ in $F$.
\end{lem}

\begin{proof}
We prove by noetherian induction the following proposition: given a closed subset $Z \subset \Spec(k)$, there exists a finite extension $F'/F$, where $F = \Frac(k)$ such that over $k'$, for any point $z' \in \Spec(k')$ over a point of $Z$, the algebra $A_{\kappa(z')}$ is split. 

If $Z = \emptyset$, then this is clear. If $Z$ is reducible, then we can apply the induction hypothesis to all irreducible components and proceed. Hence we may assume that $Z$ is irreducible. Write $Z = V(\pfrak)$ for some prime ideal $\pfrak \subset k$. Let $E/\kappa(\pfrak)$ be the splitting field of $A_{\kappa(\pfrak)}$. We claim that there exists some finite extension $F'/F$ such that for all points $\pfrak'$ in $\Spec(k')$ over $\pfrak$ the algebra $A_{\kappa(\pfrak')}$ is split. Indeed, present $E = \kappa(\pfrak)(\alpha_1,\dots, \alpha_n)$ for elements $\alpha_i$ that are integral over $k/\pfrak$ and let $h_i \in (k/\pfrak)[x]$ be their minimal polynomials. Lift those polynomials to $H_i \in k[x]$ and take $F'$ to be the splitting field of all the $H_i$. All the roots of the $H_i$ are in $k'$ hence for any $\pfrak' \in \Spec(k')$ over $\pfrak$, in the quotient $k'/\pfrak'$, all the polynomials $(h_i)$ are split and therefore $\kappa(\pfrak')$ contains a copy of $E$ and the algebra $A_{\kappa(\pfrak')}$ is split. 

Since $k$ is of finite type over $\Z$, the map $\Spec(k') \to \Spec(k)$ is finite by \cite[Tag 035B]{StacksProject}. Consider the map 
$$g : \Spec(k') \times_{\Spec(k)} Z \to Z$$ 
this is a finite map. Over the generic point $\eta_Z$ of $Z$, the algebra $A_{g^{-1}(\eta_Z)}$ is split. As this fiber is finite, by \cite[Theorem 4.6]{ThielDecompMatrices}, for any point $\eta' \in g^{-1}(\eta_Z)$ there exists some open neighborhood of $\eta'$ in $\Spec(k') \times_{\Spec(k)} Z$ over which $A$ is split. Since there are finitely many of those and since the map $g$ is finite, there exists some open subset $U\subset Z$ such that over $g^{-1}(U)$ the algebra $A$ is split. Consider now the complement $Z' = Z-U \subsetneq Z$. We can apply the induction hypothesis to find an extension $F"/F'$ which splits $A$ over $k" \otimes_k \Ocal(Z')$. It then follows that $A$ is split over $k" \otimes_k \Ocal(Z)$.

\end{proof}

\begin{corol}\label{corol:finiteness-decomposition-matrices}
Assume that $k$ is normal. There are finitely many possible decomposition matrices for $A$ along any specialization map in $\Spec(k)$. 
\end{corol}

\begin{proof}
Upon replacing $k$ by its normalization in a finite extension of $\Frac(k)$, we may assume that $A$ is split over $k$ by Lemma \ref{lem:generic-splitting}. Then the finiteness is an immediate consequence of the finiteness of the stratification in Theorem \ref{thm:stratification-splitting-dec-map}.
\end{proof}

\subsubsection{Local properties of $E_{\red}$}

Consider the algebra $E_{\red} := E_{B,\red}^0$. For $N' \geq N$, we denote by $E_{N', \red}$ the algebra obtained in the same way as $E_{B, \red}^0$ after descending the reduced Hecke category to $B_{N'}$ via Lemma \ref{lem:existence-descent-datum}. As a module, it is a finite $B_{N'}$-module and we denote by $\overline{E}_{N', \red}$ its $B_{N'}$-torsion-free quotient. Similarly, we denote by $\overline{E}_{\red}$ the torsion-free quotient of $E_{\red}$. 

\begin{lem}\label{lem:structural-properties-reduced-ring}
\begin{enumerate}
\item There is a (noncanonical) isomorphism 
$$E_{\red} \cong \overline{E}_{\red} \oplus E_{\red}[\mathrm{tors}].$$
Moreover, the algebra $\overline{E}_{\red}$ is of finite presentation and flat over $B$.
\item Let $(q,\ell, \Zl)$ be an étale realization datum. There is a canonical isomorphism 
$$\overline{E}_{\red} \otimes_B \Zl \cong \overline{E}^{0}_{q,\ell,\red}.$$
\item There exists $N' \geq N$ such that 
$$\overline{E}_{N', \red} \otimes_{B_{N'}} B \cong \overline{E}_{\red}$$
and over $B_{N'}$, $\overline{E}_{N', \red}$ is a finite projective $B_{N'}$-module.
\end{enumerate}
\end{lem}

\begin{proof}
By \cite[Theorem VI.1.7]{CahenChabert}, the ring $B$ is a Prüfer domain and its local rings are valuation rings. As $E_{\red}^0$ is of finite presentation by Lemma \ref{lem:main-lemma-generator-of-trace}, there is a decomposition of the module into a torsion part and a torsion-free part and torsion-free modules are flat, see \cite[Tag 0539]{StacksProject}. This establishes the first item. 

Since $E_{N,\red}$ descends $E_{\red}$, the quotient $\overline{E}_{\red}$ descends to some finite step $E_{N', \red} \rightarrow \overline{E}_{N',\red}$ with 
$$\overline{E}_{N', \red} \otimes_{B_{N'}} B \cong \overline{E}_{\red},$$
this establishes the third point. 

For the second point, since we have $E_{\red} \cong \overline{E}_{\red} \oplus E_{\red}[\mathrm{tors}]$ applying the functor $\otimes_B \Zl$ gives a decomposition 
$$E^0_{q,\ell, \red} \cong \overline{E}_{\red} \otimes_{B} \Zl \oplus E_{\red}[\mathrm{\tors}] \otimes_{B} \Zl$$ 
where the first summand is $\Zl$-projective and the second one is $\Zl$-torsion hence agrees with the $\ell$-torsion of $E^0_{q,\ell, \red}$. It thus follows that $\overline{E}_{\red} \otimes_B \Zl \cong \overline{E}^{0}_{q,\ell,\red}$.
\end{proof}

From now on, fix $N' \geq N$ as in Lemma \ref{lem:structural-properties-reduced-ring}. Let us fix a pair $(q, \ell)$ defining an étale realization datum. Recall from Section \ref{subsubsec:specializations} that we denote by $b_{q,\ell}$ the point $B \xrightarrow{t \mapsto q} \Zl \to \Fl$. Let us view this as a point of $\Spec(B_{N\prime})$. 

\begin{defi}[Unipotent decomposition matrix] 
We denote by $D_{q,\ell}^{\unip}$ the unipotent part of the decomposition matrix of $\Gbb(q)$ mod $\ell$, that is, the part of the decomposition matrix of $\Gbb(q)$ whose rows and columns are indexed by irreducible unipotent representations over $\Qlb$ and $\Flb$.
\end{defi}

\begin{lem}\label{lem:comparison-decomposition-matrices}
Consider the specialization of points $b_{q, \ell, \Qlb} \rightsquigarrow b_{q, \ell,\Flb}$ in $\Spec(B_{N'})$. There are canonical isomorphisms 
\begin{align*}
G_0(\overline{E}_{N',\red} \otimes_{B_{N'}, \mathrm{sp}_{q,\ell}} \Qlb) &\cong G_0(\Rep^{\unip}_{\Qlb}(\Gbb(q))), \\
G_0(\overline{E}_{N',\red} \otimes_{B_{N'}, \mathrm{sp}_{q,\ell}} \Flb) &\cong G_0(\Rep^{\unip}_{\Flb}(\Gbb(q))),
\end{align*}
and under these isomorphisms, the specialization map $\mathrm{sp}_{b_{q, \ell, \Qlb} \to b_{q, \ell,\Flb}}$ is given by the matrix $D_{q,\ell}^{\unip}$. 
\end{lem}

\begin{proof}
This is immediate from the second point of Lemma \ref{lem:structural-properties-reduced-ring} and Lemma \ref{lem:control-of-ell-torsion}.
\end{proof}

\subsubsection{Structure of the decomposition matrix}

\begin{thm}\label{thm:finiteness-decomposition-matrices}
Given $\Gbb$, there exist finitely many possible matrices for $D_{q,\ell}^{\unip}$.
\end{thm}

\begin{proof}
By Lemma \ref{lem:comparison-decomposition-matrices}, all the matrices $D_{q,\ell}^{\unip}$ arise as matrices for certain specialization maps of the algebra $\overline{E}_{N',\red}$ but by Corollary \ref{corol:finiteness-decomposition-matrices} there are finitely many such matrices. 
\end{proof}

Recall from Theorem \ref{thm:splitting-field} that the category $\Rep_{\Q(t)}(\Gbb(t))$ and hence the algebra $E_{\Q(t),\red}$ are split over $\overline{\Q}(t^{\frac{1}{2}})$. 

\begin{thm}\label{thm:structure-decomposition-matrix}
For $\ell \geq N'$, where $N'$ is as in Lemma \ref{lem:structural-properties-reduced-ring}, let $(q,\ell)$ be an étale realization datum and let $\xi$ be the root of unity corresponding to $q$ mod $\ell$. Then there exist two matrices with integer entries $D_{\xi,\ell}$ and $D_{\xi}$ depending only on $(\xi, \ell)$ and $\xi$, respectively, such that 
$$D_{q, \ell}^{\unip} = D_{\xi,\ell}D_{\xi}.$$
\end{thm} 

\begin{proof}
The argument follows that of \cite[Example 3.2]{GeckHeckeAlgRootOfUnity}. Let $q^{1/2} \in \Qlb$ be a square root of $q$ and let $\zeta \in \Zlb$ be a primitive $2d$-th root of $1$ in $\Zlb$ and let $\overline{\zeta}$ denote its image in $\Flb$. We assume that $q = \overline{\zeta}$ in $\Flb$.

Over the ring $\Zlb[t^{1/2}, t^{-1/2}]$ the algebra $\overline{E}^0_{\red}$ has split generic fiber. Consider the diagram of specialization maps 
\[\begin{tikzcd}
	& {\Qlb(t^{1/2})} & \\
	{\Qlb = b_{q,\ell,\Qlb}} && \Qlb \\
	& {\Flb = b_{q,\ell, \Flb}}
	\arrow["{t^{1/2}\mapsto q^{1/2}}"{description}, dashed, from=1-2, to=2-1]
	\arrow["{t^{1/2} \mapsto \zeta}"{description}, dashed, from=1-2, to=2-3]
	\arrow[dashed, from=2-1, to=3-2]
	\arrow[dashed, from=2-3, to=3-2]
\end{tikzcd}\]
Our choice of root of $1$ implies that this diagram is commutative. By Lemma \ref{lem:composition-decomp-matrices}, there is an equality 
$$D^{\unip}_{q,\ell}D_{t^{1/2} \mapsto q^{1/2}} = D_{\xi,\ell}D_{\xi}$$
where the two factors on the right indicate the decomposition matrices along the two right vertical specializations in the above diagram. Since $\Rep^{\unip}_{\Qlb}(\Gbb(q))$ is semisimple by Lemma \ref{lem:generic-semisimplicity}, we have $D_{t^{1/2} \mapsto q^{1/2}} = I$, where $I$ is the identity matrix. 
\end{proof}

\appendix

\section{Recollections on categorical traces}\label{appendix:recollection-categorical-traces}

In this section, we review basic definitions and properties concerning categorical traces; we refer to \cite[Section 3]{ToyModel}, \cite[Section 7]{ZhuCategoricalTame} and \cite{HSS} for more details on the theory. 

\subsection{Traces and Hochschild homology}\label{sec:recollection-hh}

Let $\Obb$ be an ambient symmetric monoidal presentable category whose tensor product commutes with colimits in both variables. 

\begin{defi}[Dualizable object]
An object $x \in \Obb$ is dualizable if there exists an object $x^{\vee}$ and maps 
$$1_{\Obb} \xrightarrow{\coev} x \otimes x^{\vee}, \qquad x^{\vee} \otimes x \xrightarrow{\ev} 1_{\Obb}$$
such that the compositions 
\begin{align*}
x &\xrightarrow{\coev \otimes \id_x} x \otimes x^{\vee} \otimes x \xrightarrow{\id_x \otimes \ev} x \\
x^{\vee} &\xrightarrow{\id_{x^{\vee}} \otimes \coev} x^{\vee} \otimes x \otimes x^{\vee} \xrightarrow{\ev \otimes \id_{x^{\vee}}} x^{\vee},
\end{align*}
are isomorphic to the identity. 
\end{defi}

\begin{defi}[Trace of an endomorphism]
Let $x \in \Obb$ be a dualizable object and $f : x \to x$ be an endomorphism. We define its trace to be 
$$\Tr(f, x) \in \End(1_{\Obb})$$
given by the composition
$$1_{\Obb} \xrightarrow{\coev} x \otimes x^{\vee} \xrightarrow{f \otimes \id_{x^{\vee}}} x \otimes x^{\vee} \xrightarrow{\mathrm{sym}} x^{\vee} \otimes x \xrightarrow{\ev} 1_{\Obb}.$$
We shall denote by $\Omega(\Obb)$ the category of pairs $(x,f)$ where $x \in \Obb$ is a dualizable object and $f : x \to x$. 
\end{defi}

\begin{example}
\begin{enumerate}
\item If $\Obb = \Prr_{\Lambda}$, then for a dualizable category $\Ccal \in \Prr_{\Lambda}$ together with an endomorphism $\Frob : \Ccal \to \Ccal$, we have 
$$\Tr(\Frob, \Ccal) \in \DD(\Lambda).$$
\item If $\Obb = 2\Prr_{\Lambda}$ denotes the $3$-category of presentable $\Lambda$-linear $2$-categories in the sense of Scholze--Stefanich \cite{ScholzeStefanich} (for some appropriate cardinal $\kappa$ fixed once and for all), then for any pair $(\Ccal,\Frob) \in \Omega(2\Prr_{\Lambda})$, we have 
$$\Tr(\Frob, \Ccal) \in \Prr_{\Lambda}.$$
\item Any monoidal category $\Acal \in \Alg_{\E_1}(\Prr_{\Lambda})$ defines an object $\Prr_{\Acal} \in 2\Prr$. Moreover, any pair $(\Acal, \Frob) \in \Omega(\Alg_{\E_1}(\Prr_{\Lambda}))$ defines an object $(\Prr_{\Acal}, \restr_{\Frob}) \in \Omega(2\Prr)$.
\end{enumerate}
\end{example}

\begin{defi}[Bar complexes] 
\begin{enumerate}
\item Let $A \in \Alg(\Obb)$. We define the bar complex $B(A) \in \Fun(\Delta^{\opp}, \Obb)$ of $A$ as the simplicial object 
$$[n] \mapsto A^{\otimes (n+2)}$$
where the face maps and degeneracy maps are induced by partial multiplications and inserting units.
\item Let $A \in \Alg(\Obb)$ and $M \in \Obb$ be an $A$-bimodule. We define the cyclic bar complex $B^{\cyc}(A,M) \in \Fun(\Delta^{\opp}, \Obb)$ to be the simplicial object 
$$[n] \mapsto A^{\otimes n} \otimes M$$
where the face maps $(d_0, \dots, d_n) : B^{\cyc}_n(A,M) \to B^{\cyc}_{n-1}(A,M)$ are defined as 
\begin{enumerate}
\item if $0 \leq i \leq n-2$, then $d_i(a_0 \otimes \dots \otimes a_{n-1} \otimes m) := a_0 \otimes \dots \otimes a_ia_{i+1} \otimes \dots \otimes a_{n-1} \otimes m$, 
\item $d_{n-1}(a_0 \otimes \dots \otimes a_{n-1} \otimes m) = a_0 \otimes \dots \otimes a_{n-2} \otimes a_{n-1}m$,
\item $d_{n}(a_0 \otimes \dots \otimes a_{n-1} \otimes m) = a_1 \otimes \dots \otimes a_{n-1} \otimes ma_0$.
\end{enumerate}
\item Let $R \to A$ be a map of algebras. The bar complex (resp. cyclic bar complex) of $A$ (resp. $(A,M)$) relative to $R$ is the simplicial object given by 
\begin{enumerate}
\item $B(A/R) : [n] \mapsto A \otimes_R \dots \otimes_R A$ with $(n+2)$ copies of $A$, 
\item $B^{\cyc}((A,M)/R) : [n] \mapsto A \otimes_{R} \dots \otimes_R A \otimes_{R \otimes R} M$, where there are $n$ copies of $A$
\end{enumerate}
In both cases, the face maps and degeneracy maps are induced by the ones of $B(A)$ and $B^{\cyc}(A,M)$. 
\end{enumerate}
\end{defi}

\begin{defi}[Hochschild homology] 
Let $A \in \Alg_{\E_1}(\Obb)$ and $M$ be an $A$-bimodule. The Hochschild homology of $A$ with values in $M$ is 
$$\HHH(A,M) = A \otimes_{A \otimes A^{\rev}} M$$
where $A^{\rev}$ denotes $A$ with reversed multiplication, i.e. $x \otimes^{\rev}y = y \otimes x$. 
Given an endomorphism $\Frob : A \to A$, we denote by $A_{\Frob}$ the $A$-bimodule equal to $A$ as an object with right module structure twisted by $\Frob$. 
\end{defi}

The next lemma is standard; see for instance \cite[Section 3.6]{ToyModel}.

\begin{lem}\label{lem:trace-is-computed-via-HH}
Let $(A, \Frob) \in \Omega(\Alg_{\E_1}(\Prr))$, let $M$ be an $A$-bimodule, and let $R\to A$ be a map of algebras. There are canonical equivalences in $\Prr$ 
\begin{enumerate}
\item $\varinjlim_{\Delta^{\opp}} B(A) \cong \varinjlim_{\Delta^{\opp}} B(A/R) \cong A$,
\item $\varinjlim_{\Delta^{\opp}} B^{\cyc}(A,M) \cong \varinjlim_{\Delta^{\opp}} B^{\cyc}((A,M)/R) \cong \HHH(A,M)$,
\item $\Tr(\restr_{\Frob}, \Prr_A) \cong \HHH(A, A_{\Frob}).$
\end{enumerate}
\end{lem}

\subsection{Colimits and $\HHH$}

\begin{lem}[\protect{\cite[Proposition 1.8]{HorelRamzi}}]\label{lem:hh-is-symmetric-monoidal}
The functor
\[
\Omega(\Alg_{\E_1}(\Prr_{\Lambda})) \longrightarrow \Prr_{\Lambda},
\qquad (A,\psi) \longmapsto \HHH(A,\psi)
\]
commutes with tensor products.
\end{lem}

\begin{lem}\label{lem:commutation-trace-sifted-colimits}
The functor $\Omega(\Alg_{\E_1}(\Prr_{\Lambda})) \to \Prr_{\Lambda}, (A, \psi) \mapsto \HHH(A, \psi)$ commutes with sifted colimits. 
\end{lem}

\begin{proof}
This is a standard argument. The functor $(A, \psi) \mapsto \HHH(A, \psi)$ factors as 
\begin{align*}
\Omega(\Alg_{\E_1}(\Prr_{\Lambda})) &\rightarrow \Fun(\Delta^{\opp}, \Prr_{\Lambda}) \to \Prr_{\Lambda} \\
(A, \Frob) &\mapsto B^{\cyc}(A, A_{\Frob}),  \Fcal \mapsto \varinjlim_{\Delta^{\opp}} \Fcal.
\end{align*}
The second functor commutes with all colimits, so it is enough to prove that the first one commutes with sifted colimits. Since colimits in $\Fun(\Delta^{\opp}, \Prr_{\Lambda})$ are computed termwise, it is enough to prove that $A \mapsto A^{\otimes n}$ commutes with sifted colimits. This is now a direct calculation. Let $(A_i)_{i\in I}$ be a sifted diagram in $\Alg_{\E_1}(\Prr_{\Lambda})$ and set $A = \varinjlim_i A_i$,
\begin{align*}
A^{\otimes n} &= \varinjlim_{i_1, \dots, i_n} A_{i_1} \otimes A_{i_2} \otimes \dots \otimes A_{i_n} \\
&= \varinjlim_{i} A_{i}^{\otimes n},
\end{align*}
where the second isomorphism comes from the fact that since $I$ is sifted, $\varinjlim_{I^n}$ and $\varinjlim_{I}$ agree.
\end{proof}

\subsection{Functoriality of traces}\label{sec:funtoriality-of-traces}

Let $\Ccal, \Dcal$ be two rigid and compactly generated monoidal categories in $\Prr_{\Lambda}$. We assume that they come equipped with monoidal endofunctors $\Frob_{\Ccal}$ and $\Frob_{\Dcal}$. 

\begin{defi}\label{defi:functoriality-of-traces}
\begin{enumerate}
\item A $(\Dcal, \Ccal)$-bimodule $M$ defines a map 
$$[M] : \Prr_{\Ccal} \to \Prr_{\Dcal}, \qquad N \mapsto M \otimes_{\Ccal} N,$$
\item we say that $M$ is right adjointable if the functor $[M]$ admits a right adjoint; we shall denote by $M^R$ the bimodule defining the right adjoint. 
\item A morphism of pairs $(\Ccal, \Frob_{\Ccal}) \to (\Dcal, \Frob_{\Dcal})$ is a pair $(M, \phi)$ where $M$ is a right adjointable $(\Dcal, \Ccal)$-bimodule and $\phi$ is a map of bimodules $\phi : M \to M$ compatible with $\Frob_{\Ccal}$ and $\Frob_{\Dcal}$. 
\end{enumerate}
\end{defi}

Given a right adjointable bimodule $M$, there are well-defined unit and counit maps 
$$\unit : \Ccal \to M^R \otimes_{\Dcal} M$$
and 
$$\counit : M \otimes_{\Ccal} M^R \to \Dcal.$$

\begin{construction}
Given a morphism of pairs $(M,\phi) : (\Ccal, \Frob_{\Ccal}) \to (\Dcal, \Frob_{\Dcal})$, there is a well-defined map between traces 
\begin{align*}
\Tr(\restr_{\Frob_{\Ccal}}, \Prr_{\Ccal}) &\xrightarrow{\unit} \Tr(\restr_{\Frob_{\Ccal}} \circ [M^R \otimes_{\Dcal} M], \Prr_{\Ccal}) \\
&\xrightarrow{\cycl} \Tr([M] \circ \restr_{\Frob_{\Ccal}} \circ [M^R], \Prr_{\Dcal}) \\
&\xrightarrow{\phi} \Tr(\restr_{\Frob_{\Dcal}} \circ [M] \circ [M^R], \Prr_{\Dcal}) \\
&\xrightarrow{\counit} \Tr(\restr_{\Frob_{\Dcal}} , \Prr_{\Dcal}).
\end{align*}
\end{construction}

\begin{rque}\label{rque:class-map}
A key example of this functoriality is given by the class map (or categorical character), see \cite[Section 3.8.1]{ToyModel}. Let $M$ be a $\Ccal$-module with a map $\Frob_{M} : M \to M$ compatible with $\Frob_{\Ccal}$ and assume that $M$ is dualizable as a $\Ccal$-module (here since we assume that $\Ccal$ is rigid, this amounts to the dualizability of $M$ as a plain category). Then the pair $(M, \Frob_M)$ can be viewed as a morphism of pairs 
$$(\DD(\Lambda), \id) \to (\Ccal, \Frob_{\Ccal}).$$
We shall denote by $\cl(M, \Frob_{M}) \in \Tr(\restr_{\Frob_{\Ccal}}, \Prr_{\Ccal})$ the image of $\Lambda$ under the induced map $\Tr(\id,\Prr_{\Lambda}) = \DD(\Lambda) \to \Tr(\restr_{\Frob_{\Ccal}}, \Prr_{\Ccal})$. 
\end{rque} 

\begin{lem}[\protect{\cite[Section 3.6.9]{ToyModel}}]\label{lem:right-adjoint-to-induction}
Let $\iota : \Ccal \to \Dcal$ be a map commuting with $\Frob_{\Ccal}$ and $\Frob_{\Dcal}$. We can view $\Dcal$ as a $(\Dcal, \Ccal)$-bimodule. Then $\Dcal$ is right adjointable and $\Dcal^R \cong \Dcal$ viewed as a $(\Ccal, \Dcal)$-bimodule. 
\end{lem}

\begin{lem}[\protect{\cite[Section 3.6.5]{ToyModel}}]\label{lem:dual-of-categorical-trace}
Under our assumption that $\Ccal$ is compactly generated, so is $\Tr(\restr_{\Frob_{\Ccal}}, \Prr_{\Ccal})$. It follows that $\Tr(\restr_{\Frob_{\Ccal}}, \Prr_{\Ccal})$ is dualizable and we have 
$$\Tr(\restr_{\Frob_{\Ccal}}, \Prr_{\Ccal})^{\vee} \cong \Tr(\restr_{\Frob_{\Ccal}}, \Prr_{\Ccal^{\rev}}).$$
Moreover, the identity functor 
$$\id_{\Tr(\restr_{\Frob_{\Ccal}}, \Prr_{\Ccal})} \in \Tr(\restr_{\Frob_{\Ccal}}, \Prr_{\Ccal}) \otimes \Tr(\restr_{\Frob_{\Ccal}}, \Prr_{\Ccal})^{\vee}$$
is given by the class of $\Ccal$ viewed as a $\Ccal \otimes \Ccal^{\rev}$-module.
\end{lem}

\subsection{A Hochschild homology calculation}\label{ssec:HHcalculation}

Let $A$ be a commutative ring, set $\RR_A = A\llbracket x_1, \dots, x_n\rrbracket$, and denote the ideal $(x_1,\dots, x_n)$ by $I_A$. Let
$$\phi : \RR_A \to \RR_A$$
be an $A$-linear $I_A$-adically continuous automorphism. Let us denote by $\overline{\phi} : I_A/I_A^2 \to I_A/I_A^2$ the induced map. 

\begin{defi}[Continuous Hochschild homology]
We define continuous Hochschild homology as 
$$\HHH^{\cont}(\RR_A, \phi) := \RR_{A} \otimes_{\RR_A \widehat{\otimes}_A \RR_A} \RR_{A, \phi}$$
where $\RR_A \widehat{\otimes}_A \RR_A$ denotes the completed tensor product, i.e. $\RR_A \widehat{\otimes}_A \RR_A \cong A\llbracket x_1,\dots, x_n, y_1, \dots, y_n\rrbracket$. 
\end{defi}

\begin{lem}\label{lem:reduction-koszul-calculation}
\begin{enumerate}
\item Consider the two maps $\RR_A \to \RR_A/I_A$ and $\id - \phi : \RR_A \to \RR_A$ given by $x_i \mapsto 0$ and $x_i \mapsto x_i - \phi(x_i)$, respectively. We have an isomorphism in $\DD(A)$: 
$$\HHH^{\cont}(\RR_A, \phi) = A \otimes_{\RR_A, \id - \phi} \RR_A.$$
\item If the sequence $(x_1 - \phi(x_1), \dots, x_n - \phi(x_n))$ is regular in $\RR_A$ then $\HHH^{\cont}(\RR_A, \phi)$ is concentrated in degree $0$.
\item If $\id - \overline{\phi}$ is an automorphism of $I_A/I_A^2$, then $\HHH^{\cont}(\RR_A, \phi) = A[0]$.
\item Assume furthermore that $(\id - \phi)(I_A) \supset I_A^N$ for some $N \geq 1$, then the natural map 
$$\HHH(\RR_A, \phi) \to \HHH^{\cont}(\RR_A, \phi)$$
is an isomorphism. 
\end{enumerate}
\end{lem}

\begin{proof}
Let $\RR_A^{(1)} = A\llbracket t_1, \dots, t_n\rrbracket$ and $\RR_A^{(2)} = A\llbracket v_1, \dots, v_n\rrbracket$ be two copies of $\RR_A$. The sequence $(t_i - v_i) \in \RR_A^{(1)} \widehat{\otimes}_A \RR_A^{(2)}$ is regular. Let $\Kos^{\bullet}(t_i - v_i)$ be the associated Koszul complex. In particular, we have $\Kos^{\bullet}(t_i - v_i) \cong \RR_A$ in $\DD(\RR_A^{(1)} \widehat{\otimes}_A \RR_A^{(2)})$. By definition of Hochschild homology, we have 
\begin{align*}
\HHH^{\cont}(\RR_A, \RR_{A, \phi}) &= \RR_A \otimes_{\RR_A^{(1)} \widehat{\otimes}_A \RR_A^{(2)}} \RR_{A,\phi} \\ 
&\cong \Kos^{\bullet}(t_i - v_i) \otimes_{\RR_A^{(1)} \widehat{\otimes}_A \RR_A^{(2)}} \RR_{A, \phi} \\
&\cong \Kos_{\RR_A}(x_i - \phi(x_i)).
\end{align*}
Similarly, let $\Kos^{\bullet}_{\RR_A}(x_i)$ be the Koszul complex of the sequence $(x_i)$, so that $\Kos^{\bullet}_{\RR_A}(x_i) \cong A$ in $\DD(\RR_A)$. We have 
$$A \otimes_{\RR_A, \id - \phi} \RR_A \cong \Kos^{\bullet}_{\RR_A}(x_i) \otimes_{\RR_A,\id - \phi} \RR_A \cong \Kos_{\RR_A}(x_i - \phi(x_i)),$$
which proves the first part of the lemma.

Item $(2)$ is clear as $\HHH^{\cont}(\RR_A, \phi) \cong \Kos^{\bullet}_{\RR_A}(x_i - \phi(x_i))$. For item $(3)$, if $\id - \overline{\phi}$ is an isomorphism, then the sequence $(x_i - \phi(x_i))$ is regular and we have 
$$\HHH^{\cont}(\RR_A, \phi) = \RR_A/(\id - \phi)[0].$$
Since $\id - \overline{\phi}$ is an isomorphism, $(\id - \phi)(I_A) = I_A$ and we have $\RR_A/(\id - \phi) \cong A$. 

For item $(4)$, first note that $\HHH(\RR_A, \phi)$ is naturally an $\RR_A \otimes_A \RR_A$-module, and consider the ideal $I_A^{(2)} = I_A \otimes_A \RR_A + \RR_A \otimes_A I_A$. The $I_A^{(2)}$-adic completion of $\HHH(\RR_A, \phi)$ is naturally identified with $\HHH^{\cont}(\RR_A, \phi)$. Now the hypothesis of item $(4)$ guarantees that $\HHH(\RR_A, \phi)$ is killed by $I_A^{(2)}$ and in particular already derived $I_A^{(2)}$-complete, hence $\HHH^{\cont}(\RR_A, \phi) \cong \HHH(\RR_A, \phi)$. 
\end{proof}

\subsection{Computing in the trace}

Let $\Hcal$ be a $\Lambda$-linear rigid monoidal category with a monoidal endomorphism $\Frob : \Hcal \to \Hcal$. The next proposition is a formal consequence of \cite[Theorem 3.8.5]{ToyModel} when $x = 1_{\Hcal}$. The general case appears in \cite{ZhuCategoricalTame}.

\begin{prop}[\protect{\cite[Corollary 7.92, Example 7.93]{ZhuCategoricalTame}}]\label{prop:computing-in-the-trace}
Let $x,y \in \Hcal^{\omega}$ be compact objects. Consider the universal trace functor 
$$\tr_{\Frob} : \Hcal \to \Tr(\restr_{\Frob}, \Prr_{\Hcal}),$$
Then there is a canonical isomorphism in $\DD(\Lambda)$
$$\Hom_{\Tr(\restr_{\Frob}, \Prr_{\Hcal})}(\tr_{\Frob}(x), \tr_{\Frob}(y)) \cong \Tr(y \otimes \Frob(-) \otimes x^{\vee}, \Hcal),$$
where $y \otimes \Frob(-) \otimes x^{\vee} : \Hcal \to \Hcal$ is the endofunctor $z \mapsto y \otimes \Frob(z) \otimes x^{\vee}$. 
\end{prop}

\bibliographystyle{alpha}
\begingroup\raggedright
\bibliography{biblio}

@misc{EteveTilting,
 author = {Arnaud Eteve},
 title = {Tilting representations of finite groups of {Lie} type},
 year = {2025},
 howpublished = {Preprint, {arXiv}:2412.13326 [math.{RT}] (2025)},
 url = {https://arxiv.org/abs/2412.13326},
 arXiv = {arXiv:2412.13326}
}

@misc{EteveFreeMono,
 author = {Arnaud Eteve},
 title = {Free monodromic {Hecke} categories and their categorical traces},
 year = {2025},
 howpublished = {Preprint, {arXiv}:2412.13323 [math.{RT}] (2025)},
 url = {https://arxiv.org/abs/2412.13323},
 arXiv = {arXiv:2412.13323}
}

@misc{EteveJordan,
 author = {Arnaud Eteve},
 title = {Jordan decomposition, categorical traces and {Deligne}-{Lusztig} theory},
 year = {2026},
 howpublished = {Preprint, {arXiv}:2501.04113 [math.{RT}] (2026)},
 url = {https://arxiv.org/abs/2501.04113},
 arXiv = {arXiv:2501.04113}
}

@incollection{BezrukavnikovRiche,
 author = {Bezrukavnikov, Roman and Riche, Simon},
 title = {A topological approach to {Soergel} theory},
 booktitle = {Representation theory and algebraic geometry. A conference celebrating the birthdays of Sasha Beilinson and Victor Ginzburg, Chicago, IL, USA, August 21--25, 2017},
 isbn = {978-3-030-82006-0; 978-3-030-82007-7},
 pages = {267--343},
 year = {2022},
 publisher = {Cham: Birkh{\"a}user},
 language = {English},
 doi = {10.1007/978-3-030-82007-7_7},
 zbMATH = {7610537},
 Zbl = {1497.18019}
}

@book{BroueMalleMichel,
 author = {Brou{\'e}, Michel and Malle, Gunter and Michel, Jean},
 title = {Repr{\'e}sentations unipotentes g{\'e}n{\'e}riques et blocs des groupes r{\'e}ductifs finis avec un appendice de {George} {Lusztig}},
 fseries = {Ast{\'e}risque},
 series = {Ast{\'e}risque},
 issn = {0303-1179},
 volume = {212},
 year = {1993},
 publisher = {Montrouge: Soci{\'e}t{\'e} Math{\'e}\-ma\-tique de France},
 language = {French; German; English},
 zbMATH = {479873},
 Zbl = {0780.00005}
}

@article{BezYun,
 author = {Bezrukavnikov, Roman and Yun, Zhiwei},
 title = {On {Koszul} duality for {Kac}-{Moody} groups.},
 fjournal = {Representation Theory},
 journal = {Represent. Theory},
 issn = {1088-4165},
 volume = {17},
 pages = {1--98},
 year = {2013},
 language = {English},
 doi = {10.1090/S1088-4165-2013-00421-1},
 zbMATH = {6183352},
 Zbl = {1326.20051}
}

@article{EteveEberhardt,
 author = {Eberhardt, Jens Niklas and Eteve, Arnaud},
 title = {Universal {Koszul} duality for {Kac}-{Moody} groups},
 fjournal = {Compositio Mathematica},
 journal = {Compos. Math.},
 issn = {0010-437X},
 volume = {162},
 number = {5},
 pages = {1053--1092},
 year = {2026},
 language = {English},
 doi = {10.1017/S0010437X26102863},
 zbMATH = {8240201}
}

@article{TaylorUniversalMonodromic,
 author = {Taylor, Jeremy},
 title = {Universal monodromic tilting sheaves},
 fjournal = {Representation Theory},
 journal = {Represent. Theory},
 issn = {1088-4165},
 volume = {29},
 pages = {871--892},
 year = {2025},
 language = {English},
 doi = {10.1090/ert/705},
 zbMATH = {8129085}
}

@incollection{RouquierCategorification,
 author = {Rouquier, Rapha{\"e}l},
 title = {Categorification of {{\(\mathfrak{sl}_2\)}} and braid groups.},
 booktitle = {Trends in representation theory of algebras and related topics. Workshop on representations of algebras and related topics, Quer\'etaro, M\'exico, August 11--14, 2004.},
 isbn = {0-8218-3818-0},
 pages = {137--167},
 year = {2006},
 publisher = {Providence, RI: American Mathematical Society (AMS)},
 language = {English},
 zbMATH = {5083072},
 Zbl = {1162.20301}
}

@article{Demazure,
 author = {Demazure, Michel},
 title = {Une nouvelle formule des caract{\`e}res},
 fjournal = {Bulletin des Sciences Math{\'e}matiques. Deuxi{\`e}me S{\'e}rie},
 journal = {Bull. Sci. Math., II. S{\'e}r.},
 issn = {0007-4497},
 volume = {98},
 pages = {163--172},
 year = {1975},
 language = {French},
 zbMATH = {3568850},
 Zbl = {0365.17005}
}

@article{HorelRamzi,
 author = {Horel, Geoffroy and Ramzi, Maxime},
 title = {A multiplicative comparison of {Mac} {Lane} homology and topological {Hochschild} homology},
 fjournal = {Annals of \(K\)-Theory},
 journal = {Ann. \(K\)-Theory},
 issn = {2379-1683},
 volume = {6},
 number = {3},
 pages = {571--605},
 year = {2021},
 language = {English},
 doi = {10.2140/akt.2021.6.571},
 zbMATH = {7395129},
 Zbl = {1479.19003}
}

@misc{StacksProject,
    shorthand    = {Stacks},
    author       = {The {Stacks Project Authors}},
    title        = {\textit{Stacks Project}},
    howpublished = {\url{https://stacks.math.columbia.edu}},
    year         = {2018},
  }

@article{Bondarko,
 author = {Bondarko, M. V.},
 title = {Weight structures vs. {{\(t\)}}-structures; weight filtrations, spectral sequences, and complexes (for motives and in general)},
 fjournal = {Journal of \(K\)-Theory},
 journal = {J. \(K\)-Theory},
 issn = {1865-2433},
 volume = {6},
 number = {3},
 pages = {387--504},
 year = {2010},
 language = {English},
 doi = {10.1017/is010012005jkt083},
 zbMATH = {5862078},
 Zbl = {1303.18019}
}

@article{PittieSteinberg,
 author = {Steinberg, Robert},
 title = {On a theorem of {Pittie}},
 fjournal = {Topology},
 journal = {Topology},
 issn = {0040-9383},
 volume = {14},
 pages = {173--177},
 year = {1975},
 language = {English},
 doi = {10.1016/0040-9383(75)90025-7},
 zbMATH = {3498187},
 Zbl = {0318.22010}
}

@article{KostantKumar,
 author = {Kostant, Bertram and Kumar, Shrawan},
 title = {{{\(T\)}}-equivariant {{\(K\)}}-theory of generalized flag varieties},
 fjournal = {Journal of Differential Geometry},
 journal = {J. Differ. Geom.},
 issn = {0022-040X},
 volume = {32},
 number = {2},
 pages = {549--603},
 year = {1990},
 language = {English},
 doi = {10.4310/jdg/1214445320},
 zbMATH = {4209084},
 Zbl = {0731.55005}
}

@misc{ZhuCategoricalTame,
 author = {Xinwen Zhu},
 title = {Tame categorical local {Langlands} correspondence},
 year = {2025},
 howpublished = {Preprint, {arXiv}:2504.07482 [math.{RT}] (2025)},
 url = {https://arxiv.org/abs/2504.07482},
 arXiv = {arXiv:2504.07482}
}

@article{ToyModel,
 author = {Gaitsgory, D. and Kazhdan, D. and Rozenblyum, N. and Varshavsky, Y.},
 title = {A toy model for the {Drinfeld}-{Lafforgue} shtuka construction},
 fjournal = {Indagationes Mathematicae. New Series},
 journal = {Indag. Math., New Ser.},
 issn = {0019-3577},
 volume = {33},
 number = {1},
 pages = {39--189},
 year = {2022},
 language = {English},
 doi = {10.1016/j.indag.2021.11.002},
 zbMATH = {7478434},
 Zbl = {1535.14033}
}

@incollection{GeckHeckeAlgRootOfUnity,
 author = {Geck, Meinolf},
 title = {Representations of {Hecke} algebras at roots of unity},
 booktitle = {S\'eminaire Bourbaki. Volume 1997/98. Expos\'es 835--849},
 pages = {33--55, ex},
 year = {1998},
 publisher = {Paris: Soci{\'e}t{\'e} Ma\-th{\'e}\-ma\-tique de France},
 language = {English},
 url = {https://eudml.org/doc/110251},
 zbMATH = {1296375},
 Zbl = {0931.20006}
}

@article{ThielDecompMatrices,
 author = {Thiel, Ulrich},
 title = {Decomposition matrices are generically trivial},
 fjournal = {IMRN. International Mathematics Research Notices},
 journal = {Int. Math. Res. Not.},
 issn = {1073-7928},
 volume = {2016},
 number = {7},
 pages = {2157--2196},
 year = {2016},
 language = {English},
 doi = {10.1093/imrn/rnv204},
 zbMATH = {6992261},
 Zbl = {1404.16013}
}

@book{CahenChabert,
 author = {Cahen, Paul-Jean and Chabert, Jean-Luc},
 title = {Integer-valued polynomials},
 fseries = {Mathematical Surveys and Monographs},
 series = {Math. Surv. Monogr.},
 issn = {0076-5376},
 volume = {48},
 isbn = {0-8218-0388-3},
 year = {1997},
 publisher = {Providence, RI: American Mathematical Society},
 language = {English},
 doi = {10.1090/surv/048},
 zbMATH = {981272},
 Zbl = {0884.13010}
}

@misc{Bourbaki456,
 author = {Bourbaki, N.},
 title = {{\'E}l{\'e}ments de math{\'e}matique. {Fasc}. {XXXIV}. {Groupes} et alg{\`e}bres de {Lie}. {Chapitres} {IV}, {V} et {VI}: {Groupes} de {Coxeter} et syst{\`e}mes de {Tits}. {Groupes} engendr{\'e}s par des r{\'e}flexions. {Syst{\`e}mes} de racines},
 year = {1968},
 language = {French},
 howpublished = {Actualit{\'e}s {Scientifiques} et {Industrielles}. 1337. {Paris}: {Hermann} \& {Cie}. 288 p. {F} 48,00 (1968).},
 zbMATH = {3297198},
 Zbl = {0186.33001}
}

@article{HSS,
 author = {Hoyois, Marc and Scherotzke, Sarah and Sibilla, Nicol{\`o}},
 title = {Higher traces, noncommutative motives, and the categorified {Chern} character},
 fjournal = {Advances in Mathematics},
 journal = {Adv. Math.},
 issn = {0001-8708},
 volume = {309},
 pages = {97--154},
 year = {2017},
 language = {English},
 doi = {10.1016/j.aim.2017.01.008},
 url = {orbilu.uni.lu/handle/10993/51963},
 zbMATH = {6686567},
 Zbl = {1361.14014}
}

@misc{ScholzeStefanich, 
author = {Scholze, Peter}, 
year = {2026}, 
title = {Geometry and higher categories}, 
howpublished = {Available on the author's website}}

@book{HTT,
 author = {Lurie, Jacob},
 title = {Higher topos theory},
 fseries = {Annals of Mathematics Studies},
 series = {Ann. Math. Stud.},
 volume = {170},
 isbn = {978-0-691-14049-0; 978-0-691-14048-3},
 year = {2009},
 publisher = {Princeton, NJ: Princeton University Press},
 language = {English},
 doi = {10.1515/9781400830558},
 zbMATH = {5497319},
 Zbl = {1175.18001}
}

@misc{HA, 
author = {Lurie,Jacob}, 
title = {Higher Algebra},
year = {2017},
howpublished = {available on the author's website}}

@article{Borel,
 author = {Borel, Armand},
 title = {Sur la cohomologie des espaces fibr{\'e}s principaux et des espaces homog{\`e}nes de groupes de {Lie} compacts},
 fjournal = {Annals of Mathematics. Second Series},
 journal = {Ann. Math. (2)},
 issn = {0003-486X},
 volume = {57},
 pages = {115--207},
 year = {1953},
 language = {French},
 doi = {10.2307/1969728},
 zbMATH = {3084047},
 Zbl = {0052.40001}
}

@article{BroueMalle,
author = {Broué, Michel and  Malle, Gunter},
journal = {Mathematische Annalen},
number = {2},
pages = {241-262},
title = {Théorèmes de Sylow génériques pour les groupes réductifs sur les corps finis.},
url = {http://eudml.org/doc/164912},
volume = {292},
year = {1992},
}

@book{LusztigBook,
 author = {Lusztig, George},
 title = {Characters of reductive groups over a finite field},
 fseries = {Annals of Mathematics Studies},
 series = {Ann. Math. Stud.},
 volume = {107},
 year = {1984},
 publisher = {Princeton University Press, Princeton, NJ},
 language = {English},
 doi = {10.1515/9781400881772},
 zbMATH = {3886137},
 Zbl = {0556.20033}
}

@misc{GaitsgoryRozemblyumVarshavsky,
 author = {D. Gaitsgory and N. Rozenblyum and Y. Varshavsky},
 title = {Applications of (higher) categorical trace {II}: {Deligne}-{Lusztig} theory},
 year = {2026},
 howpublished = {Preprint, {arXiv}:2603.26957 [math.{AG}] (2026)},
 url = {https://arxiv.org/abs/2603.26957},
 arXiv = {arXiv:2603.26957}
}

@article{BZCHN,
 author = {Ben-Zvi, David and Chen, Harrison and Helm, David and Nadler, David},
 title = {Coherent {Springer} theory and the categorical {Deligne}-{Langlands} correspondence},
 fjournal = {Inventiones Mathematicae},
 journal = {Invent. Math.},
 issn = {0020-9910},
 volume = {235},
 number = {2},
 pages = {255--344},
 year = {2024},
 language = {English},
 doi = {10.1007/s00222-023-01224-2},
 zbMATH = {7790851},
 Zbl = {1558.14025}
}

@book{GaiRozDAG1,
 author = {Gaitsgory, Dennis and Rozenblyum, Nick},
 title = {A study in derived algebraic geometry. {Volume} {I}: {Correspondences} and duality},
 fseries = {Mathematical Surveys and Monographs},
 series = {Math. Surv. Monogr.},
 issn = {0076-5376},
 volume = {221},
 isbn = {978-1-4704-3569-1; 978-1-4704-3568-4; 978-1-4704-4085-5},
 year = {2017},
 publisher = {Providence, RI: American Mathematical Society (AMS)},
 language = {English},
 zbMATH = {6780490},
 Zbl = {1408.14001}
}

@misc{SandvikRelativeSerreDuality,
 author = {Colton Sandvik},
 title = {Relative {Serre} duality for {Coxeter} groups},
 year = {2026},
 howpublished = {Preprint, {arXiv}:2604.06084 [math.{RT}] (2026)},
 url = {https://arxiv.org/abs/2604.06084},
 arXiv = {arXiv:2604.06084}
}

@article{HoLiRelativeSerre,
 author = {Ho, Quoc P. and Li, Penghui},
 title = {Relative {Serre} duality for {Hecke} categories},
 fjournal = {Selecta Mathematica. New Series},
 journal = {Sel. Math., New Ser.},
 issn = {1022-1824},
 volume = {32},
 number = {4},
 pages = {11},
 note = {Id/No 68},
 year = {2026},
 language = {English},
 doi = {10.1007/s00029-026-01175-5},
 zbMATH = {8227907}
}

@article{LusztigParabCS,
 author = {Lusztig, George},
 title = {Parabolic character sheaves. {I}.},
 fjournal = {Moscow Mathematical Journal},
 journal = {Mosc. Math. J.},
 issn = {1609-3321},
 volume = {4},
 number = {1},
 pages = {153--179},
 year = {2004},
 language = {English},
 zbMATH = {2142714},
 Zbl = {1102.20030}
}

@incollection{CurtisRepresentationsHeckeAlgebras,
 author = {Curtis, Charles W.},
 title = {Representations of {Hecke} algebras},
 booktitle = {Orbites unipotentes et repr\'esentations. I. Groupes finis et alg\`ebres de Hecke},
 year = {1988},
 publisher = {Paris: Soci{\'e}t{\'e} Math{\'e}matique de France},
 language = {English},
 zbMATH = {4204651},
 Zbl = {0729.20016}
}

@article{LusztigCoxeterOrbit,
 author = {Lusztig, G.},
 title = {Coxeter orbits and eigenspaces of {Frobenius}},
 fjournal = {Inventiones Mathematicae},
 journal = {Invent. Math.},
 issn = {0020-9910},
 volume = {38},
 pages = {101--159},
 year = {1976},
 language = {English},
 doi = {10.1007/BF01408569},
 url = {https://eudml.org/doc/142444},
 zbMATH = {3570675},
 Zbl = {0366.20031}
}

@article{LusztigRationality,
 author = {Lusztig, G.},
 title = {Rationality properties of unipotent representations.},
 fjournal = {Journal of Algebra},
 journal = {J. Algebra},
 issn = {0021-8693},
 volume = {258},
 number = {1},
 pages = {1--22},
 year = {2002},
 language = {English},
 doi = {10.1016/S0021-8693(02)00514-8},
 zbMATH = {1868077},
 Zbl = {1141.20300}
}

@book{GeckMalleBook,
 author = {Geck, Meinolf and Malle, Gunter},
 title = {The character theory of finite groups of {Lie} type. {A} guided tour},
 fseries = {Cambridge Studies in Advanced Mathematics},
 series = {Camb. Stud. Adv. Math.},
 volume = {187},
 isbn = {978-1-108-48962-1; 978-1-108-77908-1},
 year = {2020},
 publisher = {Cambridge: Cambridge University Press},
 language = {English},
 doi = {10.1017/9781108779081},
 zbMATH = {7162758},
 Zbl = {1536.20002}
}

@misc{ReevesUnipotentTrace,
 author = {Reeves, Wyatt},
 title = {Unipotent representations and the categorical trace of {Frobenius}},
 year = {2026},
 howpublished = {Preprint, arXiv:2609.39721 [math.RT]}
}

@book{FriedlanderEtaleHomotopy,
 author = {Friedlander, Eric M.},
 title = {{\'E}tale homotopy of simplicial schemes},
 series = {Annals of Mathematics Studies},
 volume = {104},
 publisher = {Princeton University Press, Princeton, N.J.; University of Tokyo Press, Tokyo},
 year = {1982}
}

@article{BrotoMollerOliver,
 author = {Broto, Carles and M{\o}ller, Jesper and Oliver, Bob},
 title = {Equivalences between fusion systems of finite groups of {Lie} type},
 journal = {J. Amer. Math. Soc.},
 volume = {25},
 number = {1},
 pages = {1--20},
 year = {2012},
 doi = {10.1090/S0894-0347-2011-00713-3}
}

@article{RouquierModularLie,
 author = {Rouquier, Rapha{\"e}l},
 title = {Modular representations of finite groups and {Lie} theory},
 journal = {J. Algebra},
 volume = {656},
 pages = {446--485},
 year = {2024},
 doi = {10.1016/j.jalgebra.2023.10.009}
}

@article{DipperJamesQSchur,
 author = {Dipper, Richard and James, Gordon},
 title = {The {$q$}-{Schur} algebra},
 journal = {Proc. London Math. Soc. (3)},
 volume = {59},
 number = {1},
 pages = {23--50},
 year = {1989},
 doi = {10.1112/plms/s3-59.1.23}
}

@article{LusztigCellsAndTensorCat,
 author = {Lusztig, G.},
 title = {Cells in affine {Weyl} groups and tensor categories},
 fjournal = {Advances in Mathematics},
 journal = {Adv. Math.},
 issn = {0001-8708},
 volume = {129},
 number = {1},
 pages = {85--98},
 year = {1997},
 language = {English},
 doi = {10.1006/aima.1997.1645},
 zbMATH = {1042603},
 Zbl = {0884.20026}
}

@article{LusztigCategoricalCentre,
 author = {Lusztig, G.},
 title = {Unipotent representations as a categorical centre},
 fjournal = {Representation Theory},
 journal = {Represent. Theory},
 issn = {1088-4165},
 volume = {19},
 pages = {211--235},
 year = {2015},
 language = {English},
 doi = {10.1090/ert/468},
 zbMATH = {6503624},
 Zbl = {1360.20041}
}

@article{EliasWilliamson,
 author = {Elias, Ben and Williamson, Geordie},
 title = {Relative hard {Lefschetz} for {Soergel} bimodules},
 fjournal = {Journal of the European Mathematical Society (JEMS)},
 journal = {J. Eur. Math. Soc. (JEMS)},
 issn = {1435-9855},
 volume = {23},
 number = {8},
 pages = {2549--2581},
 year = {2021},
 language = {English},
 doi = {10.4171/JEMS/1061},
 zbMATH = {7367669},
 Zbl = {1475.20006}
}
\endgroup

\end{document}